\documentclass[draft]{amsart}
\usepackage{amsmath,amsthm,amssymb}

\theoremstyle{plain}
\newtheorem{thm}{Theorem}[subsection]
\newtheorem{prop}[thm]{Proposition}

\newtheorem{cor}[thm]{Corollary}

\theoremstyle{definition}

\newtheorem{rem}{Remark}
\newtheorem{defn}[thm]{Definition}
\newtheorem{eg}[thm]{Example}
\newtheorem{subtitle}[thm]{}
\numberwithin{equation}{subsection}
\begin{document}

\def\a{\alpha}
\def\b{\beta}
\def\d{\delta}
\def\D{\triangle}
\def\e{\epsilon}
\def\g{\gamma}
\def\G{\Gamma}
\def\K{\nabla}
\def\l{\lambda}
\def\L{\Lambda}
\def\n{\,\vert\,}
\def\o{\theta}
\def\w{\omega}
\def\W{\Omega}
\def\ca{{\mathcal{A}}}
\def\cd{{\mathcal{D}}}
\def\cf{{\mathcal{F}}}
\def\cg{{\mathcal{G}}}
\def\ck{{\mathcal{K}}}
\def\cl{{\mathcal{L}}}
\def\cL{{\mathcal{L}}}
\def\cm{{\mathcal{M}}}
\def\cn{{\mathcal{N}}}
\def\co{{\mathcal{O}}}
\def\cp{{\mathcal{P}}}
\def\cs{{\mathcal{S}}}
\def\cu{{\mathcal{U}}}
\def\cv{{\mathcal{Vs}}}
\def\li{\langle}
\def\ri{\rangle}
\def\n{\ \vert\ }
\def\tr{{\rm tr}}
\def\bs{\bigskip}
\def\ms{\medskip}
\def\ss{\smallskip}
\def\hs{\hskip .05in}
\def\di{$\diamond$}
\def\ni{\noindent}
\def\ti{\tilde}
\def\p{\partial}
\def\Re{{\rm Re\/}}
\def\Im{{\rm Im\/}}
\def\I{{\rm I\/}}
\def\II{{\rm II\/}}
\def\diag{{\rm diag}}
\def\ad{{\rm ad}}
\def\Ad{{\rm Ad}}
\def\Iso{{\rm Iso}}

\def\h{|\det(h_{ij})|}
\def\III{{\rm III\/}}
\def\hol{{\rm Hol\/}}
\def\alg{{\rm alg\/}}

\def\R{\mathbb{R} }
\def\C{\mathbb{C}}
\def\H{\mathbb{H}}
\def\N{\mathbb{N}}
\def\Z{\mathbb{Z}}
\def\O{\\mathbb{O}}
\def\fg{\mathfrak{G}}

\newcommand{\beq}{\begin{equation}}
\newcommand{\eeq}{\end{equation}}
\newcommand{\beg}{\begin{eg}}
\newcommand{\eeg}{\end{eg}}
\newcommand{\bthm}{\begin{thm}}
\newcommand{\ethm}{\end{thm}}
\newcommand{\bprop}{\begin{prop}}
\newcommand{\eprop}{\end{prop}}
\newcommand{\bcor}{\begin{cor}}
\newcommand{\ecor}{\end{cor}}
\newcommand{\bca}{\begin{cases}}
\newcommand{\eca}{\end{cases}}
\newcommand{\brem}{\begin{rem}}
\newcommand{\erem}{\end{rem}}
\newcommand{\bpm}{\begin{pmatrix}}
\newcommand{\epm}{\end{pmatrix}}
\newcommand{\bdefn}{\begin{defn}}
\newcommand{\edefn}{\end{defn}}
\newcommand{\bsub}{\begin{subtitle}}
\newcommand{\esub}{\end{subtitle}}

\def\onn{\frac{O(2n)}{O(n)\times O(n)}}
\def\un{\frac{U(n)}{O(n)}}
\def\sun{\frac{SU(n)}{SO(n)}}

\title[Geometry of integrable systems]
{Geometries and Symmetries of Soliton equations  and\\
 Integrable Elliptic equations}

\author{ Chuu-Lian Terng$^1$}\thanks{$^1$Research supported
in  part by
NSF Grant DMS 9972172 and National Center for Theoretic Sciences, Taiwan.}
\address{Department of Mathematics\\
Northeastern University\\Boston, MA 02115}
\email{terng@neu.edu}

\begin{abstract}
We give a review of the systematic construction of hierarchies of soliton flows and integrable elliptic
equations associated to  a complex semi-simple Lie algebra and finite order automorphisms.  
For example, the non-linear
Schr\"odinger equation, the n-wave equation, and the sigma-model are soliton flows; and the equation 
for harmonic maps from the plane to a compact Lie group, for primitive maps from the plane to a
$k$-symmetric space, and constant mean curvature surfaces and isothermic surfaces in space forms are
integrable elliptic systems.  We also give a survey of  

\begin{itemize}
\item construction of solutions using loop group
factorizations, 
\item  PDEs in differential geometry that are soliton equations or elliptic
integrable systems, 
\item  similarities and differences of soliton equations and integrable elliptic systems.
\end{itemize}
\end{abstract}


\maketitle

\bigskip

\tableofcontents

\bs
\section{\bf Introduction} \label{ka}
\ms
 
We review the geometries and symmetries of both integrable evolution equations and 
elliptic partial differential equations, and also the method of loop group factorization for constructing 
their solutions.  In the classical literature, a differential equation is called ``integrable'' if it can 
be solved by quadratures.  A Hamiltonian system in $2n$-dimensions is {\it completely integrable\/} 
if it has $n$ independent commuting Hamiltonians.  By the Arnold-Liouville Theorem, such 
systems have action-angle variables that linearize the flow, and these can be found by quadrature.  
This concept of integrability can be extended to PDEs, and one class 
consists of evolution equations on function spaces that have Hamiltonian structures 
and are completely integrable Hamiltonian systems in the sense of Liouville, 
i.e., there exist action angle variables.  We call this class of equations soliton
equations.  The model examples are the Korteweg-de Vries equation, the non-linear Schr\"odinger
equation (NLS), and the Sine-Gordon equation (SGE).  Besides the Hamiltonian formulation and complete
integrability, soliton equations share many other remarkable properties including:
\begin{enumerate}
\item[\di] infinite families of explicit solutions,
\item[\di] a hierarchy of  commuting flows described by  partial
differential equations,
\item[\di] a Lax pair,
\item[\di] an algebraic-geometric description of certain solutions,
\item[\di] a scattering theory,
\item[\di] an inverse scattering transform to solve the Cauchy problem,
\item[\di] a construction of solutions using loop group factorizations (dressing actions). 
\end{enumerate}

\ms
 The existence of a Lax pair is one of the key properties of soliton equations.  
This was first constructed for the case of KdV by Lax, who observed that the KdV 
equation can be written as the condition for an isospectral deformation of the
Schr\"odinger operator on the line.  Later, this was shown to be equivalent to 
the zero curvature condition of a family of connections (\cite{AKNS74, ZakSha79}). 
 Roughly speaking, a PDE for $q:\R^n\to \R^m$ is said to have a  {\it zero curvature formulation\/} if there is
a family of connections $\o_\l$ on $\R^n$, (defined by $q$ and its derivatives, and  a holomorphic parameter
$\l$ defined in some open subset of
$\C$)  so that the condition for $\o_\l$ to be flat for all $\l$ is that
$q$ solve the given PDE.  The connection $\o_\l=\sum_{i=1}^n A_i dx_i$ is flat if  $d\o_\l=-\o_\l\wedge \o_\l$ for
all $\l$, or equivalently the $n$ operators $\{\frac{\p}{\p x_i} + A_i\n 1\leq i\leq n\}$ commute, i.e.,
$$\left[\frac{\p}{\p x_i}+A_i, \  \frac{\p}{\p x_j}+ A_j\right]=0,  \quad i\not=j.$$
We call $\o_\l$ a {\it Lax pair\/} if $n=2$, and a {\it Lax $n$-tuple\/} for general $n$.  A Lax $n$-tuple  
naturally gives rise to a  loop group factorization, which in turn provides a method for constructing explicit
solutions and symmetries of the equations. 

Another class of integrable PDEs are non-linear elliptic equations.  Although these equations do not
have Hamiltonian formulations, they do have zero curvature formulations that give rise to loop group
factorizations, and  hence the techniques developed for soliton equations can also be used to construct solutions
and  symmetries of these elliptic equations.  In particular, we can find solutions of the equation by factorizations.  
One class of model examples are the equations for harmonic maps from $\C$ to a compact Lie
group.  

Some goals of this paper are to give a brief survey of the following:

\begin{enumerate}
\item[\di] A systematic construction of integrable hierarchies associated to 
a complex semi-simple Lie algebra and finite order automorphisms.
\item[\di] Some geometric integrable PDEs arising in differential geometry.
\item[\di] Construction of solutions using loop group factorizations.  
\end{enumerate}

\ss\ni 
Another goal of this paper is to put some known results of evolution soliton equations and
integrable elliptic systems together so that we can compare and see similarities and
differences in these two theories.  

\ms
\ni {\bf $G$-hierarchy\/}
\ss
The ZS-AKNS construction of the $n\times n$-hierarchy of soliton flows works equally well when we
replace $sl(n,\C)$ by any complex, simple Lie algebra $\cg$.  In fact, let $a\in \cg$,
$\cg_a= \{y\in \cg\n [a,y]=0\}$ the centralizer of $a$,  and $\cg_a^\perp =\{\xi\in \cg\n (\xi,y)=0$ for all $y\in
\cg_a\}$.  Here $(\ , )$ is a non-degenerate ad-invariant bilinear form of $\cg$.  It can be shown that there exists a
sequence of polynomial differential operators  on the space 
$C(\R,\cg_a^\perp)$ of smooth functions from
$\R$ to
$\cg_a^\perp$,
$$\{Q_{b,j}(u)\n  b\in \cg_a, \ \cg_b=\cg_a, \ j\geq 0 \ \text{integer}\}.$$  These $Q_{b,j}(u)$ are
determined uniquely from the following recursive formula
\begin{align*}
&(Q_{b,j}(u))_x + [u, Q_{b,j}(u)]=[Q_{b,j+1}(u), a], \cr
& Q_{b,0}=b,\ \  Q_{a,1}(u)=u.\cr
\end{align*}
The $(b,j)$-flow is 
$$u_t=(Q_{b,j}(u))_x + [u,Q_{b,j}(u)],$$
which commutes with the $(b',j')$-flow.  
The hierarchy of these commuting flows is called $\mathfrak{g}$AKNS-hierarchy in \cite{Wil91}, and the
$G$-hierarchy in \cite{TerUhl98}.  

It follows from the recursive formula that $u$ is a solution of the $(b,j)$-th flow if and only if 
$$\o_\l= (a\l + u) dx + (b\l^j + Q_{b,1}(u)\l^{j-1} + \cdots + Q_{b,j}(u)) dt$$
is flat for all $\l\in\C$.  In other words, $\o_\l$ is a Lax pair of the $(b,j)$-flow.

Next we explain several invariant submanifolds of the
$G$-hierarchy in terms of finite order automorphisms of $\cg$.    

\ss
\ni {\bf $\sigma$-twisted $G$-hierarchy\/}
\ss

  If $\sigma$ is an order $k$ automorphism of the complex Lie group $G$, then the $(b,nk+1)$-flow in
the $G$-hierarchy leaves  $C(\R,\cg_a^\perp\cap\ \cg_0)$ invariant, 
where $\cg_0$ is the fixed point set of $d\sigma_e$ on $\cg$.  The hierarchy of the
restriction of these flows to $C(\R, \cg_a^\perp\cap \cg_0)$ is called the {\it $\sigma$-twisted
$G$-hierarchy\/}. The Kupershmidt-Wilson hierarchy is an example with $\cg=sl(n,\C)$ and $k=n$
(\cite{KupWil81}).

\ss
\ni {\bf $U$-hierarchy} \ss

 Suppose $\tau$ is a conjugate linear, Lie algebra involution of $\cg$, and $\cu$ is the fixed point set of
$\tau$, i.e., $\cu$ is a real form of $\cg$. Then the $(b,j)$-flow leaves $C(\R, \cg_a^\perp\cap\ \cu)$
invariant.  The hierarchy restricted to $C(\R, \cg_a^\perp\cap\ \cu)$ is called the $U$-hierarchy.  
For example, the NLS occurs as the second flow in the
$SU(2)$-hierarchy, and the $3$-wave equation as the first flow in the $SU(3)$-hierarchy.

\ss
\ni {\bf $U/U_0$-hierarchy\/}\ss

 Suppose $\tau$ is a conjugate-linear involution, and $\sigma$ is an order $k$
complex linear, Lie algebra automorphism of $\cg$ such that 
$$\sigma\tau=\tau^{-1}\sigma^{-1}.$$ 
Then the $(b,nk+1)$-flow leaves
$C(\R,\cg_a^\perp\cap \cu_0)$ invariant, where $\cu_0$ is the Lie subalgebra of $\cg$ that is fixed by
both $\sigma$ and $\tau$.  The hierarchy restricted to $C(\R, \cu_a^\perp\cap \cu_0)$ is called the
$U/U_0$-hierarchy.  For example, the
$3$rd flow in the $SU(2)/SO(2)$-hierarchy is the modified KdV equation with $k=2$. 

\ss
\ni {\bf $U/U_0$-system\/}\ss

Let $U/U_0$ be the rank $n$ symmetric space given by involutions $\tau, \sigma$ of $G$,
$\cu=\cu_0+\cu_1$ the Cartan decompostion, and $\ca$ a maximal abelian subspace of $\cu_1$.  Let
$\{a_1, \cdots, a_n\}$ be a basis of $\ca$.  By putting the $(a_1, 1), \cdots, (a_n,1)$-flows in the
$U/U_0$-hierarchy together, we get the {\it
$U/U_0$-system\/} for maps $v:\R^n\to \cu_{\ca}^\perp\cap \cu_1$:
\begin{equation*}
[a_i, v_{x_j}]-[a_j, v_{x_i}]= [[a_i, v], [a_j,v]], \quad i\not=j,
\end{equation*}
where $\cu_\ca=\{y\in \cu\n [y, \xi]=0 \ \forall\  \xi\in \ca\}$.
Note that $v$ is a solution of the $U/U_0$-system if and only if 
$$\o_\l = \sum_{j=1}^n \ (a_i\l+ [a_i, v])\ dx_i$$
is flat for all $\l\in\C$, i.e., $\o_\l$ is a Lax $n$-tuple of the $U/U_0$-system.  

\ss
\ni {\bf The $-1$-flow associated to $U$\/}\ss

Let $a, b\in \cu$ such that $[a,b]=0$.  The $-1$-flow associated to $U$ is the following system for $g:\R^2\to U$:
$$(g^{-1}g_x)_t=[a, g^{-1}bg],$$
with constraint $g^{-1}g_x\in \cu_a^\perp$.  The $-1$-flow has a Lax pair 
$$(a\l + g^{-1}g_x) dx + \l^{-1}
g^{-1}bg \ dt.$$

\ss
\ni {\bf Elliptic $(G,\tau)$-systems} \ss

 The $m$-th elliptic $(G,\tau)$-system is the equation for $(u_0,\cdots,
u_m):\C\to  \oplus_{i=0}^m\cg$ so that 
$$\o_\l= \sum_{j=0}^m \l^{-j} u_j dz + \l^j \tau(u_j) d\bar z$$ is flat for all $\l\in S^1$.   The first
$(G,\tau)$-system is the equation for harmonic maps from $\R^2$ to $U$, where $U$ is the fixed point set of
$\tau$.    

\ss\ni {\bf Elliptic $(G,\tau,\sigma)$-systems}\ss

Suppose $\sigma$ is an order $k$ automorphism of $G$ such that 
$$\sigma\tau=\tau\sigma.$$  Let $\cg_j$ denote the eigenspace of $\sigma_\ast$ on $\cg$ with
eigenvalue $e^\frac{2\pi i j}{k}$.
We call the $m$-th elliptic $(G,\tau)$-system with constraints $u_i\in
\cg_{-i}$ the $m$-th elliptic $(G,\tau,\sigma)$-system.   Solutions of the first $(G,\tau,\sigma)$-system
is the equation for primitive maps studied by Burstall and Pedit \cite{BurPed94}. 

\ss
\ni {\bf Dressing actions\/}\ss

To explain the symmetries and the construction of solutions of integrable systems, we need
the dressing action of Zakharov and Shabat \cite{ZakSha79}.   Suppose $G_+, G_-$ are subgroups of $G$ and
the multiplication map from $G_+\times G_-$ to $G$ is a bijection.  
Then every $g\in G$ can be factored uniquely as $g=g_+g_-$ with $g_+\in G_+$
and $g_-\in G_-$.  Moreover, the space of right cosets
$G/G_-$  can be identified with $G_+$, so the canonical action of $G_-$ on $G/G_-$ by 
left multiplication, $g_-\cdot (gG_-)= g_-gG_-$, induces
an action $\ast$ of $G_-$ on $G_+$. The action $\ast$ is called the {\it dressing action\/}. 
The dressing action can be computed by factorization.  In fact, $g_-\ast g_+=\ti g_+$, where
$g_-g_+=\ti g_+\ti g_-$ with $\ti g_+\in G_+$ and $\ti g_-\in G_-$.   If the multiplication map
from $G_+\times G_-$ to $G$ is one-to-one but only onto an open, dense subset of $G$, then the dressing
action $\ast$ is a local action, but the corresponding Lie algebra action is global.  

\ss
\ni {\bf Iwasawa and Gauss factorizations\/}\ss

There are two well-known factorizations associated to a complex simple Lie group $G$.  The {\it Iwasawa
factorization\/} is $G=KAN$, where $K$ is a maximal compact subgroup of $G$, $A$ is abelian, and
$N$ is nilpotent.  We also refer to $G=KB$ as the Iwasawa factorization of $G$, where $B=AN$ is a
Borel subgroup.  Let $\ca$ be a Cartan subalgebra of $\cg$, $\cn_+, \cn_-$ the spaces spanned
by all positive and negative roots respectively, and $A, N_+, N_-$ the
corresponding Lie subgroups of $G$.  Then the multiplication map from $N_-\times A\times N_-$ to
$G$ is one to one and onto an open dense subset of $G$.  The set $N_-AN_+$ is called a {\it big
cell\/} of $G$.  The so-called {\it Gauss factorization\/} associated to $G$ refers to the fact that 
any $g$ in the big cell can be factorized uniquely as $n_-an_+$ with $n_\pm\in N_\pm$ and $a\in A$.    
For example, for $G=SL(n,\C)$, let $K=SU(n)$, $B_n$ the subgroup of upper triangular matrices with 
real diagonal, $A_n$ the subgroup of diagonal matrices, and $N_+(n), N_-(n)$ the subgroups of strictly 
upper and lower triangular matrices.  Then the Iwasawa factorization of $SL(n,\C)$ is $KB_n$, and the  
 Gauss factorization for the big cell is  $N_-(n)A_nN_+(n)$.  

\ss
\ni {\bf Loop group factorizations\/}\ss

We review three types of loop group factorizations that are needed for the study of symmetries of
soliton equations and elliptic integrable systems.   Let $L(G)$ denote the group of smooth
$f:S^1\to G$, $L_+(G)$ the subgroup of
$f\in L(G)$ that are the boundary values of a holomorphic map defined on $\n\l\n <1$, and $L_-(G)$ the subgroup of
$f\in L(G)$ that can be extended holomorphically to $\n\l\n>1$ in $S^2$ and $f(\infty)=e$.  
Let $U$ be a maximal compact subgroup of $G$, and $L_e(U)$ the subgroup of $f\in L(G)$ such that the image of $f$
lies in $U$ and $f(1)=e$ the identity of $G$.
\begin{enumerate}
\item[\di]
The  {\it Gauss loop group factorization\/} (or  the {\it Birkhoff factorization\/}) states that  there is an open dense
subset $L'$ of $L(G)$ such that any $g\in L'$  can be
factored uniquely as $g_+g_-$ with $g_\pm\in L_\pm(G)$. 
\item[\di] The {\it
Iwasawa loop group factorization\/}, proved in \cite{PreSeg86}, states that the multiplication map from
$L_e(U)\times L_+(G)$ to $L(G)$ is a bijection. 
\item[\di] Let $\e>0$, $\co_\e=\{\l\in \C\n \n\l\n<\e\}$, and 
$\co_{\frac{1}{\e}}=\{\l\in S^2=\C\cup \{\infty\}\n \n\l\n >1/\e\}$.   Let  $\C^*=\{\l\in \C\n
\l\not=0\}$, and
$\W^\tau(G)$ the group of holomorphic maps
$f: (\co_\e\cup \co_{1/\e})\cap \C^*\to G$ that satisfies the {\it $(G,\tau)$-reality condition\/}
$$\tau(f(1/\bar\l))= f(\l),$$
$\W^\tau_+(G)$ the subgroup of $f\in \W^\tau(G)$ that
extends holomorphically to $\C$, and $\W^\tau_-(G)$ the subgroup of $f\in \W^\tau(G)$ that
extends holomorphically to $\co_\e\cup \co_{1/\e}$.   McIntosh
proved (\cite{McI94}) that the multiplication map from $\W_-^\tau(G)\times \W_+^\tau(G)$ to
$\W^\tau(G)$ is a bijection.  
\end{enumerate}

\ni  
These loop group factorizations play central roles in the study of integrable PDEs.  

\ss\ni {\bf Solutions of soliton flows via loop group factorizations\/}\ss

Let $\co=\{\l\in\C\n \n\l\n>1/\e\}$, and $\L^\tau(G)$ the group of holomorphic maps
 $f:\co\to G$ that satisfies the $U$-reality condition $$\tau(f(\bar\l))=f(\l),$$ 
where $\tau$ is the involution on $\cg$ that defines the real form $\cu$.  Note that $f(r)\in U$ for
real $r$.   Let $\L^\tau_+(G)$
denote the subgroup of $f\in \L^\tau(G)$ that extend holomorphically to $\C$, and $\L^\tau_-(G)$
the subgroup of $f\in \L^\tau(G)$ that extend holomorphically to $1/\e<\n\l\n \leq \infty$ and $f(\infty)=e$.   The
Gauss loop group factorization implies that the multiplication map from 
$\L^\tau_+(G)\times \L^\tau_-(G)$ to $\L^\tau(G)$ is one-to-one and its image is open and dense.   

The Lax pair $\o_\l$ of a soliton flow in the 
$U$-hierarchy is a flat $\cg$-valued connection $1$-form that satisfies the $U$-reality condition
$\tau(\o_{\bar\l})=\o_\l$.  So $\o_\l (x,t)$ can be viewed as a map form
$(x,t)\in \R^2$ to the Lie algebra of $\L^\tau_+(G)$.  
Therefore the trivialization
$E_\l(x,t)$ of $\o_\l(x,t)$ can be viewed as a map from $\R^2$ to $\L^\tau_+(G)$.  
Given $g_-\in \L_-^\tau(G)$, let $\ti E(x,t)$ denote the dressing action of
$g_-$ on $E(x,t)$, i.e., $\ti E(x,t)$ is obtained using the Gauss loop group factorization to factor $g_-E(x,t)= \ti
E(x,t) \ti g(x,t)$ with $\ti E(x,t)\in
\L_+^\tau(G)$ and $\ti g(x,t)\in \L_-^\tau(G)$ for each $(x,t)$.   It can be shown that
$\ti E(x,t)$ is again a trivialization of some solution of the soliton flow.  This defines an action of $\L_-^\tau(G)$
on the space of solutions.  Moreover, 
$0$ is a solution.   If $g_-\in \L_-^\tau(G)$ is rational, then $g_-\ast 0$ can be computed explicitly and
is a rational function of exponentials.  These are the pure soliton solutions.  For general
$g_-\in \L_-^\tau(G)$, $g_-\ast 0$ is a local analytic solution of the soliton flow.  Algebraic geometric
solutions are included in the orbit $\L_-^\tau(G)\ast 0$.  To construct general rapidly decaying solutions for
the flows in the
$U$-hierarchy, we need a new type of loop group factorization.  Namely, factor
$fg$ as $\ti g \ti f$, where $f, \ti f\in L_+(G)$ so that $f_b, \ti f_b$ equal to the identity $e\in G$ at
$\l=-1$ up to infinite order and $g, \ti g$ are loops in $U$ that have essential singularity at $\l=-1$.  Here
$f_b(\l)$ and $\ti f_b(\l)$ denote the $B$-component of $f(\l)$ and $\ti f(\l)$ in the Iwasawa factorization
$G=UB$ for each $\l$.  

\ss\ni {\bf Solutions of elliptic systems via loop group factorizations\/}\ss

The Lax pair $\o_\l$ of the $m$-th $(G,\tau)$-system satisfies the
$(G,\tau)$-reality condition 
$$\tau(g(1/\bar \l))= g(\l).$$  The trivialization $E$ of $\o_\l$ is
a map from $\C$ to $\W^\tau_+(G)$.  It follows from the McIntosh loop group factorization that the dressing action
of
$\W_-^\tau(G)$ induces an action on the space of solutions of the $(G,\tau)$-systems.   Since there are constant 
solutions for the $(G,\tau)$-system,  the $\W^\tau_-(G)$-orbits through these constant solutions give rise
to a class of solutions.  But these are not all the solutions.  The $(G,\tau)$- reality
condition implies that the restriction of the trivialization $E$ of a solution to the unit circle
$\n\l\n=1$ lies in  $U$, i.e., $E$ can be viewed as a map from $\C$ to $L(U)$. 
Dorfmeister, Pedit and Wu (\cite{DPW98}) use meromorphic
maps and  the Iwasawa loop group factorization $L(G)=L_e(U)L_+(G)$ to give a method of constructing
all local solutions of the $(G,\tau)$-systems.  This is the so-called the {\it Weierstrass
representation\/}  or the {\it DPW method\/}.

\ms

Although methods of constructing solutions for both the $U$-hierarchy and the elliptic
$(G,\tau)$-systems are similar in spirit, initial data and techniques used are somewhat different. 
Moreover, while there is a canonical choice of initial data used in the factorization method to solve
soliton flows,  there is no clear canonical choice of meromorphic data for the $(G,\tau)$-hierarchy.  Since the
$(G,\tau)$-hierarchy contains the equation for harmonic maps from a domain of $\R^2$ to $U$,  the main
interest has been to understand the relation between the initial meromorphic data of the factorization
method and the global geometry.  For example, find properties of meromorphic data which corresponds
to harmonic maps from a complete surface $M$ to $U$.   This has been done when
$M$ is $S^2$ and more generally for harmonic maps of finite uniton numbers (\cite{Uhl89, BurGue97,
Gue01}), and also when
$M$ is $T^2$ (\cite{PinSte89, BFPP93}).  For a detailed survey of results concerning harmonic
maps, loop groups, and integrable systems, we refer the reader to \cite{Gue97}. 

We now turn to examples of integrable PDEs arising from geometry of maps.  When we study a
geometric problem concerning maps $f$ from a manifold $M$ to a homogeneous
space $U/U_0$, it is often useful to find a good  lifting $\ti f:M\to U$ and write down the geometric
condition imposed on the map $f$ in terms of the flat $\cu$-valued $1$-form $\ti f^{-1}d\ti f$.  
If there is a natural holomorphic deformation $F_\l:M\to G$ of such maps so that $F_0=\ti f$ and the
flatness of $F_\l^{-1}dF_\l$ for all $\l$ is equivalent to the flatness of $\ti f^{-1}d\ti f$ in some
natural coordinate system on $M$, then the corresponding geometric PDE is often an integrable system
with a zero curvature formulation.  

\ss\ni {\bf Integrable systems in differential geometry\/}\ss

One of the main interests in classical differential geometry is to find natural geometric conditions for
surfaces in $\R^3$ so that there are many explicit solutions and deformations.  It is now known that the Gauss-Codazzi
equations for surfaces with constant mean
curvature, constant Gaussian curvature, and isothermic surfaces in $\R^3$ studied by classical differential geometers
are integrable systems  and B\"acklund and Ribaucour transformations can be constructed naturally
using loop group factorizations (cf. \cite{Bob91,Cie97A,CieGolSym95, HerPed97, PinSte89,  TerUhl00a,
TerUhl00b}). 

In this paper, we review some relations between the following geometric problems and their
corresponding integrable systems:
 \begin{enumerate}
\item[(i)] The Gauss-Codazzi equations of $n$-submanifolds with
constant sectional curvature in $\R^m$, $S^m$ and hyperbolic space $\H^m$ are 
the $U/U_0$-system associated to certain real Grassmannian manifolds $U/U_0$ (cf.
\cite{BrDuPaTe02, DajToj95a, FerPed96b, Ten85, Ter97}).
\item[(ii)] The Gauss-Codazzi equations of flat Lagrangian submanifolds of $\C P^n$ is  the
$SU(n+1)/SO(n+1)$-system.
\item[(iii)]  Indefinite affine spheres in $\R^3$ are given by solutions of the $-1$-flow in the
$SL(3,\R)/\R^+$-hierarchy (\cite{Bob99}).
 \item[(iv)]  Solutions of the $-1$-flow in the $U/U_0$-hierarchy give rise to harmonic maps from the
Lorentz space $\R^{1,1}$ to $U/U_0$. (These are called sigma-models by physicists.)
\item[(v)]  The first elliptic $(G,\tau)$-system is the equation for harmonic maps from
$\R^2$ to $U$. The first elliptic $(G,\tau,\sigma)$-system is the equation for harmonic maps from $\R^2$
to the symmetric space $U/U_0$ if the order of $\sigma$ is two (\cite{BFPP93}), where $U_0$ is the
fixed point set of
$\sigma$ in $U$. 
\item[(vi)] The equation for minimal surfaces in $\C P^2$ is the first
$(SL(3,\C),\tau,\sigma)$-system, where $\tau,\sigma$ gives the $3$-symmetric space $SU(3)/T^2$
(\cite{Bur95, BoPeWo95}).  
\item[(vii)] Equations for minimal Lagrangian surfaces in $\C P^2$, minimal Legendre surfaces in $S^5$,
and minimal Lagrangian cones in $\R^6=\C^3$ are given by the first $(SL(3,\C),\tau,\sigma)$-system,
where
$\tau,\sigma$ give the $6$-symmetric space $SU(3)/ SO(2)$ (\cite{McI02}).   
\item[(viii)] The equation for Hamiltonian stationary surfaces in $\C P^2$ 
is  the second elliptic system associated to the $4$-symmetric space $SU(3)/
SU(2)$ (\cite{HelRom00}).
\end{enumerate}

\ms
Note that there may be several
geometric problems associated to one integrable system.  For example:
\begin{enumerate}
\item[\di] The SGE is the equation for surfaces in $\R^3$ with
Gaussian curvature $K=-1$, and is also the equation for harmonic maps from $\R^{1,1}$ to $S^2$.  The reason
here is that if $M$ is a surface in $\R^3$ with $K=-1$, then the second fundamental form $\II$ of $M$  is
conformally equivalent to the flat Lorentzian metric and the Gauss map $\nu:M\to S^2$ is harmonic when $M$ is
equipped with metric $\II$.
\item[\di] The $U(n)/O(n)$-system is the equation for flat Lagrangian submanifolds in $\R^{2n}$ that lie in
$S^{2n-1}$, is the equation for flat Lagrangian submanifolds in $\C P^{n-1}$, and is also the equation for flat
Egoroff metrics.  These three geometries are related as follows: the preimage of a flat Lagrangian submanifold in 
$\C P^{n-1}$ via the Hopf fibration $\pi:S^{2n-1}\to \C P^{n-1}$ is a flat Lagrangian submanifold in $\R^{2n}$ that
lies in
$S^{2n-1}$, and  the induced metrics on these flat Lagrangian submanifolds are flat Egoroff metrics.  
\end{enumerate}
\ms

Most of the integrable geometric PDEs mentioned above are either the $U/U_0$-system, the $-1$-flow, or the
$(G,\tau,\sigma)$-systems.  We would like to end this introduction by proposing a program:  Find
geometric problems whose equations are given by the
$U/U_0$-system, the $-1$-flow, or the $m$-th $(G,\tau,\sigma)$-system.   We explain the program
briefly for the 
$U/U_0$-system.   Let $U/U_0$ be a rank $n$ symmetric space, and $\o_\l$ be the corresponding family of flat
connections for a solution of the $U/U_0$-system.  We want to
find a gauge transformation $\phi$ and a value $\l=\l_0$ so that the gauge transformation $\phi\ast
\o_{\l_0}$ represents the pull back of the Maurer-Cartan form of a map from $\R^n$ to
some symmetric space $N$, whose holonomy group is $U$ or $U_0$. 
All the examples given in this paper have $U=O(m)$
or $SU(m)$.  We believe success of this program for general compact Lie group $U$ should provide new
natural classes of submanifolds in symmetric spaces and in homogeneous Riemannian manifolds with exceptional
holonomy. 

\ms

The author would like to thank Martin Guest for many helpful comments and suggestions.  
  
\bs

\section{\bf Soliton equations} \label{kb}

We review the method of constructing a hierarchy of $n\times n$ soliton flows
developed by Zakharov-Shabat \cite{ZakSha79} and Ablowitz-Kaup-Newell-Segur \cite{AKNS74}.  Their
method works equally well if we replace the algebra of $n\times n$ matrices by a general semi-simple, complex
Lie algebra $\cg$ (cf. \cite{Sat84, TerUhl98, Wil91}).  We also review the construction of
new hierarchies of flows by restricting  the
$G$-hierarchy to submanifolds naturally associated to finite order automorphisms of
$\cg$. 
Many interesting equations in differential geometry and mathematical physics are flows in these
restricted hierarchies.    

\ms
\subsection{\bf The G-hierarchy} 		\label{kc}					
\hfil\break

\ss

Let $\li\ , \ \ri$ be a non-degenerate, ad-invariant bilinear form on $\cg$, $a\in \cg$,
$\cg_a$ the centralizer of $a$ in $\cg$, and $\cg_a^\perp=\{\xi\in \cg\n \li \xi, \cg_a\ri=0.\}$.
Let $\cs(\R,\cg_a^\perp)$ denote the space  of rapidly decaying maps from $\R$ to $\cg_a^\perp$.

There is a unique family of $\cg$-valued maps $Q_{b,j}(u)$ parametrized by $\{b\in \cg\n \cg_b=\cg_a\}$
and  positive integer $j$ that satisfies the following conditions:
\begin{equation} \label{ad}
(Q_{b,j}(u))_x + [u, Q_{b, j}(u)] = [Q_{b,j+1}(u),a],
\quad Q_{b,0}(u)=b,
\end{equation}
\begin{equation} \label{db}
\sum_{j=0}^\infty Q_{b,j}(u)\l^{-j} \ {\rm is\ conjugate
\ to \ } b \
{\rm as\ an \ asymptotic\ expansion\/}.
\end{equation}
These conditions imply that $Q_{b,j}(u)$ is a polynomial in
$u, \p_x u, \ldots, \p_x^{j-1}u$ (cf., \cite{Sat84, TerUhl98}).    
 The {\it $G$-hierarchy\/} is a family of evolution equations on $\cs(\R, \cg_a^\perp)$
parametrized by
$(b,j)$, where $b\in
\cg_a$ such that $\cg_b=\cg_a$ and $j$ is a positive integer.  The $(b,j)$-flow  is

\begin{equation} \label{co}
u_t= (Q_{b,j}(u))_x+ [u,Q_{b,j}(u)]=[Q_{b,j+1}(u),a].
\end{equation}

Recall that  a
$\cg$-valued connection $1$-form $\o=\sum_{i=1}^n A_i(x)dx_i$ is {\it flat\/} if 
$$d\o= -\o\wedge\o,$$ i.e., 
$$-(A_i)_{x_j}+ (A_j)_{x_i}+ [A_i, A_j]=0, \ \  1\leq i< j\leq n.$$
The flatness of $\o$ is equivalent to the solvability of the following linear system:
\begin{equation} \label{gn}
E_{x_i}= E A_i, \ \  \ 1\leq i\leq n.
\end{equation}
Note that \eqref{gn} can also be written as $E^{-1} dE= \o$. 

\begin{defn} Let $\o$ be a flat $\cg$-valued connection $1$-form on $\R^n$.  A map $E:\R^n\to
G$ is called {\it a trivialization of $\o$\/} if $E^{-1}dE= \o$.  A trivialization $E$ of $\o$ is called the {\it
frame\/} of
$\o$ if $E$ satisfies the initial condition
$E(0)=e$, where $e$ is the identity of $G$. 
\end{defn}

 The recursive formula \eqref{ad} implies that $u$ is a solution of the
$(b,j)$-flow \eqref{co} if and only if
\begin{equation} \label{ex}
\o_\l=(a\l + u)\ dx + (b\l^j + Q_{b,1}(u)\l^{j-1} +
\cdots + Q_{b,j}(u))\ dt
\end{equation}
is a flat $\cg$-valued connection 1-form on the $(x,t)$ plane
for all $\l\in \C$.  In other words, the
$(b,j)$-flow has a Lax pair. 

\ms
 The Cauchy problem with rapidly decaying initial data for the $(b,j)$-flow \eqref{co} in the
$G$-hierarchy is solved by the inverse scattering method (cf. \cite{BeaCoi84}).  

\begin{thm} \label{beals}
 (\cite{BeaCoi84}). Suppose $a\in \cg$ such that $\cg_a$ is a maximal abelian subalgebra $\ca$ of
$\cg$.   Then there is an open dense subset
$\cs_0$ of $\cs(\R,\ca^\perp)$ such that if $u_0\in\cs_0$, then the Cauchy problem for the $(b,j)$-flow in the
$G$-hierarchy,
$$\begin{cases}
u_t= (Q_{b,j}(u))_x+[u, Q_{b,j}(u)], &\cr
u(x, 0)= u_0(x),&\cr
\end{cases}$$
  has a unique solution $u$.  Moreover, $u(x,t)$ is
defined for all
$(x,t)\in\R^2$ and $u(\cdot, t)\in \cs(\R, \ca^\perp)$.  
\end{thm}

The following is well-known, and the proof can be found in many places (cf. \cite{AblCla91, TerUhl98}).

\begin{thm} \label{hf}
 Let $X_{b,j}$ denote the vector field on $\cs(\R, \ca^\perp)$ defined by 
\begin{equation} \label{dx}
X_{b,j}(u)=(Q_{b,j}(u))_x+[u, Q_{b,j}(u)].
\end{equation}
Then  
$[X_{b,j}, X_{b',j'}]=0$ for all $b, b'\in \ca$ and positive integers $j, j'$. 
In other words, the $(b,j)$-flow commutes with the $(b',j')$-flow. 
\end{thm}

\begin{eg}\label{bn}  {\it The $SL(2,\C)$-hierarchy\/}.  

  Let $G=SL(2,\C)$, $a=\diag(i,-i)$.  Then $\cg_a=\ca=\C a$, and
$$\cg\cap\ca^\perp=\left\{\begin{pmatrix}
0&q\cr r&0\cr
\end{pmatrix} 
\bigg| \  q,r\in \C\right\}.$$
Let $u=\begin{pmatrix}
0&q\cr r&0\cr
\end{pmatrix}$.  Use \eqref{ad} and \eqref{db} and a direct computation to get the
first three terms of $Q_{a,j}(u)$:
\begin{align*}
&Q_{a,1}(u)= u, \quad  Q_{a,2}(u)=\begin{pmatrix}
\frac{iqr}{2}& \frac{iq_x}{2}\cr -\frac{ir_x}{2}&
-\frac{iqr}{2}\cr
\end{pmatrix}, \cr 
&  Q_{a,3}(u)=\begin{pmatrix}
\frac{qr_x-q_xr}{4}& -\frac{q_{xx}}{4} + {q^2r\over 2}\cr
-{r_{xx}\over 4} +{qr^2\over 2}& - {qr_x-q_xr\over 4}\cr
\end{pmatrix}.\cr
\end{align*}
Then the $(a,j)$-flow, $j=1,2, 3$, in the $SL(2,\C)$-hierarchy is the following evolution for $q$ and $r$:
\begin{align*}
&q_t=q_x, \ \ r_t= r_x,\ \ \cr
&q_t= {i\over 2} (q_{xx} - 2q^2 r), \ \  r_t= -{i\over 2} (r_{xx} - 2qr^2),\cr
&q_t= -{q_{xxx}\over 4} + {3\over 2} q r q_x, \ \  r_t= -{r_{xxx}\over 4} + {3\over 2} qr r_x.\cr
\end{align*}
\end{eg}

\ms
\subsection{\bf The $U$-hierarchy} \label{kd}
\hfil\break

Let $\tau$ be an involution of $G$ such that its differential at the identity $e$ (still denoted by
$\tau$) is a conjugate linear involution on the complex Lie algebra $\cg$, and $U$ the fixed point set of
$\tau$.  The Lie algebra $\cu$ of $U$ is a real form of
$\cg$.  Let $\cu_a$ denote the centralizer of $a$ in $\cu$, and $\cu_a^\perp$ the orthogonal complement of
$\cu_a$ in $\cu$.  Note that $\cu_a^\perp = \cg_a^\perp\cap \cu$.    It is known that the
$(b,j)$-flow in the $G$-hierarchy leaves $\cs(\R, \cu_a^\perp)$-invariant (for more detail see 
\cite{TerUhl00a}). The restriction of the flow \eqref{co} to $\cs(\R, \cu_a^\perp)$ is the
$(b,j)$-flow in the {\it $U$-hierarchy\/}.  The Lax pair $\o_\l$ defined by \eqref{ex} is a $\cg$-valued
$1$-form, and $\o_\l$ satisfies the {\it
$U$-reality condition\/}:
\begin{equation}\label{fp}
\tau(\o_{\bar\l})= \o_\l.
\end{equation}
Note that $\xi=\sum_j \xi_j \l^j$ satisfies the $U$-reality condition if and only if $\xi_j\in \cu$ for all
$j$.

\begin{eg}\label{es}
{\it The $SU(2)$-hierarchy\/}. 

Let $\tau$ be the involution of $sl(2,\C)$ defined by
$\tau(\xi)= -\bar\xi^t$. Then the fixed point set of $\tau$ is the real form $\cu= su(2)$.  Let
$a=\diag(i, -i)\in \cu$. Then  $\cu_a=\ca=\R a$, and 
$$\ca^\perp\cap \cu=\left\{\begin{pmatrix}
0&q\cr -\bar q&0\cr
\end{pmatrix}
\bigg| \ q\in \R\right\}.$$ So
$C(\R,\ca^\perp\cap\cu)$ can be identified as $C(\R,\C)$, and the $SU(2)$-hierarchy is the
restriction of the $SL(2,\C)$-hierarchy to the subspace $r=-\bar q$.   The first three flows in the
$SU(2)$-hierarchy are
\begin{align*}
&q_t=q_x, \cr 
&q_t= {i\over 2} (q_{xx} + 2 \n q\n^2 q),\cr 
&q_t=-{1\over 4}  (q_{xxx} + 6 \n q\n^2 q).\cr
\end{align*}
Note that the $(a,2)$-flow in the $SU(2)$-hierarchy is the NLS.  
\end{eg}

\ms
\subsection{\bf The $\sigma$-twisted $G$-hierarchy}\label{ke}
\hfil\break

Let $\sigma$ be an order $k$ group automorphism of $G$ such that its differential at the identity $e$
(still denoted by
$\sigma$) is a complex linear Lie algebra homomorphism of $\cg$.   Let 
$$\cg = \cg_0 +\cdots + \cg_{k-1},$$
where $\cg_j$ is the eigenspace of $\sigma$ with eigenvalue $e^{2j\pi i\over k}$.   
Note that $\cg_i=\cg_j$ if $i\equiv j$ mod $k$, and 
$$[\cg_j,\cg_r]\subset \cg_{j+r}.$$
 Let $\ca$ be a
maximal abelian subspace in $\cg_1$, and $a\in\ca$ {\it regular in $\cg_1$\/}, i.e.,
$$\{x\in \cg_1\n [x,a]=0\} = \ca.$$
  It is known (cf.  \cite{TerUhl00a}) that if the image of
$u$ lies in $\cg_0\cap\cg_a^\perp$, then 

\begin{equation}\label{fo}
Q_{b,j}(u)\in \cg_{1-j}.
\end{equation}
  Since $a\in \cg_1$,  the right hand side of the
$(b,mk+1)$-flow 
$$(Q_{b,mk+1}(u))_x + [u, Q_{b,mk+1}(u)]= [Q_{b,mk+2},a]\in 
\cg_{-mk}=\cg_0.$$  In other words, the $(b,mk+1)$-flow leaves $\cs(\R,
\cg_a^\perp\cap\cg_0)$ invariant.  {\it The $\sigma$-twisted $G$-hierarchy\/} is the restriction of the
$(b,mk+1)$-flow in the
$G$-hierarchy to $\cs(\R, \cg_a^\perp\cap \cg_0)$ for $m=1, 2, \cdots$.   

It follows from \eqref{fo} that
the Lax pair of the $(b, mj+1)$-flow in the $\sigma$-twisted $G$-hierarchy satisfies the {\it $(G,\sigma)$-reality
condition\/}:
\begin{equation}\label{hd}
\sigma(\o_{e^{-{2\pi i\over k}}\l}) = \o_\l.
\end{equation}
Note that $\xi=\sum_j \xi_j\l^j$ satisfies the $(G,\sigma)$-reality condition if and only if $\xi_j\in
\cg_j$ for all $j$.

\begin{eg} \label{fb} {\it Kupershmidt-Wilson hierarchy\/} (\cite{KupWil81}).  

Let  $G= SL(n,\C)$, and  $\sigma$ the order
$n$ automorphism of
$SL(n,\C)$ defined by $\sigma(g)= C^{-1}gC$,
where $C=e_{21} +e_{32} +\cdots + e_{n, n-1} + e_{1n}$ is the permutation matrix $(12\cdots n)$. Here $e_{ij}$
denote the elementary matrix of $gl(n)$.  The eigenspace $\cg_k$ of $\sigma$ on $sl(n,\C)$ with eigenvalue
$\a=e^{2\pi i k\over n}$ is the space of all 
$y=(y_{ij})$ such that
$y_{i+1,j+1}=\a^k y_{ij}$ for all $1\leq i,j\leq n$.
Let $a=\diag(1,\a,\cdots, \a^{n-1})\in \cg_1$, and $\ca=\C a$.  Then $\ca$ is a maximal abelian
subalgebra of
$\cg_1$.  
 The $(SL(n,\C),\sigma)$-hierarchy is the restriction of the $(jn+1)$-th flow in the
$sl(n,\C)$-hierarchy to $\cs(R, \cg_0\cap \cg_a^\perp)$.   For example, for
$n=2$, 
$$\cg_0\cap \cg_a^\perp= \left\{
\begin{pmatrix}
0&q\cr q&0\cr
\end{pmatrix}
\bigg|\ \  q\in \C\right\}.$$
The first flow is the translation $q_t= q_x$,   and the third flow is the complex modified KdV
$$q_t = {1\over 4} (q_{xxx}- 6q^2 q_x).$$
\end{eg}

\ms
\subsection{\bf The $U/U_0$-hierarchy} \label{kf}
\hfil\break

Let $\tau$ be a conjugate linear involution of $\cg$, $\cu$ its fixed point set, and
$\sigma$ a complex linear, order $k$ automorphism of $\cg$ such that 
$$\tau\sigma= \sigma^{-1}\tau^{-1}=\sigma^{-1}\tau.$$
Let $\cg_j$ denote the eigenspace of $\sigma$ with eigenvalue $e^{2\pi i j\over k}$.  
We claim that $\tau(\cg_j)\subset \cg_j$. To see this, let $\xi_j\in \cg_j$.  Then
$$\sigma(\tau(\xi_j))= \tau(\sigma^{-1}(\xi_j))= \tau(\a^{-j} \xi_j) = \overline{\a^{-j}}\tau( \xi_j)=
\a^j\tau(\xi_j),$$ where $\a=e^{2\pi i\over k}$, proving the claim.  
Let $\cu_j= \cg_j\cap \cu$.  Then we have
$$\cu=\cu_0+\cdots + \cu_{k-1}.$$
Let $\ca\subset \cu_1$ be a maximal abelian subspace in $\cu_1$.  An element $a\in \ca$ is
{\it regular in $\cu_1$\/} if 
$$\{\xi\in \cu_1\n [\xi,a]=0\}=\ca.$$
Let $b\in \ca$, and $u\in\cg_a^\perp\cap\cu_0$.   Then
$Q_{b,j}(u)\in \cu_{1-j}$ for all $j\geq 0$ (\cite{TerUhl00a}).  So the
 $(b, mk+1)$-flow in the $\sigma$-twisted $G$-hierarchy leaves $\cs(\R,\cg_a^\perp\cap \cu_0)$ invariant.  The
restriction of these flows to $\cs(\R, \cg_a^\perp\cap \cu_o)$ is called the  {\it $U/U_o$-hierarchy\/}. 

The Lax pair $\o_\l$ of the $(b,mk+1)$-flow in the $U/U_0$-hierarchy satisfies the {\it
$U/U_o$-reality condition\/}:
\begin{equation}\label{cb}
\tau(\o_{\bar\l})= \o_\l, \ \ \sigma(\o_{ \l})= \o_{e^{2\pi i\over k}\l}.
\end{equation}
Note that $\xi=\sum_j \xi_j \l^j$ satisfies the $U/U_0$-reality condition if and only if  $\xi_j\in \cu_j$
for all $j$.  

When the order of $\sigma$ is $2$, the condition $\tau\sigma=\sigma^{-1}\tau^{-1}$ implies that 
 $\tau$ and $\sigma$ commute, $U/U_0$ is a
symmetric space, and  $\cu=\cu_0+\cu_1$ is a Cartan decomposition for the symmetric space
$U/U_0$.  

\begin{eg}\label{et}  {\it The $SU(2)/SO(2)$-hierarchy\/}. 

Let $\tau(\xi)=-\bar\xi^t$ and $\sigma(\xi)= -\xi^t$ be involutions of
$sl(2,\C)$ that give the symmetric space $SU(2)/SO(2)$.  Then 
$$\cu_0= so(2), \ \ \cu_1=\{i \xi\n \xi \ {\rm is \ a } \  2\times
2 \ {\rm real\  symmetric\  matrix}\}$$ with $SU(2)/SO(2)$ the corresponding symmetric space. Let
$a=\diag(i,-i)$. Then $\cg_a^\perp\cap \cu_0= so(2)$. So
$\cs(\R, \cg_a^\perp\cap\cu_0)$ can be identified as $\cs(\R, \R)$.  The $(a,1)$- and $(a,3)$-flow in the
$SU(2)/SO(2)$-hierarchy are
\begin{align*}
&q_t= q_x,\cr 
&q_t=- {1\over 4}(q_{xxx} + 6q^2 q).\cr
\end{align*}
Note that the $(a,3)$-flow is  {\it modified KdV equation\/} (mKdV).
\end{eg}

\begin{eg}\label{eu}  {\it The $SU(n)/SO(n)$-hierarchy\/}. 

 Let $\tau$ and $\xi$ be involutions of
$sl(n,\C)$ defined by
$$\tau(\xi)=-\bar\xi^t, \ \ \sigma(\xi)= -\xi^t.$$
Then $\tau\sigma=\sigma\tau$, $\cu=su(n)$,  $\cu_0= so(n)$, and $\cu_1$ is the space of
$iY\in su(n)$, where $Y$ is real and symmetric. 
The corresponding symmetric space is 
${U\over U_0}={SU(n)\over SO(n)}$.    Let $\ca$ denote the space of diagonal
matrices in $su(n)$. Then $\ca$ is a maximal abelian linear subspace in $\cu_1$, and $\ca^\perp\cap
\cu_0=so(n)$.  An element 
$a=i\diag(a_1, \cdots, a_n)$ is {\it regular in $\cu_1$\/} if $a_1, \cdots, a_n$ are distinct. 
Let $a\in \ca$ be a regular element, and $b=i\diag(b_1, \cdots, b_n) \in\ca$.   The
$(b,1)$-flow in the ${SU(n)\over SO(n)}$-hierarchy on
$\cs(\R,so(n))$ is the {\it reduced
$n$-wave equation\/}
$$(u_{ij})_t= {b_i-b_j\over a_i-a_j}\ (u_{ij})_x + \sum_k \left({b_k-b_j\over a_k-a_j}
-{b_i-b_k\over a_i-a_k}\right)\ u_{ik}u_{kj}, \quad i\not=j.$$
\end{eg}

\begin{eg} \label{bo} Let $U/U_0$ be a symmetric space, $\cu=\cu_0+\cu_1$ a Cartan
decomposition, 
$\ca$ a maximal abelian subspace in $\cu_1$, $a\in \ca$ regular, and $b\in \ca$. Note
$\ad(a)^{-1}$ maps
$\cu_a^\perp\cap \cu_0$ and $\cu_a^\perp\cap\cu_1$ isomorphically onto $\cu_a^\perp\cap \cu_1$
and $\cu_a^\perp\cap \cu_0$ respectively.    So
$\ad(b)\ad(a)^{-1}(\cu_a^\perp\cap \cu_0) \subset \cu_a^\perp\cap\cu_0$.   The recursive formula
\eqref{ad} implies that
\begin{equation}\label{ada}
Q_{b,1}(u)=\ad(b)\ad(a)^{-1}(u).
\end{equation}
   So  the
$(b,1)$-flow in the
$U/U_0$-hierarchy is the equation for maps $u:\R^2\to \cu_a^\perp\cap \cu_0$:
\begin{equation}\label{ac}
u_t= \ad(b)\ad(a)^{-1}(u_x) + [u, \ad(b)\ad(a)^{-1}(u)].
\end{equation}
This is the {\it reduced $n$-wave equation associated to $U/U_0$\/}, which  has a Lax pair
\begin{equation}\label{ae}
\o_\l=(a\l + u)dx + (b\l + \ad(b)\ad(a)^{-1}(u)) dt.
\end{equation}
\end{eg}

\ms
\subsection{\bf The $U/U_0$-system}  \label{kg}
\hfil\break

Let $\tau$ be a conjugate linear involution of $\cg$, $\sigma$ a complex linear involution of $\cg$ such that
$\tau\sigma=\sigma\tau$, $\cu$ the fixed point set of $\tau$, and $U_0$ the subgroup of $U$ fixed by
$\sigma$.  Let $\cu=\cu_0+\cu_1$ denote the Cartan decomposition of the symmetric space $U/U_0$.  Let
$\ca$ be a maximal abelian linear subspace of $\cu_1$, and $a_1, \cdots, a_n$ a basis of
$\ca$.    {\it The $U/U_0$-system\/}  is   the following system for
$v:R^n\to \cu_{\ca}^\perp\cap\cu_1$:
\begin{equation}\label{aa}
[a_i, v_{x_j}] - [a_j, v_{x_i}]=[[a_i, v], [a_j,v]], \quad i\not=j.
\end{equation}
It has a Lax $n$-tuple,
\begin{equation}\label{ax}
\o_\l=\sum_{i=1}^n (a_i\l + [a_i, v]) dx_i,
\end{equation}
which satisfies the $U/U_0$-reality condition \eqref{cb}. Moreover, the following statements are
equivalent for smooth map
$v:\R^n\to
\cu_{\ca}^\perp\cap \cu_1$:
\begin{enumerate}
\item[(i)] $v$ is a solution of   the $U/U_0$-system \eqref{aa},
\item[(ii)] $\o_\l$ defined by \eqref{ax}  is a flat
$\cg$-valued connection $1$-form on $R^n$ for all $\l\in C$,
\item[(iii)] $\o_r$ is flat for some $r\in \R$. \end{enumerate}

\ms
We claim that the $U/U_0$-system is independent of the choice of basis of $\ca$.  If  $b_1, \cdots, b_n$
is a basis of $\ca$, then there exists a constant matrix $(c_{ij})$ such that   $a_i =\sum_{j=1}^n c_{ij}
b_j$.  The $U/U_0$-system defined by the new base $b_1, \cdots, b_n$ is
$$[b_i, v_{y_j}] - [b_j, v_{y_i}]=[[b_i, v], [b_j,v]].$$ This is the same system as \eqref{aa} if we change
coordinates $y_i= \sum_{j=1}^n c_{ji} x_j$. 

The $U/U_0$-system is given by the first commuting $n$-flows in the $U/U_0$-hierarchy, i.e., 

\begin{prop} \label{bp} (\cite{Ter97}).  With the same notation as above, let $a_1, \cdots, a_n$ be a
basis of $\ca$ such that
$a_1, \cdots, a_n$ are regular. Let $a=a_1$. Then $v:\R^n\to \cu_{\ca}^\perp\cap \cu_1$ is a
solution of the $U/U_0$-system \eqref{aa} if and only if $u(x)=[a,v(x)]$ satisfies the $(a_j, 1)$-flow in
the $U/U_0$-hierarchy, 
$$u_{x_j}= \ad(a_j)\ad(a)^{-1}(u_{x_1}) + [u, \ad(a_j)\ad(a)^{-1}(u)],$$
for all $1\leq j\leq n$.  
\end{prop}

As a consequence of Theorem \ref{beals} and Proposition \ref{bp} we
have

\begin{cor} \label{hx} (\cite{Ter97}).  Suppose $a=a_1\in\ca$ is regular in $\cu_1$.  Then there exists an
open dense subset
$\cs_0$ of $\cs(\R, \cu_a^\perp \cap \cu_1)$ such that given any $v_0\in
\cs_0$ there exists a unique solution $v$ of \eqref{aa} defined for all $x\in \R^n$ such that
$v(x_1,0,\cdots, 0)= v_0(x_1)$ and $v(\cdot, x_2, \cdots, x_n)\in \cs(\R, \cu_a^\perp \cap  \cu_1)$.  
\end{cor}

Next we give some examples. 

\begin{eg}\label{hy} The $U$-system.   

Let $\tau$ be a conjugate linear involution of $\cg$, and $\cu$ the
fixed point set of $\tau$.  Let $\tau_2(x,y)=(\tau(x),
\tau(y))$ and $\sigma(x,y)=(y,x)$ be involutions of $\cg\times\cg$. Then
$\tau_2\sigma=\sigma\tau_2$, and the corresponding symmetric space is
$(U\times U)/\D(U)\simeq U$, where $\D(U)$ is the diagonal group $\{(g,g)\n g\in U\}$.  The
 $(U\times U)/\D(U)$-system is the $U$-system \eqref{aa} for maps $v:\R^n\to
\ca^\perp\cap \cu$, where
$\ca$ is a maximal abelian subalgebra of $\cu$ and 
$\{a_1, \cdots, a_n\}$ is a basis of $\ca$.
\end{eg}

\begin{eg}\label{ai} The  ${O(2n)\over O(n)\times O(n)}$-system.   

Here $U/U_0$ is the
symmetric space
${O(2n)\over O(n)\times O(n)}$, 
$G=O(2n,\C)$, $\tau(g)= \bar g$, $\sigma(g)= \I_{n,n} \ g\ \I_{n,n}^{-1}$, where $I_{n,n}$ is the
diagonal matrix with $a_{ii}=1$ for $1\leq i\leq n$ and $a_{ii}=-1$ for $n+1\leq i\leq 2n$.  So
$\cu= so(2n)$,
$\cu_0=so(n)+ so(n)$,  and 
$$\cu_1=\left\{\begin{pmatrix}
0& \xi\cr -\xi^t& 0\cr
\end{pmatrix}
 \bigg| \ \xi \in
gl(n,\R)\right\}.$$ The linear subspace $\ca$ spanned by 
$$\{a_i= -e_{i, n+i}+e_{n+i, i}\n 1\leq i\leq n\}$$
is a maximal abelian subspace of $\cu_1$, and 
$$\cu_1\cap\ca^\perp=\left\{\begin{pmatrix}
0& F\cr -F^t& 0\cr
\end{pmatrix}\ \bigg|  F=(f_{ij})\in gl(n,\R),
f_{ii}=0 \ {\rm for\ } 1\leq i\leq n.\right\}.$$
The corresponding $U/U_0$-system \eqref{aa} written in terms of $F=(f_{ij})$ is
\begin{equation}\label{ag}
\begin{cases}(f_{ij})_{x_i} + (f_{ji})_{x_j} + \sum_k f_{ki} f_{kj} =0, &\text{ if
$i\not=j$,}\cr (f_{ij})_{x_j} + (f_{ji})_{x_i} + \sum_k f_{ik} f_{jk} =0, & \text{if $i\not=j$,}\cr
(f_{ij})_{x_k}= f_{ik}f_{kj}, & \text{if $i, j, k$ are distinct.}\cr
\end{cases}
\end{equation}
The Lax $n$-tuple $\o_\l$ \eqref{ax}, written in matrix form is
\begin{equation}\label{dy}
\o_\l= 
\begin{pmatrix}
\d F^t- F\d & -\l\d \cr \l\d & -F^t\d + \d F\cr
\end{pmatrix}, \quad {\rm
where\ } \d=
\diag(dx_1, \cdots, dx_n).
\end{equation} 
Note that the first and the third equations of \eqref{ag} imply that $\d F^t- F\d$ is flat, and the
second and third equations of \eqref{ag} imply that $-F^t\d + \d F$ is flat.  
\end{eg}

\begin{eg}\label{am}  The ${U(n)\over O(n)}$-system. 

 Here $\cg=gl(n,\C)$, $\tau(\xi)=-\bar \xi^t$, and
$\sigma(\xi)=-\xi^t$.  Then  $\cu= u(n)$, $\cu_0= o(n)$, and 
$$\cu_1=\{i F\n F=(f_{ij})\in gl(n,\R), f_{ij}=f_{ji}\}.$$
The linear subspace $\ca$ spanned by 
$$\{a_j= i e_{jj} \n 1\leq j\leq n\}$$
is a maximal abelian subspace of $\cu_1$, and 
$$\cu_1\cap\ca^\perp=\{i F\n F=(f_{ij})\in gl(n,\R), \ f_{ij}=f_{ji},\ f_{ii}=0 \ {\rm for\ } 1\leq
i, j\leq n\}.$$ The corresponding $U/U_0$-system written in terms of $F$ is the restriction of
 system \eqref{ag} to the linear subspace of $F=(f_{ij})$ such that $f_{ij}=f_{ji}$. So the
${U(n)\over O(n)}$-system is  the system for symmetric $F=(f_{ij})$:
\begin{equation}\label{an}
\begin{cases}
(f_{ij})_{x_i} + (f_{ij})_{x_j} + \sum_k f_{ik} f_{jk} =0, & \text{if
$i\not=j$,}\cr (f_{ij})_{x_k}= f_{ik}f_{kj}, &\text{ if $i, j, k$ are distinct.}\cr
\end{cases}
\end{equation}
Or equivalently, $[\d, F]=\d F- F\d$ is flat. 
\end{eg}

\begin{eg} \label{bt}  The ${SU(n)\over SO(n)}$-system.  

Here $\cg= sl(n,\C)$, 
$\tau(\xi)= -\bar
\xi^t$, and $\sigma(\xi)= -\xi^t$ for $\xi\in \cg$.   Then
$\cu=su(n)$, $\cu_0= so(n)$, and 
$$\cu_1=\{i F\n F=(f_{ij})\in gl(n,\R), f_{ij}=f_{ji}, \ \sum_{i=1}^n f_{ii}=0\}.$$
The linear subspace $\ca$ spanned by 
$$\{b_j= i(e_{jj}- e_{11}) \n 2\leq j\leq n\}$$
is a maximal abelian linear subspace of $\cu_0$, and 
$$\ca^\perp\cap \cu_1=\{i F\n F=(f_{ij})\in gl(n,\R), f_{ij}=f_{ji}, f_{ii}=0 \ \ {\rm for\ \ } 1\leq i,j\leq
n\}.$$
The ${SU(n)\over SO(n)}$-system is 
\begin{equation}\label{bu}
[b_i, F_{t_j}]-[b_j, F_{t_i}]= [b_i, F], [b_j, F]]. \quad 2\leq i\not=j\leq n.
\end{equation}
\end{eg}

\begin{eg} \label{eb} The $U/U_0= \frac{GL(2,\H)}{
(\R^+\times SU(2))^2}$-system (\cite{Bob00p}).  

Here $\cg=gl(4,\C)$.  For $X\in gl(4,\C)$,
write
$$X=\begin{pmatrix}X_1&X_2\cr X_3&X_4\cr\end{pmatrix}, \ \  X_i\in gl(2,\C).$$
Let $\tau$ be the involution of $\cg$ defined by
$$\tau(X)= \begin{pmatrix}J&0\cr 0& J\cr\end{pmatrix}\overline{ X}
\begin{pmatrix}J^{-1}&0\cr 0&
J^{-1}\cr\end{pmatrix},\
\
\ {\rm where
\ } J=\begin{pmatrix}0&1\cr -1&0\cr\end{pmatrix}.$$  The fixed point set $\cu$ of $\tau$ is the
subalgebra of
$X\in
\cg$  such that $J\bar X_iJ^{-1}= X_i$ for all $1\leq i\leq 4$, i.e., $X_i$ lies in the fixed point set of the
involution of $gl(2,\C)$ defined by $\tau_0(Y)= J\bar Y J^{-1}$.   A direct computation
implies that the fixed point set of $\tau_0$ is 
$$\left\{\begin{pmatrix}z_1& z_2\cr -\bar z_2 & \bar z_1\cr\end{pmatrix}\ \bigg| \ z_1, z_2\in
\C\right\} =
\R\times su(2).$$  Note that
$\R\times su(2)$ is isomorphic to  the quaternions $\H$ as  associative algebras via the
isomorphism 
$${\bf i\ }=\begin{pmatrix} i&0\cr 0&-i\cr\end{pmatrix}, \ {\bf j}\ = \begin{pmatrix} 
0&1\cr -1&0\cr\end{pmatrix}, \
{\bf k}\ =\begin{pmatrix}
0&i\cr i&0\cr\end{pmatrix}, \ {\bf 1}\ =\begin{pmatrix}1&0\cr 0&1\end{pmatrix}.$$
(It is easy to check that ${\bf ij\/} = {\bf k\/}$, ${\bf jk\/}= {\bf i\/}$, ${\bf ki\/}={\bf j\/}$).  
 So we can view $\cu=gl(2,\H)$, i.e., the algebra of $2\times 2$ matrices with entries in
 the quaternions $\H$.  

Let $\sigma$ be the
involution on $gl(4,\C)$ defined by $$\sigma(Y)= \begin{pmatrix}{\bf 1\/}& 0 \cr 0& -{\bf
1\/}\cr\end{pmatrix} Y
\begin{pmatrix}{\bf 1\/}& 0 \cr 0& -{\bf 1\/}\cr\end{pmatrix}^{-1}.$$  Then $\sigma\tau=\tau\sigma$,
and
\begin{align*}
\cu_0&=\left\{\begin{pmatrix}
P&0\cr 0&Q\cr\end{pmatrix}
\bigg| \ P, Q\in \R\times su(2)\right\},\cr
\cu_1&=\left\{\begin{pmatrix} 
0&P\cr Q&0\cr\end{pmatrix}\bigg| \ P, Q\in \R\times su(2)\right\}.
\end{align*}
Let
\begin{equation}\label{by} 
a_1=\begin{pmatrix}
0&{\bf k\/}\cr {\bf k\/}&0\cr
\end{pmatrix} \ \ a_2=\begin{pmatrix}
0 &{\bf
j\/}\cr -{\bf j\/}& 0\cr\end{pmatrix}.
\end{equation}
The space $\ca$ spanned by $a_1$ and $a_2$ is a maximal abelian linear subspace of $\cu_1$, and
$\ca^\perp\cap \cu_1$  is the space of all matrices of the form
\begin{equation}\label{bx}
v=\begin{pmatrix}
0& \begin{pmatrix}p_1& p_2\cr -\bar p_2& \bar p_1\cr\end{pmatrix}\cr
\begin{pmatrix}q_1&
\bar p_2\cr -p_2& \bar q_1\cr\end{pmatrix} 
&0\cr\end{pmatrix}.
\end{equation} 
The ${GL(2,\H)\over (\R\times SU(2))^2}$-system is
\begin{equation}\label{dq}
\begin{cases}
(p_2)_y+ i(p_2)_x= -\n p_1\n^2 +\n q_1\n^2,&\cr
(\bar q_1+p_1)_x = 2i \ (\bar p_2-p_2) (\bar q_1 -p_1),&\cr
(\bar q_1-p_1)_y = -2 (p_2+\bar p_2)(\bar q_1 + p_1).&\cr
\end{cases}
\end{equation}
Its  Lax pair \eqref{ax}  is
\begin{align*}
\o_\l &=\begin{pmatrix}
\begin{pmatrix}Z& W\cr -\bar W&\bar Z \end{pmatrix} 
& {\bf k\/}\l\cr {\bf k\/} \l& 
\begin{pmatrix}-\bar Z& \bar W\cr -W& - Z\end{pmatrix}\end{pmatrix}
 \ dx\cr &\ \  +
\begin{pmatrix}
\begin{pmatrix}-iZ& Y\cr -\bar Y& i\bar Z\cr
\end{pmatrix}
& {\bf  j\/}\l \cr -{\bf j\/} \l & 
\begin{pmatrix}-i \bar Z& -\bar Y\cr Y& iZ\cr
\end{pmatrix}\cr
\end{pmatrix}\ dy,\cr
\end{align*}
 where 
\begin{equation}\label{dp}
Z=-2ip_2, \ \ W= i(\bar q_1-p_1),  \ \  Y= (\bar q_1 + p_1).
\end{equation}
Let $$p_2= \b_1+ i\b_2, \ \ \ \b_1, \b_2\in \R.$$  

Equate the imaginary part of the first equation in
system \eqref{dq} to get 
$$(\b_2)_y + (\b_1)_x=0.$$
So there exists $u$ such that $\b_1= -{u_y\over 8}$ and $\b_2={u_x\over 8}$, and hence
$$p_2={-u_y+ iu_x\over 8}.$$
Substitute this into \eqref{dq} to get 
\begin{equation}\label{bz}
\begin{cases}
u_{xx} + u_{yy} = 8(\n p_1\n^2 - \n q_1\n^2), &\cr
(\bar q_1+ p_1)_x = {u_x\over 2} (\bar q_1 - p_1),&\cr
(\bar q_1 -p_1)_y= {u_y\over 2} (\bar q_1 + p_1).&\cr
\end{cases}
\end{equation}

System \eqref{bz} is gauge equivalent to system \eqref{dq}.  To see this, we recall that $v$ is a solution
of the
$U/U_0$-system \eqref{aa} if and only if  $\o_0$
 is flat. So $\o_0$ is a $(\R\times su(2))\times
(\R\times su(2))$-valued flat connection $1$-form.  The $\R\times \R$-component of $\o_0$ is
$$\o_0^0=\begin{pmatrix}
2\b_2{\bf 1\/}& 0\cr 0& -2\b_2{\bf 1\/}\cr\end{pmatrix}\ dx + \begin{pmatrix}
-2\b_1{\bf 1\/}&
0\cr 0&  2\b_1 {\bf 1\/}\cr\end{pmatrix}\ dy= \begin{pmatrix}
{du\over 4}\ {\bf 1\/}&0\cr 0&
-{du\over 4}\ {\bf 1\/}\cr\end{pmatrix}.$$ Let
$$g=\begin{pmatrix} 
e^{u/4}{\bf 1\/}& 0\cr 0&e^{-u/4}{\bf 1\/}\cr\end{pmatrix}.$$
The gauge transformation of $\o_0$ by $g$ is
\begin{align*}
&g\ast \o_0  = \begin{pmatrix}
\tau_1&0\cr 0& \tau_2\cr
\end{pmatrix}, \ \ {\rm
where\ \ } \cr & \tau_1= \begin{pmatrix}
{i(u_y dx - u_x dy)\over 4} &
i(\bar q_1-p_1)dx + (\bar q_1+ p_1) dy\cr  i(q_1-\bar p_1)dx -(q_1+\bar p_1) dy& -{i(u_y dx - u_x
dy)\over 4}\cr\end{pmatrix},
\cr &\tau_2= 
\begin{pmatrix}{i(u_y dx - u_x dy)\over 4} & -i(q_1-\bar p_1) dx- (q_1 +\bar p_1) dy\cr -i(\bar
q_1-p_1) dx + (\bar q_1+p_1) dy& -{i(u_y dx - u_x dy)\over 4} \end{pmatrix}. 
\end{align*}
The connection $g\ast\o_0$ is flat if and only if $(u, p_1, q_1)$ is a solution of \eqref{bz}.  So
$g\ast \o_\l$ is a Lax pair of \eqref{bz}.   In other words, system
\eqref{bz} is gauge equivalent to  the ${GL(2,\H)\over (\R^+\times SU(2))^2}$-system. 

Suppose $(u,p_1,q_1)$ is a solution of \eqref{bz} and $q_1$ is real.  Let $p_1=B_1+iB_2$.  Equate the
imaginary part of the second and third equations of \eqref{bz} to get
$$(B_2)_x=-u_xB_2/2, \ \  (B_2)_y=-u_yB_2/2.$$
So $B_2= ce^{-u/2}$ for some constant $c$, and \eqref{bz} becomes the following system for real
functions $u, q_1, B_1$. 
\begin{equation}\label{fr}
\begin{cases}
u_{xx}+u_{yy}=8(c^2e^{-u} + B_1^2-q_1^2),&\cr 
(q_1+ B_1)_x= {u_x\over 2} (q_1-B_1),&\cr
(q_1-B_1)_y= {u_y\over 2} (q_1+B_1).&\cr\end{cases}
\end{equation}
If $p_1$ is also real, i.e., $c=0$ in \eqref{fr}, then system \eqref{fr} becomes the following system for
real functions $u, q_1, p_1$:
\begin{equation} \label{fs}
\begin{cases}
u_{xx}+u_{yy}=8(p_1^2-q_1^2),&\cr 
(q_1+ p_1)_x= {u_x\over 2} (q_1-p_1),&\cr
(q_1-p_1)_y= {u_y\over 2} (q_1+p_1).&\cr
\end{cases}
\end{equation}
\end{eg}

\ms
\subsection{\bf The $-1$-flow\/} \label{kh}
\hfil\break

Let $\cu$ be the real form defined by the involution $\tau$ on $\cg$,  $a, b\in \cu$ such that $[a,b]=0$.  The {\it
$-1$-flow associated to
$U$ defined by $a,b$\/} is the following system for $u:\R^2\to \cu\cap\cu_a^\perp$:
\begin{equation}\label{ew}
u_t=[a, g^{-1}bg], \ \  {\rm where\ } g^{-1}g_x= u.
\end{equation}
This system has a Lax pair 
$$\o_\l= (a\l + u) dx + \l^{-1}g^{-1} bg\ dt, \ \  {\rm where \ } g:\R^2\to U, \ \ 
g^{-1}g_x=u.$$\label{fg}
Note that \eqref{fg} satisfies the $U$-reality condition \eqref{fp}.  

\begin{thm} \label{ey} (\cite{Ter97}).  The $-1$-flow \eqref{ew} commutes with all the flows in the
$U$-hierarchy.  
\end{thm}

Let $\sigma$ be an order $k$ automorphism of $\cg$ such that $\tau\sigma=\sigma^{-1}\tau^{-1}$ and
$\cu=\cu_0+\cdots +\cu_{k-1}$ as in section \ref{kf} Let  $a\in \cu_1$ a regular element, and
$b\in \cu_{k-1}$ such that
$[a,b]=0$, then the right hand side of \eqref{ew} is a vector field on $C(\R, \cu_0\cap\cu_a^\perp)$.
The flow \eqref{ew} restricted to the space $C(\R, \cu_0\cap\cu_a^\perp)$ of smooth maps from $\R$
to
$\cu_0\cap \cu_a^\perp$ is  {\it the
$-1$-flow associated to $U/U_0$\/}, and $\o_\l$ defined by \eqref{fg} is its Lax pair that satisfies the
$U/U_0$-reality condition \eqref{cb}.   

We can also write the $-1$-flow \eqref{ew} associated to $U$ ($U/U_0$ resp.) as an equation for
$g:\R^2\to U$ ($g:\R^2\to U_0$ respectively):
\begin{equation} \label{cu}
(g^{-1}g_x)_t= [a, g^{-1} bg],
\end{equation}
where $g^{-1}g_x\in\cu_a^\perp$ ($\in  \cu_0\cap \cu_a^\perp$ respectively).  
Its Lax pair  is
\begin{equation}\label{cv}
\o_\l= (a\l + g^{-1}g_x)dx + \l^{-1} g^{-1} b g dt.
\end{equation}

\ms

\begin{eg}\label{ec} The $-1$-flow associated to $SU(2)/SO(2)$. 

Let $a=\diag(i, -i)$, and $b=-{a\over 4}$.  Let
$g=\begin{pmatrix}
\cos {q\over 2}& \sin  {q\over 2}\cr -\sin {q\over 2} & \cos {q\over
2}\cr \end{pmatrix}$.  Then $u=g^{-1}g_x={1\over 2}\ \begin{pmatrix}
0&q_x\cr
-q_x&0\cr \end{pmatrix}$, and the $-1$-flow \eqref{cu} associated to
$SU(2)/SO(2)$  is the sine-Gordon equation (SGE):
$$q_{xt}= \sin q,$$
and its Lax pair is
$$\o_\l= \begin{pmatrix} i\l & {q_x\over 2}\cr -{q_x\over 2}& -i\l\cr 
\end{pmatrix} dx -{i\l\over 4}
\begin{pmatrix} \cos q&
\sin q\cr
\sin q&\cos q\cr \end{pmatrix}dt.$$
\end{eg}

\ms
\beg\label{cm} The $-1$-flow associated to $U/U_0 = SL(3,\C)/\R^+$.
\eeg

Here $\R^+$ is the subgroup  $$\R^+=\{\diag(r,
r^{-1},1)\n r>0\}$$ of $SL(3,\R)$, 
$G=SL(3,\C)$, $\tau(g)=\bar g$, 
and $\sigma$ is the order $6$ automorphism of
$SL(3,\C)$ defined by 
$$\sigma(g)= C\ (g^t)^{-1} C^{-1}, \ \  {\rm where\ \ } C=\begin{pmatrix} 0& \a^2&0\cr \a&0&0\cr
0&0&1\cr \end{pmatrix}, \
\  {\rm and\ }\   \a=e^{2\pi i\over 3}.$$ The induced automorphism  $\sigma$ on $sl(3,\C)$ is 
$$\sigma(A)= -CA^t C^{-1}.$$
Note that the order of $\sigma$ is $6$,  $\sigma$ is complex linear on $sl(3,\C)$, and
$\sigma\tau=\tau^{-1}\sigma^{-1}$.  Let $\b=e^{2\pi i\over 6}= e^{\pi i\over 3}$.  
A direct computation implies that $Y_j$ lies in the eigenspaces $\cg_j$ of $\b^j$ if 
\begin{align*}
&Y_0=\begin{pmatrix} s&0&0\cr 0&-s&0\cr 0&0&0\cr \end{pmatrix},\ \ 
Y_1=\begin{pmatrix} 0&0&s_1\cr s_2&0&0\cr 0&s_1&0\cr \end{pmatrix},\cr
& Y_2=\begin{pmatrix} 0&0&0\cr 0&0&s\cr -s&0&0 \end{pmatrix}, \ \ Y_3=\begin{pmatrix} s&0&0\cr
0&s&0\cr 0&0&-2s \end{pmatrix},\cr
&Y_4=\begin{pmatrix} 0&0&s\cr 0&0&0\cr 0&-s&0 \end{pmatrix},  \ \  Y_5= \begin{pmatrix}
0&s_1&0\cr 0&0&s_2\cr s_2&0&0
\end{pmatrix}.\cr
\end{align*}
The fixed point set of $\tau$ is $\cu= sl(3, \R)$, and $\cu_j= sl(3,\R)\cap \cg_j$.  

Let 
$$a=\begin{pmatrix} 0&0&1\cr 1&0&0\cr 0&1&0 \end{pmatrix}\ \in \cu_1, \ \
b=\begin{pmatrix} 0&1&0\cr 0&0&1\cr 1&0&0 \end{pmatrix}\
\in \cu_{-1}=\cu_5.$$ Note that $[a,b]=0$.    
The fixed point set $U_0$ of $\sigma$ on $U$  is the abelian group 
 $$U_0=\{\diag(r, r^{-1},1)\n r>0\}.$$  
A smooth map $g:\R^2\to U_0$ is of the form
$$g=\begin{pmatrix} e^{w} &0&0\cr 0& e^{-w}&0\cr 0&0&1 \end{pmatrix}$$
for some smooth function $w$. 
So $g^{-1}g_x= w_x\diag(1, -1,0)$, and 
$$g^{-1}bg =\begin{pmatrix} 0&e^{-2w} & 0\cr 0&0&e^w\cr  e^w&0&0 \end{pmatrix}.$$
Hence the $-1$-flow \eqref{cu} associated to ${SL(3,\R)/U_0}$ is 
$$w_{xt}\diag(1,-1,0)= (e^w-e^{-2w}) \diag(1, -1, 0), $$
i.e., the {\it Tzitzeica equation\/}:
\begin{equation} \label{ck}
w_{xt}= e^w-e^{-2w}.
\end{equation}
The corresponding Lax pair $\o_\l$ \eqref{cv} is
\begin{equation} \label{cj}
\o_\l= \begin{pmatrix} w_x& 0 &\l\cr \l & -w_x & 0\cr 0 & \l &0 \end{pmatrix}\ dx +
\l^{-1}\begin{pmatrix} 0& e^{-2w}& 0\cr 0& 0 & e^w\cr e^w &0&0 \end{pmatrix} \ dt.
\end{equation}
Note that $\o_\l$ satisfies the ${SL(3,\R)\over \R^+}$-reality condition:
\begin{equation}\label{cy}
\overline{\o_{\bar\l}}= \o_\l, \ \  -C\o_\l^t C^{-1} = \o_{\b\l}, \ \ {\rm where\ \ }
\b=e^{2\pi i\over 6}= e^{\pi i\over 3}.
\end{equation}

\beg\label{go}
The $-1$-flow associated to $SL(n,\R)/(\R^+)^{n-1}$.
\eeg

 Let $\cg=sl(n,\C)$, $\tau(\xi)=\bar \xi$, and $\sigma(\xi)= C\xi C^{-1}$, where 
$$C=\diag(1, \a, \cdots, \a^{n-1})$$ and $\a=e^{2\pi i\over n}$.  Then the order of $\sigma$ is $n$
and $\tau\sigma = \sigma^{-1}\tau^{-1}$. The fixed point set $\cu$ of $\tau$ is $sl(n,\R)$.   The
eigenspace
$\cg_j$ of $\sigma$ with eigenvalue $e^{2\pi j i\over n}$ is  spanned by $\{e_{i+j, i}\n i=1, \cdots, n\}$,
where $e_{ij}$ is the elementary matrix and $e_{ij}=e_{i'j'}$ if $i\equiv i'$ and $j\equiv j'$ mod $n$. 
Let $\cu_j=\cg_j\cap \cu$, 
$a=e_{21}+e_{32}+\cdots + e_{n, n-1} + e_{1n}$, and
$b=e_{12}+e_{23}+\cdots + e_{n-1, n}+ e_{n1}$.  Then $a\in \cu_1$, $b\in \cu_{-1}$, and $[a,b]=0$. 
Let 
$g=\diag(e^{u_1},\cdots, e^{u_n})$ with $\sum_i u_i=0$.  So $g^{-1}g_x=  \diag((u_1)_x, \cdots,
(u_n)_x)$.  The
$-1$-flow \eqref{cu} associated to $G, \tau, \sigma$ is  the {\it $2$-dimensional periodic Toda lattice\/}
(\cite{McI94}):
$$(u_i)_{xt}= e^{u_i - u_{i-1}} - e^{u_{i+1}- u_i},  \ \ 1\leq i\leq n,$$
where $u_{n+1}=u_1$ and $u_0=u_n$.

\ms
\subsection{\bf The hyperbolic systems} \label{ki}
\hfil\break

Let $\cu$ be the real form defined by the involution $\tau$ on $\cg$.  
{\it The hyperbolic $U$-system\/} is the following system for $(u_0,u_1, v_0,v_1):\R^2\to
\prod_{i=1}^4\cu$:
\begin{equation} \label{ez}
\bca
(u_1)_t=[u_1,v_0], &\cr (u_0)_t= (v_0)_x + [u_1,v_1]+[u_0,v_0],&\cr
(v_1)_x=-[u_0, v_1].
\eca
\end{equation}
 It has a Lax pair 
\begin{equation}\label{he}
\W_\l = (u_1\l + u_0) dx + (\l^{-1} v_1+v_0) dt,
\end{equation}
which satisfies the $U$-reality condition \eqref{fp}.

Let $\sigma$ be an order $k$ automorphism of $\cg$ such that $\tau\sigma=\sigma^{-1}\tau$, and
$\cu=\cu_0+\cdots +\cu_{k-1}$, where $\cu_j$ is the intersection of $\cu$ and the $e^{2\pi i j\over
k}$-eigenspace of $\sigma$.   {\it The hyperbolic
$U/U_0$-system\/} is the restriction of the hyperbolic $U$-system \eqref{ez} to $(u_0,u_1,v_0,v_1)\in
\cu_0\times \cu_1\times \cu_0\times \cu_{-1}$.   The corresponding Lax pair \eqref{he} satisfies the
$U/U_0$-reality condition \eqref{cb}.

\bs

\section{\bf  Geometries associated to soliton equations} \label{kj}

We give geometric interpretations of certain soliton equations.  For example, solutions of the
$O(2n)/O(n)\times O(n)$-system give rise  to orthogonal coordinates of $\R^n$ and flat submanifolds
in $\R^{2n}$, solutions of the
$U(n)/O(n)$-system give  Egoroff flat metrics and flat Lagrangian
submanifolds of $\R^{2n}=\C^n$, a subclass of solutions of the $GL(2,\H)/(SU(2)\times
\R^+)^2$-system give rise to Bonnet pairs in $\R^3$, and solutions of the $-1$-flow associated to $SL(3,\R)/\R^+$ 
(the Tzitzeica equation) is the Gauss-Codazzi equation for affine spheres in the affine $3$-space.  If $U/U_0$ is a rank
$n$ symmetric space, then we can associate to each solution of the
$U/U_0$-system a flat $n$-submanifold in $U/U_0$ and a flat $n$-submanifold in the tangent space of
$U/U_0$.   We also give a brief review of the relation between harmonic maps from $\R^{1,1}$ to $U$
and solutions of the hyperbolic $U$-system.  

\ms
\subsection{\bf The method of moving frames}  \label{kk}
\hfil\break

Let $(N,g)$ be an $(n+k)$-dimensional Riemannian manifold, $\K$ the Levi-Civita connection of $g$, 
and $X:M^n\to N$ an immersion.  We set up some notation next. The first fundamental form $\I$ is the
induced metric.  Let
$\xi$ be a normal vector field on $M$, $v$ a tangent vector field, $(\K_v\xi)^t$ and $(\K_v
\xi)^\nu$ the tangential and normal components of $\K_v\xi$ respectively. The induced normal
connection on the normal bundle $\nu(M)$ is defined by
$$\K^\nu_v\xi =(\K_v\xi)^\nu.$$  The second fundamental form $\II$ is  a smooth section of
$S^2(T*M)\otimes \nu(M)$ defined by
$$\II_\xi (v_1, v_2)= -g(\K_{v_1}\xi, v_2).$$
Next we express $\I, \II, \K^\nu$ using a moving frame.   Let
$(e_1, \cdots, e_{n+k})$ be a local orthonormal frame on $M$ such that $e_1, \cdots, e_n$ are
tangent to
$M$.  We use the following index convention:
 $$1\leq A,B,C\leq n+k, \ \ 1\leq i, j,k\leq n, \ \  n+1\leq \a,\b,\g\leq n+k.$$
Let $w_A$ denote the dual coframe of $e_A$, and write 
$$\K e_A= \sum_B w_{BA}e_B, \ \  w_{AB}+w_{BA}=0.$$
Then we have 
\begin{align*}
&dX=\sum_i w_i e_i,\cr
&dw_A=-\sum_B w_{AB}\wedge w_B,\ \  \ w_{AB}+w_{BA}=0,\cr
&dw_{AB}=-\sum_C w_{AC}\wedge w_{CB}+\sum_{CD} \ti R_{ABCD}\ w_C\wedge w_D,\cr
\end{align*}
where $\ti R_{ABCD}$ is the coefficients of the Riemann tensor of $g$.  
The first fundamental of $M$ is
$$\I= \sum_i w_i^2.$$
Let 
$$w_{i\a}=\sum_j h^\a_{ij} w_j.$$
Since $w_\a=0$, $-\sum_i w_{\a i}\wedge w_i =0$.  This implies that $h_{ij}^\a = h_{ji}^\a$. 
The second fundamental form and the normal connection are 
\begin{align*}
&\II=\sum_{i,\a}w_i w_{i\a}e_\a =\sum_{i,j,\a} h_{ij}^\a w_iw_je_\a,\cr
&\K^\nu(e_\a)=\sum_i w_{\b\a} e_\b\cr
\end{align*}
respectively.  The normal curvature is the curvature of $\K^\nu$, i.e., 
$$\W^\nu_{\a\b}= dw_{\a\b}+\sum_\g w_{\a\g}\wedge w_{\g\b}.$$
The normal bundle is {\it flat\/} if the normal connection is flat, i.e., $\W^\nu_{\a\b}=0$ for all $\a,\b$. 
The Levi-Civita connection of $\I$ is $(w_{ij})$, and the curvature  is 
$$\sum_{kl} R_{ijkl}w_k\wedge w_l = dw_{ij}+\sum_k w_{ik}\wedge w_{kj} = -\sum_\a
w_{i\a}\wedge w_{\a j} + \ti R_{ijkl}.$$ 
  Note that given $\I, \II, \K^\perp$ is the same as
given $w_i, w_{i\a}, w_{\a\b}$. Moreover,  the Levi-Civita
connection of $\I$ can be obtained by solving
$$dw_i=-\sum_j w_{ij}\wedge w_j, \ \ \  w_{ij}+w_{ji}=0.$$  The Gauss-Codazzi equation is
\begin{equation}\label{fy}
\begin{cases}
dw_{ij}+\sum_k w_{ik}\wedge w_{kj}= -\sum_\a w_{i\a}\wedge w_{\a j} + \ti
R_{ijkl}w_k\wedge w_l,& \\ 
dw_{i\a}=-\sum_A w_{iA} \wedge w_{A\a}, & \\ 
dw_{\a\b}+ \sum_\g w_{\a\g}\wedge w_{\g\b} = -\sum_i w_{\a i}\wedge w_{i\b} + \sum_{ij}\ti
R_{\a\b i j} w_i\wedge w_j. &\\
\end{cases}
\end{equation}
The Fundamental theorem for submanifolds states that $\I, \II$ and
$\K^\nu$ together with the Gauss-Codazzi equation \eqref{fy} determine the submanifold $M$ up to
isometries of $N$. 

The mean curvature vector field is defined as the trace of $\II$ with respect to $\I$, i.e.,
$$H=\tr_{\I} \II = \sum_{i\a} h^\a_{ii} e_\a.$$
The normal bundle $\nu(M)$ is said to be {\it non-degenerate\/} if the dimension of the space
of all shape operators of $M$, $\{A_v\n v\in \nu(M)_p\}$ is equal to dim$(M)$ for all $p\in
M$. 

If  $X:M\to N$ is a submanifold of a space form $N^{n+k}$, then the frame $F=(X, e_1, \cdots,
e_{n+k})$ given above is a lift of $X$ to $\Iso(N)$ and the Gauss-Codazzi equation for $M$ is exactly
the flatness of $F^{-1}dF$.  When $M$ satisfies certain geometric conditions, we often can find 
special coordinates and frames $F$ on $M$ so that $F^{-1}dF$ takes a special simple form.  If
moreover, such submanifolds admit a natural holomorphic deformation, then the Gauss-Codazzi
equation for $M$ is likely to be an integrable system.  On the other hand, if the Lax $n$-tuple of the
$U/U_0$-system can be interpreted as the connection $1$-form of a submanifold, then we can read its
geometry from the Lax $n$-tuple.  This gives a natural method to find interesting submanifolds whose equations
are integrable.  We have had some success when $U$ is an orthogonal group or a unitary group.  But very little is
known for other simple Lie group $U$.  

\ms
\subsection{\bf  Orthogonal coordinate systems and the $\onn$-system}  \label{kl}
\hfil\break

A local coordinate system $(x_1,
\cdots, x_n)$ of $\R^n$ is called an {\it orthogonal coordinate system\/}  if the flat metric
written in this coordinate system is diagonal, i.e., of the form
$\sum_{i=1}^n g_{ii}(x) dx_i^2$.  The theory of orthogonal coordinate systems of $\R^n$ was studied
extensively by classical differential geometers (cf. Darboux \cite{Dar10}).  

An elementary computation gives:

\bprop\label{bq}
 The Levi-Civita connection 1-form $(w_{ij})$ of the metric
$ds^2=
\sum_{i=1}^n b_i^2 dx_i^2$  is 
$$w_{ij}= {(b_i)_{x_j}\over b_j} dx_i - {(b_j)_{x_i}\over b_i} dx_j.$$
\eprop

So the Levi-Civita connection $1$-form $w$ of  $ds^2=\sum_{i=1}^n b_i^2
dx_i^2$ written in matrix form is
$$w=(w_{ij}) = \d F - F^t \d, \ \ \ {\rm where\ \ } f_{ij}= 
\bca
{(b_i)_{x_j}\over b_j},& \text{if
$i\not=j$}\cr f_{ii}=0, & \text{if $1\leq i\leq n$.}\cr
\eca
$$

\ms
Let $gl_\ast(n)$ denote the space of $\xi=(\xi_{ij})\in gl(n,\R)$ such that $\xi_{ii}=0$ for $1\leq
i\leq n$.  
Recall that  $F=(f_{ij}):\R^n\to gl_\ast(n)$ is a solution of the ${O(2n)\over O(n)\times O(n)}$-system
\eqref{ag} if and only if both
$\d F - F^t\d$ and $\d F^t- F\d$ are flat connection $1$-forms.  Note that both connections have the
same form as the Levi-Civita connection of an orthogonal metric.  In the rest of the section, we try to
answer the following question: Are there orthogonal coordinate systems of $\R^n$ whose Levi-Civita
connections are $\d F-F^t\d$ and
$\d F^t- F\d$?  

Given $F:\R^n\to gl_\ast(n)$, there is a diagonal metric whose Levi-Civita connection
$1$-form is $$w=\d F- F^t\d$$ if and only if there exist positive functions
$b_1, \cdots, b_n$ so that 
\begin{equation}\label{al}
(b_i)_{x_j}= f_{ij} b_j, \quad i\not=j.
\end{equation}
However, if system \eqref{al} is solvable for $b_1, \cdots, b_n$, then the mixed derivatives must be
equal.  This implies that
\begin{equation}\label{ak}
(f_{ij})_{x_k}= f_{ik}f_{kj},  \quad \   i, j, k \ {\rm distinct.\/}
\end{equation}
It is a classical result that this condition is also sufficient for \eqref{al} to be solvable:

\bthm\label{ap}
 Given a smooth function  $F=(f_{ij}):\R^n\to gl_\ast(n)$,
system \eqref{al} is solvable for
$(b_1, \cdots, b_n)$ if and only if $F$ satisfies \eqref{ak}.  Moreover, given $n$ smooth one variable
functions $b_1^0, \cdots, b_n^0$, there exists a unique local solution $(b_1, \cdots, b_n)$ of 
\eqref{al} such that $b_i(0, \cdots, x_i, 0, \cdots, 0) = b_i^0(x_i)$. 
\ethm

\bcor\label{br}
 The space of local $n$-dimensional orthogonal metrics that have the same
Levi-Civita connection $1$-form is parametrized by $n$ smooth positive functions of one variable.  
\ecor
If $F$ is a solution of the $\onn$-system \eqref{ag}, then $F$ is a solution of \eqref{ak}.  So by Theorem
\ref{ap} we can construct orthogonal coordinates of
$\R^n$, whose Levi-Civita connections are  $\d F-F^t\d$ and $\d F^t-F\d$.  Therefore we have 
   
\bprop\label{bc}
  Let $F=(f_{ij})$ be  a solution of the $\onn$-system \eqref{ag},  $\tau_1= \d
F^t - F \d$, $\tau_2= \d F- F^t\d$, and $a_1^0,
\cdots, a_n^0, b_1^0, \cdots, b_n^0$  smooth positive functions of
one variable.  Then there exist unique flat local orthogonal metrics $g_1= \sum_{i=1}^n a_i^2(x)
dx_i^2$ and  $g_2= \sum_{i=1}^n b_i(x)^2 dx_i^2$ such that
\begin{enumerate}
\item[(i)] $a_i(0, \cdots, x_i, 0, \cdots)= a_i^0(x_i)$ and 
$b_i(0, \cdots, x_i, 0, \cdots)= b_i^0(x_i)$,
\item[(ii)] the Levi-Civita connection $1$-form for $g_1$ and $g_2$ are 
$\tau_1$ and $\tau_2$ respectively, 
\item[(iii)] there exist $O(n)$-valued maps $A=(\xi_1, \cdots, \xi_n)$ and $B=(\eta_1, \cdots,
\eta_n)$ such that $A^{-1}dA=\tau_1$ and $B^{-1} dB=\tau_2$,
\item[(iv)]there exist $\phi$ and
$\psi$ defined on a neighborhood of the origin in $\R^n$  such that 
$$d\phi= \sum_i a_i
\eta_i dx_i, \ \ d\psi=\sum_i b_i \xi_idx_i,$$ 
\item[(v)] $\phi$ and $\psi$ are local orthogonal
coordinates on $\R^n$ with  Levi-Civita connection $\tau_1$ and $\tau_2$ respectively.
\end{enumerate}
\eprop

The next theorem states that a subclass of  orthogonal coordinate systems of $\R^n$
can be obtained using trivialization of $\tau_1$ and $\tau_2$. 

\bthm\label{az}
  Let $F=(f_{ij})$ be a solution of \eqref{ag}, and  $A=(a_{ij})$, $B=(b_{ij})$
smooth $O(n)$-valued maps defined on an simply connected domain $\co$ of $\R^n$ satisfying
\begin{equation}\label{ar}
A^{-1}d A= \d F^t- F\d,\ \  B^{-1}dB= \d F- F^t \d.
\end{equation}
If $a_{mj}, b_{mj}$ never vanishes on $\co$ for all $1\leq  j\leq n$, then:
\begin{enumerate}
\item[(i)] $ds_m^2 = a_{m1}^2 dx_1^2 + \cdots + a_{mn}^2 dx_n^2$ is a flat metric with $\d
F- F^t \d$ as its Levi-Civita connection,
\item[(ii)] $d\ti s_m^2 = b_{m1}^2 dx_1^2 + \cdots + b_{mn}^2 dx_n^2$ is a flat metric with
$\d F^t- F\d$ as its Levi-Civita connection,
\item[(iii)]  there exists a smooth map $X:\co\to gl(n,\R)$ such that 
$$dX=B\d A^t,$$
\item[(iv)] the $m$-th column $X_m$ and the $m$-th row $Y_m$ of $X$ are local orthogonal
coordinates for $\R^n$ such that the standard metric on $\R^n$ written in these coordinates are
$ds_m^2$  and $d\ti s_m^2$ respectively. 
\end{enumerate} 
\ethm

\proof 
Let $\xi_i$ denote the $i$-th column of $A$. W claim that 
\begin{equation}\label{ds}
(\xi_j)_{x_k}= f_{jk} \xi_k, \  \  j\not=k.
\end{equation}
Note that \eqref{ar} gives
\begin{equation}\label{ah}
\xi_i \cdot d\xi_j = f_{ji} dx_i - f_{ij} dx_j, \ \  i\not=j.
\end{equation}
This implies
$$(\xi_j)_{x_k} \cdot \xi_i =0, \ \ {\rm if\ } i, j, k \ \ {\rm distinct\/}. $$
Since $\xi_j\cdot \xi_j=1$,  $(\xi_j)_{x_k}\cdot \xi_j=0$. By \eqref{ah}, $\xi_k\cdot
(\xi_j)_{x_k}= f_{jk}$.  This proves \eqref{ds}.  
Equate each coordinate of \eqref{ds} to get
$$(a_{mj})_{x_k}= f_{jk} a_{mk}, \ \ 1\leq m\leq n,  \ j\not=k.$$
By Proposition \ref{bq}, the Levi-Civita connection of $ds_m^2$ is 
$${(a_{mj})_{x_k}\over a_{mk}} dx_j - {(a_{mk})_{x_j}\over a_{mj}} dx_k = f_{jk} dx_j - f_{kj}
dx_k,$$ i.e., the Levi-Civita connection $1$-form for $ds^2_m$ is $\d F- F^t \d$.  
 This proves (i).  Similar argument gives (ii). 

Since $F$ is a solution of \eqref{ag}, 
$$\o_\l= \begin{pmatrix} \d F^t - F\d& -\d \l\cr \d \l & \d F - F^t \d \end{pmatrix}$$ is flat.  Let
$h=\begin{pmatrix} A&0\cr 0& B\cr \end{pmatrix}$.  Then $h^{-1}dh = \o_0=\begin{pmatrix} \d
F^t-F\d&0\cr 0&\d F- F^t
\d \end{pmatrix}$.   The gauge transformation of $\o_\l$ by
$h$ is 
$$\Theta_\l=h\o_\l h^{-1}- dh h^{-1}=
\bpm
0&-\l A\d B^t\cr \l B \d A^t & 0\cr
\epm.$$
Since $\o_\l$ is flat for all $\l\in \C$, $\Theta_\l= h\ast \o_\l$ is flat for all $\l$, i.e.,
$d\Theta_\l= -\Theta_\l\wedge\Theta_\l$.   This gives
$$\l \bpm 
0& -d\zeta^t\cr d\zeta &0\cr
\epm
 +\l^2 \bpm
-\zeta^t\wedge \zeta &0\cr 0 &
-\zeta  \wedge \zeta^t\cr
\epm = 0,$$
where $\zeta= B\d A^t$.  
Compare coefficients of $\l$ to get $d\zeta=0$.  Since $\co$ is simply connected, there exists
$X$ such that 
\begin{equation}\label{at}
dX= B\d A^t.
\end{equation}
  This proves (iii).  

Equate the $m$-th column and row of \eqref{at} to get
\begin{equation}\label{as}
dX_m= B(a_{m1}dx_1, \cdots, a_{mn}dx_n)^t.
\end{equation}
Recall that $A=(\xi_1, \cdots, \xi_n)$.  Write equation \eqref{as}
using columns of $A$ and $B$ to get
$$dX_m= \sum_{i=1}^n a_{mi}dx_i \eta_i.$$
Let $\eta_i$ denote the $i$-th row of $X$, and $Y=X^t$.  Then $dY= A\d B^t$ and 
$dY_m=\sum_{i=1}^n b_{mi} dx_i \xi_i$.  
This proves (iv).  \qed

\ms
\subsection{\bf  Flat submanifolds and the $\onn$-system} \label{km}
\hfil\break

The $\onn$-system can also be viewed as the Gauss-Codazzi equations for flat $n$-dimensional
submanifolds in $\R^{2n}$ with flat and non-degenerate normal bundles. 
In fact, there is an isomorphism from the space of local
$n$-dimensional flat submanifolds in
$\R^{2n}$ with flat and non-degenerate normal bundle modulo rigid motions to the space of  
$(F,c_1,\cdots, c_n)$, where $F$ is a local solution of the $\onn$-system \eqref{ag} and $c_1, \cdots,
c_n$ are positive functions of one variable.  We state this more precisely in the following two known
theorems (cf. \cite{Ter97}):

\bthm\label{ao}
 Let $M^n$ be a
$n$-dimensional flat submanifold of
$\R^{2n}$ with flat and non-degenerate normal bundle.  Then there exist
local coordinates $x_1, \cdots$, $x_n$, parallel normal frame $e_{n+1}, \cdots, e_{2n}$, an
$O(n)$-valued map $A=(a_{ij})$, and a map
$b=(b_1, \cdots, b_n)$ such that the fundamental forms of $M$ are
\begin{equation}\label{av}
\bca
\I= \sum_{i=1}^n b_i^2 dx_i^2, &\cr
\II= \sum_{i, j=1}^n b_ia_{ji}dx_i e_{n+j}.&\cr
\eca
\end{equation}
Moreover, let $f_{ij}=(b_i)_{x_j}/b_j$ for $1\leq i\not=j\leq n$, $f_{ii}=0$ for $1\leq i\leq n$,
and $F=(f_{ij})$.   Then $F$ is a solution of the $\onn$-system \eqref{ag}.  
\ethm

\bthm\label{fc}
 Let $F$ be a solution of the ${O(2n)\over O(n)\times O(n)}$-system \eqref{ag}, and
$b_1^{0}, \cdots$, $b_n^{0}$ be $n$ smooth positive functions  of one variable.  Then  there exist an
open subset $\co$ of the origin in $\R^n$, smooth maps $A:\co\to O(n)$ and 
\begin{equation}\label{aq}
\phi=\bpm
g & X\cr 0 &1\cr
\epm: \ \co\to GL(2n+1, \R)
\end{equation}
 with $g:\co\to O(2n)$,
$X:\co\to \R^{2n}$, and $b_1, \cdots, b_n:\co\to \R$ such that 
\begin{align*}
&A^{-1}dA =\d F^t -F\d,\cr
&\phi^{-1}d\phi=\bpm
0& -A\d & 0\cr \d A^t & \d F- F^t \d
&\varpi \cr 0 & 0 & 0\cr
\epm,\cr  &b_i(0,\cdots, x_i,0, \cdots)= b_i^0(x_i), \quad 1\leq i\leq n,
\end{align*}
where $\varpi=(b_1dx_1, \cdots, b_ndx_n)^t$.
Moreover,
\begin{enumerate} 
\item[(i)] $X$ is an immersion of a flat
$n$-dimensional submanifolds of $\R^{2n}$ with flat and non-degenerate normal bundle, 
\item[(ii)] $g=(e_{n+1}, \cdots, e_{2n}, e_1, \cdots, e_n)$ is a local orthonormal frame for $X$ such
that $e_{n+1}, \cdots, e_{2n}$ are parallel normal field,
\item[(iii)] $b_i(0,\cdots, x_i, 0, \cdots, 0)= b_i^0(x_i)$ for $1\leq i\leq n$,
\item[(iv)] the fundamental forms of the immersion $X$ are given as in \eqref{av},
\item[(v)] the Levi-Civita connection for the induced metric is $\d F- F^t\d$.   
\end{enumerate}
\ethm

 \'E Cartan proved that a flat $n$-dimensional submanifold can not be
locally isometrically immersed in $S^{n+k}$ if $k< n-1$, but can be locally isometrically immersed into
$S^{2n-1}$.   Moreover,  the normal bundle of a flat $n$-dimensional submanifold of
$S^{2n-1}$ is flat, and is non-degenerate viewed as a submanifold of $\R^{2n}$.  By Theorem \ref{ao},
flat $n$-dimensional submanifolds in $S^{2n-1}$ give rise to solutions of the $\onn$-system
\eqref{ag}.  This gives the following theorem of Tenenblat (\cite{Ten85}):

\bthm\label{au}
 (\cite{Ten85}).  Let $X:M^n\to S^{2n-1}$ be an immersion of a flat submanifold.  Then
there exist local coordinates $x_1, \cdots, x_n$, parallel normal frame $e_{n+2}, \cdots,
e_{2n}$, and a smooth
$O(n)$-valued map $A=(a_{ij})$ such that
$$\I= \sum_{i=1}^n a_{1i}^2 dx_i^2, \ \ \II=\sum_{i=1, j=2}^n a_{1i} a_{ji} dx_i e_{n+j}.$$
Set $f_{ij}= (a_{1i})_{x_j}/a_{1j}$ for $i\not=j$, $f_{ii}=0$, and $F=(f_{ij})$.  Let
$e_i=X_{x_i}/a_{1i}$ for $1\leq i\leq n$, and $g=(X, e_{n+2}, \cdots, e_{2n}, e_1, \cdots, e_n)$.  
Then $F$ is a solution of the ${O(2n)\over O(n)\times O(n)}$-system \eqref{ag} and 
\begin{equation}\label{aha}
g^{-1}dg=\bpm
0& -A\d\cr \d A^t & -F^t\d + \d F\cr
\epm.
\end{equation} 
Conversely, let
$F=(f_{ij})$ be a solution of \eqref{ag}, $A=(a_{ij})$ an $O(n)$-valued map such that $A^{-1}dA= \d
F^t- F\d$ and  $g\in O(2n)$ a solution of \eqref{aha}.   If 
$a_{ij}>0$ for all $1\leq  j\leq n$ on an open subset $\co$ of $\R^n$, then
the $i$-th column of $g$ is an immersion of flat submanifold in
$S^{2n-1}$ and the corresponding  solution of \eqref{ag} is $F$. 
\ethm

\bcor\label{bb}
 Let $M^n\subset S^{2n-1}$ be a flat submanifold, $\xi$ a parallel
normal field such that the shape operator $A_\xi$ is non-degenerate, then $\xi$ is
an immersion of a flat $n$-submanifold in $S^{2n-1}$.   Moreover, the solution of \eqref{ag}
corresponding to
$\xi$ is the same as the one corresponding to $M$.  
\ecor

\ms
\subsection{\bf  Egoroff metrics and the $\un$-system}  \label{kn}
\hfil\break

The $\un$-system is the restriction of the $\onn$-system to the subspace of symmetric real $n\times
n$ matrices $F$.  We have seen that each solution $F$ of the $\onn$-system gives rises to a
$o(n)$-connection of some flat diagonal metric.  In this section, we show that such diagonal metric takes
a special form:

\bprop \label{bh}
 Let $ds^2=\sum_{i=1}^n b_i^2 dx_i^2$ be a metric,  $f_{ij}=
(b_i)_{x_j}/b_j$ for $1\leq i\not= j\leq n$, and $F=(f_{ij})$.  Then $F=F^t$ if and only if there exists a
function $\phi$ such that $b_i^2= \phi_{x_i}$ for all $1\leq i\leq n$. 
\eprop

\proof Since $f_{ij} = {(b_i)_{x_j}\over b_j}$, $F^t=F$ if and only if
$${(b_i)_{x_j}\over b_j} = {(b_j)_{x_i}\over b_i}, \quad i\not=j.$$
This is equivalent to $(b_i^2)_{x_j}= (b_j^2)_{x_i}$ for all $i\not=j$.  \qed

\bdefn\label{ff}
An {\it Egoroff\/} metric is a flat metric of the form $$\sum_{i=1}^n \phi_{x_i}
dx_i^2$$ for some smooth function $\phi$.  
\edefn

\ms
It follows from Proposition \ref{bc}, Theorem \ref{az}, and Proposition \ref{bh} that:

\bthm \label{cc}
  Let $F$ be a solution of the $\un$-system \eqref{an}, and $a_1$, $\cdots$,
$a_n$ smooth positive functions of one variable.  Then there exists a smooth function $\phi$ defined
on a simply
 connected open subset $\co$ of $\R^n$ such that $$\phi_{x_i}(0, \cdots, 0, x_i, 0, \cdots, 0)=
a_i^2(x_i)$$ for
$1\leq i\leq n$ and the Levi-Civita connection for $\sum_{i=1}^n \phi_{x_i} dx_i^2$ is $[\d, F]$. 
Moreover,  let $A=(a_{ij})$ be an $O(n)$-valued map such that $A^{-1}dA= [\d, F]$.  Then:
\begin{enumerate}
\item[(i)] $ds_m^2=\sum_{i=1}^n a_{mi}^2 dx_i^2$ is an Egoroff metric with $[\d,F]$ as its
Levi-Civita connection,
\item[(ii)] there exists a smooth map $X$ from $\co$ to the space of symmetric matrices such that
$dX= A\d A^t$, 
\item[(iii)] the $m$-th column $X_m$ of $X$ is a local orthogonal coordinate system for $\R^n$ and
the flat metric of $\R^n$ written in this coordinate system is $ds_m^2$ as in (i). 
\end{enumerate}
\ethm

\ms
\subsection{\bf   Flat Lagrangian submanifolds and the $\un$-system}  \label{ko}
\hfil\break

In this section, we explain the relation between solutions of the $\un$-system and the
Gauss-Codazzi equations for flat, Lagrangian submanifolds of $\R^{2n}$. If these submanifolds also
lie in $S^{2n-1}$, then they are invariant under the $S^1$-action of the Hope fibration. 
Hence the projection of these submanifolds are flat Lagrangian submanifolds of $\C P^{n-1}$.  For
more detail of the geometry of flat Lagrangian submanifolds of $\C P^{n-1}$ see \cite{DajToj95b}. 

 Let $\li \ , \ \ri$ and $w$ be the standard inner product
and symplectic form on $\C^n =\R^{2n}$ respectively, i.e., 
$$\li X, Y\ri = {\rm Re\/} (\bar X^t Y), \quad w(X,Y)= \Im(\bar X^tY), \ \ X, Y\in \C^n.$$
Write $Z\in \C^n$ as $Z=X+iY\in \R^n + i\R^n$, and $A\in gl(n,\C)$ as $A=B+iC$ with $B,C\in
gl(n,\R)$.  Then
$A\in gl(n,\C)$ is identified as $\bpm
B& -C\cr C& B\cr
\epm$ in $gl(2n,\R)$.  This
identifies $u(n)$ as the following subalgebra of $o(2n)$:
$$u(n)= \left\{\bpm
B& -C\cr C& B\cr
\epm
\in o(2n)\ \bigg| \ B\in o(n), C\in
gl(n,\R) \ {\rm symmetric\/} \right\}.$$
The standard complex structure on $\R^{2n}$ is 
$$J\bpm
X\cr Y\cr
\epm = \bpm
-Y\cr X\cr
\epm.$$

\bdefn \label{bd}
An $n$-dimensional submanifold $M$ of $\C^n=\R^{2n}$ is {\it
Lagrangian\/} if $w(v_1, v_2)=0$ for all $v_1, v_2\in TM$, or equivalently, $J(TM)= \nu(M)$.
\edefn

\ss

The Proposition below follows from the definition of Lagrangian submanifold:

\bprop \label{be}
  Let $X:M^n\to \R^{2n}$ be a Lagrangian submanifold, and
$(e_1,$ $ \cdots, e_n)$ a local orthonormal tangent frame.  Then
$(Je_1, \cdots, J e_n)$ is a orthonormal normal frame.  Moreover, if $g=(Je_1, \cdots, Je_n, e_1,
\cdots, e_n)$, then $g^{-1}dg$ is a $u(n)$-valued $1$-form, i.e., it is of the form
$\bpm
\xi& -\eta\cr \eta&\xi\cr
\epm$,
where $\xi$ is an $o(n)$-valued 1-form and $\eta$ is $1$-form with value in the space of symmetric
matrices.   Conversely, if $M^n$ has a local orthonormal frame $g=(e_{n+1}, \cdots, e_{2n}, e_1,
\cdots, e_n)$ such that $e_1, \cdots, e_n$ are tangent to $M$ and $g^{-1}dg$ is $u(n)$-valued
$1$-form, then $M$ is Lagrangian.  
\eprop

\bprop \label{bf}
  Let $F=(f_{ij})$ be the solution of the ${O(2n)\over O(n)\times
O(n)}$-system \eqref{ag} corresponding to the flat
$n$-submanifold $M$ of  $\R^{2n}$ with flat and non-degenerate normal bundle as in Theorem
\ref{ao}.  Then the following statements are equivalent:
\begin{enumerate}
\item[(i)] $F$ is a solution of the ${U(n)\over O(n)}$-system \eqref{an}, 
\item[(ii)] $F=F^t$,
\item[(iii)] $M$ is Lagrangian.  
\end{enumerate}
\eprop

\proof It is obvious that (i) and (ii) are equivalent.  

 Let $x_1, \cdots, x_n$, $b_1, \cdots, b_n$, $e_{n+1}, \cdots, e_{2n}$,  and
$A=(a_{ij})$ be as in Theorem
\ref{ao}. Let $e_i={X_{x_i}\over b_i}$, and $g=(e_{n+1}, \cdots, e_{2n}, e_1, \cdots, e_n)$. 
Then 
$$g^{-1}dg= \bpm
0& -A\d\cr \d A^t& [\d, F]\cr
\epm.$$

To prove (ii) implies (iii), let
$$\phi=(\ti e_{n+1}, \cdots, \ti
e_{2n}, e_1,\cdots, e_n)=g\bpm
A&0\cr 0&\I\cr
\epm.$$ 
Then $\ti e_{n+1}, \cdots, \ti e_{2n}$ are normal to $X$, and 
\begin{equation}\label{bs}
\phi^{-1}d\phi= \bpm
[\d, F]& -\d\cr \d & [\d, F]\cr
\epm.
\end{equation}
Proposition \ref{be} implies that  $M$ is Lagrangian.  
 
(iii) implies (ii):
Since $M$ is Lagrangian and $(e_{n+1}, \cdots, e_{2n})$ is an orthonormal normal frame,
$(Je_{n+1},
\cdots Je_{2n})$ is an orthonormal tangent frame for
$M$.  So there exists an $O(n)$-valued map $h$ such that 
$$\ti g:= (e_{n+1}, \cdots, e_{2n}, Je_{n+1}, \cdots, Je_{2n})= g\bpm
I& 0\cr 0& h^{-1}\cr
\epm.$$
Then
$$\ti g^{-1} d\ti g= \bpm
0& -A\d h^{-1}\cr h\d A^t& h(\d F- F^t \d) h^{-1} - dh h^{-1}\cr
\epm.$$
But $\ti g^{-1}d \ti g$ is $u(n)$-valued.   Hence

\begin{equation}\label{aw}
\bca
h(\d F- F^t \d) h^{-1} - dh h^{-1}=0, \cr
 A\d h^t= h\d A^t\cr
\eca
\end{equation}
The second equation of \eqref{aw} gives $a_{ik}h_{jk}= h_{ik}a_{jk}$ for all $i, j, k$.  This
implies that the $i$-th rows ($i$-th columns resp.) of $A$ and $h$ are proportional.  Since both $A$
and $h$ are in $O(n)$, we have
$h=A$.  So  $h^{-1} dh = A^{-1}dA=  \d F- F^t \d$. 
But $A^{-1}dA= \d F^t- F\d$.  Hence 
$$\d F^t- F\d= \d F- F^t \d.$$  Equate the $ij$-th entry of the above equation to get
$f_{ji} dx_i - f_{ij} dx_j = f_{ij} dx_i - f_{ji} dx_j$. So $F$ is symmetric.  \qed 

\ms

As a consequence of Proposition \ref{bf}, Theorem \ref{au}, and Corollary \ref{bb}, we have

\bcor \label{bg}
  Let $F$ be a solution of the ${U(n)\over O(n)}$-system \eqref{an},
$A=(a_{ij})$ an $O(n)$-valued map satisfying $A^{-1}dA= [\d, F]$,  and
$\ti g$ an $U(n)$-valued map satisfying 
\begin{equation*}
\ti g^{-1}d\ti g =\bpm
0& -A\d A^t\cr A\d A^t & 0\cr\epm.
\end{equation*}
(Here $U(n)$ is embedded as a subgroup of $O(2n)$).  Let $e_{m+n}$ denote the $m$-th column of
$\ti g$ for $1\leq m\leq n$.   If $a_{m1}, \cdots, a_{mn}$ never vanishes in an open subset $\co$ of
$\R^n$, then $e_{n+m}:\co\to S^{2n-1}$ is an $n$-dimensional immersed flat submanifold of
$S^{2n-1}$ that is Lagrangian in $\R^{2n}$.  Conversely, if $M^n$ is a flat submanifold of $S^{2n-1}$
that is Lagrangian in $\R^{2n}$, then $F$ defined in Theorem \ref{au} is a solution of the
$\un$-system \eqref{an}.
\ecor

\bprop\label{bj}
Let $F$ be a solution of the ${U(n)\over O(n)}$-system \eqref{an},
$M^n$ a flat submanifold of $S^{2n-1}$ corresponding to $F$ as in Corollary \ref{bg}, and
$\pi:S^{2n-1}\to \C P^{n-1}$  the Hopf fibration.  Then $M=\pi^{-1}(\pi (M))$ and
$\pi(M)$ is a flat Lagrangian submanifold of
$\C P^{n-1}$. 
\eprop

\proof Let $S^1$ acts on $\R^{2n}=\C^n$ by 
$$e^{is}\cdot (z_1, \cdots, z_n) = (e^{is} z_1, e^{is}z_2, \cdots, e^{is}z_n).$$  This action leaves
$S^{2n-1}$ invariant, the orbit space $S^{2n-1}/S^1$ is $\C P^{n-1}$, and the projection
$\pi:S^{2n-1}\to \C P^{n-1}$ is the Hopf fibration.  

 It  suffices to show that $M$ is invariant under the
$S^1$-action on $S^{2n-1}$. 
Let $X$ be the immersion,  $(x_1, \cdots, x_n)$,  $g= (X, e_{n+2}, \cdots, e_{2n},
e_1, \cdots, e_n)$, and  $A=(a_{ij})$ as in Theorem \ref{au}.
   First we change coordinates $x_1, \cdots, x_n$  to $t_1, \cdots, t_n$ such that
$$\bca
x_1= t_1-t_2- \cdots - t_n,&\cr
 x_j= t_j + t_1, &\text{ $ 2\leq	 j\leq n$}.\cr
\eca$$
Then ${\p\over \p t_1}= {\p\over \p x_1} + \cdots {\p\over \p x_n}$. 
Since 
$$A^{-1}{\p A\over \p t_1}= [\d, F] \left({\p\over \p t_1}\right)=[\I_n, F]=0,$$ we have ${\p
A\over
\p t_1}=0$.  Here $\I_n$ is the identity $n\times n$ matrix.  
  Let $\ti g=g\bpm
\I_n&0\cr 0& A^t\cr
\epm$.  Since $A^{-1}dA=[\d, F]$, 
$$\ti g^{-1} d\ti g= \bpm
0& - A\d A^t\cr A\d A^t & 0\cr
\epm.$$
So we have
\begin{equation} \label{bi}
\ti g^{-1}{\p \ti g\over \p t_1}=\bpm 
0&-I\cr I&0\cr
\epm
\end{equation}
This implies that 
$$\ti g(t_1, \cdots, t_n)= e^{it_1}\ti g(1, t_2, \cdots, t_n).$$
But the first column of $g$ and of $\ti g$ are the immersion $X$.   So $X$ is invariant
under the $S^1$-action on $S^{2n-1}$.\qed

  It follows from elementary submanifold theory that $\ti M$ is a flat Lagrangian
submanifold  of $\C P^{n-1}$ if and only if $\pi^{-1}(\ti M)$ is a flat submanifold of $S^{2n-1}$ that is
Lagrangian in $\R^{2n}$.  Hence the Gauss-Codazzi equations for flat, Lagrangian
submanifolds of $\C P^{n-1}$ is the $\un$-system \eqref{an}, or equivalently the ${SU(n)\over
SO(n)}$-system. 

\ms
\subsection{\bf Bonnet pairs in $\R^3$ and the ${GL(2,\H)\over (SU(2)\times \R^+)^2}$-system} \label{kp}
\hfil\break

Let $X:M\to \R^3$ be an immersion.  Locally, there exists a conformal coordinate
system
$(x,y)$, i.e., the induced metric is of the form
$\I= e^u(dx^2+dy^2)$ for some smooth function $u$.  Let $H$ denote the mean curvature
function of $M$.  Since $\II - H\I$ is traceless, there is a smooth complex valued function $h=h_1+i
h_2$ such that
$$\II-H \ \I = h_1(dx^2 - dy^2) - 2h_2 dxdy = \Re(hdz^2),$$
where $z= x+iy$.  The two fundamental forms of $M$ are
\begin{equation}\label{dt}
\bca
\I= e^u(dx^2+dy^2), &\cr \II=H\ \I + \Re(h(dx+idy)^2) \cr \ \quad=
(He^u+h_1)dx^2 - 2h_2 dxdy +(He^u - h_1)dy^2. &\cr
\eca
\end{equation}
Let $e_1=Xe^{-{u\over 2}}$, $e_2=Xe^{- {u\over 2}}$, and
$e_3=e_1\times e_2$, where $\times$ is the cross-product.   Let $w_1, w_2, w_3$ be the dual coframe:
$$w_1= e^{u\over 2} dx,\ \  w_2= e^{u\over 2} dy, \ \  w_3=0.$$
Let $g=(e_1, e_2,e_3)$, and $(w_{ij})=g^{-1}dg$, i.e., 
$$de_i =\sum_{j=1}^3 w_{ji} e_j,  \ \  1\leq i\leq 3.$$
Then
\begin{equation}\label{dw}
\begin{cases}
w_{12}={1\over 2}( u_ydx- u_x
dy),\cr 
w_{13}= (He^{{u\over 2}}+ h_1 e^{-{u\over 2}}) dx - h_2 e^{-{u\over 2}}  dy, \cr 
w_{23}= -h_2 e^{-{u\over 2}}dx + (He^{u\over 2} - h_1 e^{-{u\over 2}}) dy,
\end{cases}
\end{equation}
The Gauss-Codazzi equations for $M$ express the flatness of $(w_{ij})$, i.e.,
$$dw_{ij}= -\sum_{k=1}^3 w_{ik}\wedge w_{kj}, \ \ i\not=j.$$  Write this equation in terms of $u,
H, h=h_1+i h_2$ to get
\begin{equation}\label{du}
\bca
u_{xx}+u_{yy} = -2\left(H^2 e^u -(h_1^2+ h_2^2) e^{-u}\right),&\cr
(He^{u\over 2}+ h_1 e^{-{u\over 2}})_y + (h_2 e^{-{u\over 2}})_x = {1\over 2}\left(u_y
(He^{u\over 2}-h_1 e^{-{u\over 2}}) - u_x h_2e^{-{u\over 2}}\right),&\cr
(He^{u\over 2} -h_1 e^{-{u\over 2}})_x + (h_2e^{-{u\over 2}})_y = {1\over
2}\left(u_x(He^{u\over 2}+ h_1 e^{-{u\over 2}})-u_y h_2e^{-{u\over 2}}\right).&\cr
\eca
\end{equation}

\ms

A surface $M$ in $\R^3$ is called {\it isothermic\/} if there exists a conformal line of curvature
coordinate system, i.e., there is a coordinate system $(x,y)$ such that both $\I$ and $\II$ are
diagonalized, or equivalently,  $h_2=0$ in \eqref{dt}.  In this case, the Gauss-Codazzi equations
\eqref{du} become
\begin{equation*}
\bca
u_{xx}+u_{yy} = -2\left(H^2 e^u -h_1^2 e^{-u}\right),&\cr
(He^{u\over 2}+ h_1 e^{-{u\over 2}})_y = {1\over 2}u_y
(He^{u\over 2}-h_1 e^{-{u\over 2}}),&\cr
(He^{u\over 2} -h_1 e^{-{u\over 2}})_x  = {1\over
2}u_x(He^{u\over 2}+ h_1 e^{-{u\over 2}}).&\cr
\eca
\end{equation*}
This implies that $(u,p_1,q_1)$ is a solution of \eqref{fs}, where $p_1={ih_1\over 2} \ e^{-u/2}$ and 
$q_1={H\over 2} \ e^{u/2}$.  

A pair of surfaces $(M,\ti M)$ in $\R^3$  is called  a {\it Bonnet pair\/} if  there is an isometry
$f:M\to \ti M$ so that $\ti H=H\circ f$, where $\ti H$ and $H$ are the mean curvature functions of $M$ and
$\ti M$ respectively and $H$ is not a constant function.  The following is a consequence of the Gauss-Codazzi
equation (cf.  \cite{Bob00p}):

\bprop \label{en} 
 (\cite{Bob00p}). Let $(M, \ti M)$ be a Bonnet pair in $\R^3$. Then  away from umbilic
points there exist a conformal coordinate system 
$(x,y)$, and smooth real functions $u$, $h_1$ and $h_2$ such that
the two fundamental forms for
$M, \ti M$ are given as follows:
\begin{subequations}\label{dv}
\begin{gather}
\bca
\I= e^u(dx^2+dy^2), &\cr \II=H\ \I + \Re(h(dx+idy)^2) \cr \ \quad=
(He^u+h_1)dx^2 - 2h_2 dxdy +(He^u-h_1)dy^2, &\cr
\eca \label{dv1} \\
 \bca
\ti\I= e^u(dx^2+dy^2), &\cr \ti\II=H\ \ti\I + \Re(\bar h(dx+idy)^2) \cr\ \quad =
(He^u+h_1)dx^2 + 2h_2 dxdy +(He^u-h_1)dy^2. &\cr
\eca \label{dv2}
\end{gather}
\end{subequations}
\eprop

Since  both $(u, H, h_1, h_2)$ and $(u,H,h_1, -h_2)$ are solutions of \eqref{du}, we get
\begin{equation} \label{fq}
\bca
u_{xx}+ u_{yy} = 2(-H^2 e^u + (h_1^2+h_2^2) e^{-u}),&\cr
(he^{u/2}+h_1e^{-u/2})_y= {u_y\over 2}(He^{u/2}-h_1e^{-u/2},&\cr
(h_2e^{-u/2})_x=-{u_x\over 2}h_2 e^{-u/2}, &\cr
(He^{u/2}-h_1e^{-u/2})_x= {u_x\over 2}(He^{u/2}+h_1e^{-u/2},&\cr
(h_2e^{-u/2})_y=-{u_y\over 2}(h_2e^{-u/2}.&\cr
\eca
\end{equation}
Note that the third and the fifth equations of \eqref{fq} imply $(h_2)_x=(h_2)_y=0$.  So
$h_2$ is a constant.  So we have

\bthm \label{eo}
 (\cite{Bob00p}).  Let $(M,\ti M)$ be a Bonnet pair in $\R^3$, $(u, H, h_1, h_2)$ the
corresponding solution of \eqref{fq}.  Then $h_2$ is a constant. Moreover,  
set $p_1= {i\over 2} (h_1-ih_2) e^{-{u\over 2}}$, and $ q_1= {1\over 2} H e^{{u\over 2}}$.  
Then $(u,p_1,q_1)$ is a solution of the ${GL(2,\H)\over (\R^+\times SU(2))^2}$-system \eqref{bz}. 
 Conversely, if $(u,p_1, q_1)$ is a solution of system \eqref{bz} and $q_1$ is real, then there is a Bonnet
pair with fundamental forms given by \eqref{dv}, where $H=2q_1 e^{-u/2}$ and $h_1-ih_2=
-2i\ p_1e^{u/2}$.  
\ethm

\subsection{\bf  Curved flats in symmetric spaces}  \label{kq}
\hfil\break

Let $U/U_0$ be a rank $n$ Riemannian symmetric space, $\sigma$ the
corresponding involution on $U$, $\cu=\cu_0+\cu_1$ the eigendecomposition of $d\sigma_e$ on
$\cu$ with eigenvalue $1$ and $-1$, $\ca$ a maximal abelian linear subspace of $\cu_1$, and $a_1,
\cdots, a_n$ an orthonormal basis of $\ca$.
In this section, we associate to each solution of
the $U/U_0$-system \eqref{aa} a flat submanifold in $\cu_1$. We also review the construction of curve
flats in
$U/U_0$ given by Ferus and Pedit \cite{FerPed96a}.

\bthm \label{cz}
    Let $v:\R^n\to \cu_1\cap \ca^\perp$ be a solution of
the $U/U_0$-system \eqref{aa},
and $E(x,\l)$ the frame of the corresponding Lax $n$-tuple $\o_\l$ \eqref{ax}, i.e.,
$$E^{-1}dE= \o_\l= \sum_{i=1}^n (a_i\l + [a_i, v])dx_i, \ \ E(x,0)=e,$$ 
Set $Y= {\p E\over \p \l}E^{-1} \big|\ _{\l=0}$. 
Then  $Y$ is an immersed flat submanifold in $\cu_1$ such that the tangent plane of $Y$ is a maximal abelian
subalgebra of $\cu_1$ at every point.  Conversely, locally all such flat submanifolds in $\cu_1^0$ can be constructed
this way, where $\cu_1^0$ is the subset of regular points in $\cu_1$.  
\ethm

\proof  Write $E_\l(x)= E(x,\l)$.  Since $E^{-1}dE= \sum_i (a_i\l + [a_i, v]) dx_i$, a direct
computation gives
\begin{align*}
dY&= \left({\p\over \p\l} (dE)\right) E^{-1}-{\p E\over \p \l} E^{-1}dE
E^{-1}\bigg|_{\l=0}\cr &= {\p\over \p\l} \left(E\left(\sum_i a_i\l 
+[a_i,v]\right)\ dx_i\right)E^{-1}\bigg|_{\l=0}\cr
& \quad -Y E_0\left(\sum_i[a_i,v]\ dx_i\right )E_0^{-1}\cr
&=\sum_i(Y E_0[a_i,v]E_0^{-1} + E_0a_iE_0^{-1} -YE_0[a_i,v]E_0^{-1})\ dx_i \cr &=\sum_{i=1}^n
E_0a_i E_0^{-1}dx_i.\cr
\end{align*}
 Because $\o_\l$ satisfies the $U/U_0$-reality condition,
$\tau(E_\l)=E_{\bar\l}$ and 
$\sigma(E_\l)=E_{-\l}$.  So $E_0\in U_0$.  Let $$e_i= E_0\ a_i\  E_0^{-1}.$$    Since $a_1, \cdots, a_n$
are orthonormal and
$E_0(x)\in U_0$,
$\{e_i\n 1\leq i\leq n\}$ is an orthonormal tangent frame of $Y$.  Hence $Y$ is an immersion, and  the
induced metric is
$\sum_{i=1}^n dx_i^2$.

Since
$\ad(a_1)^2, \cdots, \ad(a_n)^2$ are commuting symmetric operators, there exist a set $\Lambda$
of linear functionals of $\ca$,  an orthonormal common basis $\{p_\a\n \a\in \Lambda\}$ for
$\ca^\perp\cap\cu_1$, an orthonormal basis
$\{k_\a\n \a\in \Lambda\}$ for $\ck\cap \ck_a^\perp$ such
that 
$$\ad(a)(p_\a)= \a(a) =k_\a, \ \  \ad(a)(k_\a)= -\a(a) p_\a,$$ for all $1\leq i\leq n$ and $\a\in
\Lambda$.  Then
$e_\a= E_0p_\a E_0^{-1}$ is an orthonormal normal frame for the immersion $Y$ in $\cu_1$.   
Write the solution $v$ of the $U/U_0$-system as $v(x)=\sum_\a v_\a(x) p_\a$ with respect to the
decomposition
$\ca^\perp\cap \cu_1= \sum_\a \R\  p_\a$.  
 Since $dE_0= E_0 \sum_i [a_i, v]dx_i$, a direct computation gives
\begin{align*}
de_i&=E_0[E_0^{-1}dE_0, a_i]E_0^{-1}
=\sum_{i,j} E_0[[a_j,v],a_i]E_0^{-1} dx_j, \cr &= -\sum_{j,\a} v_\a\a(a_i) \a(a_j) dx_j e_\a.\cr
\end{align*}
Hence $$w_{i\a}= v_\a\a(a_i)\sum_j
\a(a_j) dx_j.$$  So the normal curvature 
$\sum_\a w_{i\a}\wedge w_{j\a} $ is zero.  

To prove the converse, let $M$ be a flat submanifold of $\cu_1$ such that $TM_p$ is a maximal abelian subalgebra
of $\cu_1$.  Let $x$ be a local flat, orthonormal coordinate of $M$, and $e_i={\p\over \p x_i}$ the orthonormal
frame. Let $\ca$ be a maximal abelian subspace of $\cu_1$.  Then evey maximal abelian subspace of $\cu_1$ is of
the form $k\ca k^{-1}$ for some $k\in U_0$.  Since
$M\subset \cu_1^0$, we may assume that there exist  $\ca$-valued maps
$\xi_i$ and $U_0$-valued map $g$ such that $e_i= g\xi_i g^{-1}$.  It follows from $(de_i, e_j)=0$ and $( , )$ is
ad-invariant that $d\xi_i=0$ for all $i$.  So $\xi_i=a_i$ is constant.  In other words, $e_i= ga_i g^{-1}$.  Note that if
$h$ is $A$-valued map, then $e_i = ga_i g^{-1}= gha_i h^{-1}g^{-1}$.  Choose an $A$-valued map $h$ so that
$(dh)h^{-1} = -\pi(g^{-1}dg)$, where $\pi$ is the projection onto $\cu_{\ca}^\perp\cap \cu_0$.  Let $\ti g= gh$.
Then $e_i = \ti g a_i \ti g^{-1}$ and $\ti g^{-1} d\ti g\in \cu_{\ca}^\perp\cap \cu_0$.  Let $X$ be the immersion
of $M$ into $\cu_1$. Then  
$$dX= \sum_{i=1}^n e_i dx_i = \sum_{i=1}^n \ti g a_i \ti g^{-1} dx_i.$$
Hence $(\ti g a_i \ti g^{-1})_{x_j}= (\ti g a_j \ti g^{-1})_{x_i}$ for all $i,j$.  This implies 
$$[\ti g^{-1} \ti g_{x_j}, a_i]= [ \ti g^{-1} \ti g_{x_i}, a_j].$$
So there exists $\ca^\perp\cap \cu_1$-valued map $v$ such that $\ti g^{-1} \ti g_{x_i} = [a_i, v]$. But this means
$\sum_{i=1}^n [a_i, v] dx_i$ is flat and $v$ is a solution of \eqref{aa}.  \qed

\ms

Given an involution $\sigma$ of $U$, there is a natural $U$-action on $U$ defined by $g\ast x=
gx\sigma(g)^{-1}$.  The orbit at
$e$ is 
$$M=\{g\sigma(g)^{-1}\n g\in U\}.$$ 
Since the isotropy subgroup at $e$ is $U_0$, the orbit $M$ is diffeomorphic to
$U/U_0$.  Next we claim that $M$ is totally geodesic.  To see this, note
that the map $f(g)=(\sigma(g))^{-1}$ is an isometry of
$U$.  So the fixed point set $F$ of
$f$ is a totally geodesic submanifold of $U$.  Note that   
$df_e= -d\sigma_e$.  So $TF_e= \cu_1$, and the dimension of $F$ is equal to
$\dim(\cu_1)$.  But $M$ is fixed by $f$ and $TM_e=\{x-d\sigma_e(x)\n x\in \cu\} =\cu_1$. So $M$
is an open subset of $F$.  This proves the claim.  This is the classical Cartan
embedding of the symmetric space $U/U_0$ in $U$ as a totally geodesic submanifolds. 

Note that $U_0$ acts on $U/U_0$ ($g\cdot (hU_0)= ghU_0$).  An element $x\in U/U_0$ is {\it regular\/} if the
$U_0$-orbit at $x$ is a principal orbit.   

\bthm \label{da}
(\cite{FerPed96a}). With the same assumption as in Theorem \ref{cz}, and set
$$\psi(x)= E(x,1)E(x,-1)^{-1}.$$
Then $\psi$ is an immersed flat submanifold of the symmetric space $U/U_0$ which is tangent to a flat of $U/U_0$
at every point.  Conversely, locally all such flat submanifolds in $N'$ can be constructed this way, where $N'$ is the
open dense subset of regular points in $U/U_0$.
\ethm

\proof  The reality condition implies that $E(x,1)\in U$ and
$$\psi(x)= E(x,1)E(x,-1)^{-1}= E(x,1) \sigma(E(x,1))^{-1}.$$
So the image of $\psi$ lies in the symmetric space $U/U_0=\{g\sigma(g)^{-1}\n g\in U\}$.  
A direct computation gives 
$$\psi^{-1} d\psi= 2\sum_{i=1}^n E_{-1}a_i E_{-1}^{-1} dx_i.$$  Thus $\psi$ is a flat
immersion into $U/U_0$ and $2(x_1, \cdots, x_n)$ is an orthonormal coordinate for the induced metric.  
The rest of the theorem can be proved in a similar manner as for Theorem \ref{cz}.  \qed

\ss
Ferus and Pedit called flat submanifolds obtained in Theorem \ref{da} {\it curved flats\/}. 

\ms
\subsection{\bf  Indefinite affine spheres in $\R^3$ and the $-1$-flow} \label{kr}
\hfil\break

Affine geometry (cf. \cite{NomSas94}) studies the geometry of hypersurfaces in $\R^{n+1}$ invariant
under the affine transformations $x\mapsto Ax+v$, where $A\in SL(n+1,\R)$ and $v\in \R^{n+1}$. There are
three local affine invariants, the affine metric, the Fubini  cubic form, and the third fundamental form. These
invariants satisfy certain integrability conditions, the Gauss-Codazzi equations.  We first give a brief description
of these invariants for affine surfaces in
$\R^3$, then explain the relation between the Tzitzeica equation and indefinite affine spheres. Recall that the
Tzitzeica equation is the $-1$-flow associated to $SL(3,\R)/\R^+$ (see Example \ref{cm}).

Let $X:M\to \R^3$ be a surface with non-degenerate second fundamental form, $g=(e_1, e_2, e_3)$ a
local frame on
$M$ such that
$e_1, e_2$ are tangent to
$M$, $e_3$ is transversal to $M$, and $$\det(e_1, e_2, e_3)=1.$$  Let $w^i$ denote the dual
coframe of $e_i$, i.e.,
$$dX=w^1 e_1 + w^2 e_2.$$
 Let $(w^i_j)$ denote the $sl(3,\R)$-valued $1$-form $g^{-1}dg$, i.e., 
$$de_i = \sum_{j=1}^3 w^j_i e_j, \ \  1\leq i\leq 3.$$
Then we have the structure equation:
\begin{equation} \label{ci}
\bca
dw^i = -\sum_{j=1}^3{w^i_j} \wedge w^j = \sum_{j=1}^3w^j\wedge w^i_j,&
\text{$1\leq i\leq 2$}\cr  dw_i^j = -\sum_{k=1}^3w^j_k \wedge w_i^k. &\text{ $1\leq i, j\leq 3$.}\cr
\eca
\end{equation}
Since $w^3=0$ on $M$, 
\begin{equation} \label{cp}
w_i^3=\sum_{j=1}^2 h_{ij} w^j, \ \  h_{ij}=h_{ji}.
\end{equation}
A direct computation shows that the quadratic form 
\begin{equation} \label{cn}
ds^2=|\det(h_{ij})|^{-{1\over 4}} \sum_{ij} h_{ij} w^i w^j
\end{equation}
is invariant under change of affine frames, and it is called the {\it affine metric\/} of $M$.  An affine
surface is called {\it definite\/} or {\it indefinite\/} if the affine metric is definite or indefinite
respectively.   

We can choose a vector field $e_3$ transversal to $M$ so that 
\begin{equation} \label{ce}
w_3^3 + {1\over 4} d (\log \h)=0.
\end{equation}
Then 
$$\nu= \h^{1\over 4} e_3$$ is  an affine invariant.  The vector field $\nu$ is called the {\it
affine normal\/} of $M$.  

Take the exterior differentiation of \ref{cp} to get
\begin{equation} \label{cr}
\sum_{j=1}^2 dh_{ij} + h_{ij} w_3^3 + \sum_{k=1}^ 2 (h_{ik} w_j^k + h_{kj}
w_i^k)\wedge w^j=0
\end{equation}
for $1\leq j\leq 2$.  Define $h_{ijk}$ by
\begin{equation} \label{cq}
\sum_{k=1}^2 h_{ijk} w^k = dh_{ij} + h_{ij} w_3^3 + \sum_{k=1}^ 2 h_{ik} w_j^k +
h_{kj} w_i^k.
\end{equation}
Then \ref{cr} implies that $h_{ijk}= h_{ikj}$.  But $h_{ij}= h_{ji}$.  So $h_{ijk}$ is symmetric in $i,j,k$. 
{\it The Fubini-Pick cubic form\/},
$$J= \sum_{i,j,k} h_{ijk} w^i w^j w^k,$$
 is an affine invariant.  

Exterior differentiate \eqref{ce} to get
$$\sum_i w_3^i \wedge w_i^3=0.$$
Write 
$$w_3^i= \sum_j \ell^{ij} w_j^3.$$
The {\it third fundamental form\/},
$$\III= h^{1\over 4} w_3^i w_i^3,$$
is also an affine invariant.  The trace of \III with respect to
the affine metric $ds^2$ is the {\it affine mean curvature\/}
$$L={1\over 2} \h^{1\over 4} \sum_{ij} h_{ij}\ell^{ij}.$$

The three affine invariants $ds^2, J$ and $\III$ are completely determined by $w^i$ and $w^A_B$,
which satisfy the Gauss-Codazzi equations for affine surfaces.  
Conversely, suppose $ds^2, J$ and $\III$ are given and satisfy the Gauss-Codazzi equations.  Then $h_{ij},
h_{ijk}$, $w^i, w_i^3$, $w_3^i$, and $w^3_3$ can be computed from these three invariants. Moreover, we can
find
$w_i^j$ by  solving the linear system consisting of \eqref{cq} and the first equation of \eqref{ci}.  Then
the Gauss-Codazzi equations, written in terms of $w^i, w^A_B$, are  \eqref{ci}, i.e., the connection 
$$\W=\bpm
w^A_B & \tau\cr 0&0\cr
\epm$$
is flat, where $\tau=(w^1, w^2, 0)^t$.  Hence there exists 
$$\psi=\bpm
g& X\cr 0&1\cr
\epm$$ such that $\psi^{-1}d\psi= \W$, where $g=(e_1, e_2, e_3)\in
SL(3,\R)$ and
$X\in
\R^3$.   It follows that  
$X$ is an immersion, $e_1, e_2$ are tangent to $X$, $e_3$ is the affine normal, and
$ds^2, J$ and $\III$ are the affine metric, Fubini-Pick form, and the third fundamental form for $X$
respectively.  This is the fundamental theorem of affine surfaces in $\R^3$. 
\ms

A surface is called a {\it proper affine sphere\/} if there exists $p_0\in \R^3$ such that the affine
normal line $p+t\nu(p)$ passes through $p_0$ for all $p\in M$.   We explain below the well-known fact
(cf. \cite{Bob99}) that  the equation for proper affine spheres with indefinite affine metric is the Tzitzeica
equation \eqref{ck}. 

 Let $w$ be a solution of the Tzitzeica equation \eqref{ck}, and $\o_\l$ the
corresponding Lax pair defined by \eqref{cj}, and $E(x,t,\l)$ the solution of 
$$E^{-1}dE=\o_\l, \quad E(0,0,\l)=e.$$ (Here $e$ is the identity matrix in $SL(3,\R)$.)   Fix a
non-zero $r\in \R$, let $e_i(x,t)$ denote the $i$-th column of $E(x,t,r)$.  We claim that $X=-e_3$ is an
immersed indefinite affine sphere.  To see this, we first note that 
$$\o_r=\bpm
w_x dx& r^{-1} e^{-2w} dy& r dx\cr r dx& -w_x dx & r^{-1} e^w dy\cr r^{-1} e^w
dy& r dx & 0\cr
\epm
= E(r)^{-1} dE(r).$$  Since $\o_r$ is $sl(3,\R)$-valued flat $1$-form, $E(r)$ is a map from
$\R^2$ to $SL(3,\R)$.  Fix $r$, and let $e_i$ denote $e_i(r)$.   Equate each column of 
$dE(r)= E(r)
\o_r$  to get
$$\bca
dX= -de_3=- r dx \ e_1 - r^{-1} e^w dy\  e_2,&\cr
de_1= w_x dx\ e_1 + r dx\ e_2 + r^{-1} e^w dy\ e_3,&\cr
de_2= r^{-1} e^{-2w} dy\ e_1 - w_x dx\ e_2+ r dx\ e_3.&\cr
\eca$$
This implies that $e_1, e_2$ are tangent to $X$, 
\begin{align*}
&w^1= -rdx, \ \ w^2= -r^{-1} e^w dy,\cr
& w_1^3 = r^{-1} e^w dy= -w^2, \ \  w_2^3 = r dx=-w^1,\cr
& w_3^1= -w^1, \ \  w_3^2 = -w^2, \ \  w_3^3= 0.
\cr
\end{align*}
So $h_{11}= h_{22}=0$, $h_{12}=-1$ and the affine metric is $2e^w dxdy$.  Since
$\det(h_{ij})= -1$ and $w_3^3 =0$,  \eqref{ce} is satisfied.  Hence the affine normal is 
$$\nu= |\det(h_{ij})|^{1\over 4} e_3=e_3.$$  But
$X= -e_3$ implies that all affine normal lines pass through the origin.  In other words, $X$ is an
indefinite proper affine sphere. 

Conversely, suppose $X$ is an indefinite proper affine sphere in $\R^3$.  We want
to show that there exist a special coordinate system and a special affine frame so that the
Gauss-Codazzi equation for
$X$ as an affine sphere is the Tzitzeica equation.  First note that there exist a local asymptotic
coordinate system
$(x,y)$ and a smooth function
$w$ such that the affine metric is
$$ds^2= 2 e^{2w} dxdy.$$
Let $e_1=X_x$, $e_2= X_y$, and $e_3$ parallel to the affine normal
such that $$\det(e_1, e_2, e_3)=1.$$   Then 
$$w^1= dx, \ \ w^2= dy, \ \  w_1^3= e^{2w} dy, \ \  w_2^3= e^{2w} dx.$$
So $\det(h_{ij})= -e^{4w}$.  
   We may assume that all affine normal lines pass through the origin.   So
$X=fe_3$ for some function $f$.  Exterior differentiation of
$X=fe_3$ gives
$$w^1e_1+ w^2e_2 = df e_3 + f w_3^3 e_3 + f(w_3^1 e_1 + w_3^2 e_2).$$
Equate the coefficients of $e_3$ to get $df + fw^3_3=0$. Since $e_3$ is parallel to the affine
normal,  $w_3^3$ satisfies \eqref{ce}.  Therefore $f= c\n \det(h_{ij})\n^{1/4}= c e^w$ for some
constant $c$.  By scaling, we may assume $c=1$.   Equate coefficients of
$e_1$ and
$e_2$ to get 
$$w_3^1= e^{-w}dx, \ \   w_3^2= e^{-w} dy.$$
Therefore $\ell^{11}=\ell^{22}=0$ and $\ell^{12}=\ell^{21}=e^{-3w}$.  So the affine mean
curvature $L=1$.  Use \eqref{cr} to get
$$w_1^2= -{h_{111}\over 2} e^{-2w} dx, \ \  w_2^1 = -{h_{222}\over 2} e^{-2w} dy.$$
Use $dw^i = -\sum_j w^i_j \wedge w^j$ to conclude
$$w_1^1= w_x dx, \ \  w_2^2 = w_y dy.$$
Substitute $w_A^B$ into $dw^i_j = -\sum_A w^i_A\wedge w^A_j$ for $i=1, j=2$ and $j=1, i=2$ to get
$$(h_{111} e^{-w})_y=0, \ \  (h_{222}e^{-w})_x=0.$$
So $h_{111}= u_1(x) e^w$ and $h_{222}= u_2(y) e^w$ for some smooth function $u_1, u_2$ of one
variable.  By making coordinate change $(\ti x, \ti y)$ so that  $\ti x$ to a function of $x$ and $\ti y$
to a function of $y$, we may assume that 
$$w^1_2= e^{-w} dy, \ \ w^2_1= e^{-w} dx.$$
To summarize, we have shown that 
$$g^{-1}dg= \bpm
w_xdx& e^{-w} dy & e^{-w} dx\cr e^{-w} dx & w_y dy& e^{-w} dy\cr e^{2w} dy&
e^{2w} dx & -dw\cr
\epm,$$
where $g=(e_1, e_2, e_3)$.  Change the frame $g$ to $\ti g=g\diag(1, e^{-w}, e^w)$.
Then $$\ti g^{-1} d\ti g= \bpm
w_xdx& e^{-2w} dy& dx\cr dx& -w_x dx& e^wdy\cr e^wdy &
dx&0\cr
\epm.$$ This is the Lax pair $\o_\l$
\eqref{cj} at $\l=1$.  So $w$ is a solution of the Tzitzeica equation. 

\ms
\subsection{\bf The $-1$ flow, hyperbolic system, and the sigma model}  \label{ks}
\hfil\break

Let $\R^{1,1}$ denote the Lorentz space equipped with metric $2 dxdt$.  
 In this section, we discuss the relation between  harmonic maps from $\R^{1,1}$
 to Lie group $U$ and solutions of the
$-1$-flow and the hyperbolic $U$-system.     

 First we recall a theorem of Uhlenbeck (\cite{Uhl89}):

\bthm \label{ct}
 (\cite{Uhl89}). Let $s:\R^{1,1}\to U$ be a smooth
map, $A= {1\over 2} s^{-1}s_x$, and $B={1\over 2} s^{-1} s_y$.  Then the following statements are
equivalent:
\begin{enumerate}
\item[(i)] $s$ is harmonic,
\item[(ii)] $A_t= -B_x= [A,B]$,
\item[(iii)] $\W_\l= (1-\l)A\  dx + (1-\l^{-1}) B\  dt$ is flat for all $\l\in \C\setminus 0$.  
\end{enumerate}
\ethm

\bcor \label{ft}
 (\cite{Uhl89}).  Suppose $\o_\l=(1-\l) Adx + (1-\l^{-1})B dt$ is flat for all $\l\in
\C\setminus 0$, and
$E_\l$ satisfying $E_\l^{-1} dE_\l = \o_\l$.  Then $s=E_{-1}$ is a harmonic map from $R^{1,1}$ into
$U$ such that $s^{-1} ds= 2A dx+ 2B dt$. 
\ecor 

The following Proposition is well-known:

\bprop \label{ep}
 Let $i:N_0\to N$ be a totally geodesic
submanifold of $N$.   A smooth map $s:M\to N_0$ is a harmonic map if and only if
$i\circ s:M\to N$ is a harmonic map. 
\eprop

\bprop \label{fz}
(\cite{Ter97}).  Let $\cu$ be the real form of $\cg$ defined by the involution $\tau$, 
$a, b\in\cu$ such that $[a,b]=0$,
$g:\R^2\to U$ a solution of the $-1$-flow \eqref{cu} associated to $U$, and $E_\l(x,t)= E(x,t,\l)$ the
frame for the corresponding Lax pair $\o_\l$ \eqref{cv}, i.e., 
$$E^{-1} dE= (a\l + g^{-1}g_x)dx + \l^{-1} g^{-1}bg \ dt, \ \  E(x,t,0)=e.$$   Then
$s=E_{-1}E_1^{-1}$ is a harmonic map from $\R^{1,1}$ to $U$.   Moreover, if $\sigma$ is an order
$k=2m$ automorphism such that $\tau\sigma=\sigma^{-1}\tau$, $a\in \cu\cap\cg_1$, and $
b\in\cu\cap
\cg_{-1}$, then
$s=E_{-1}E_1^{-1}$ is a harmonic map from $\R^{1,1}$ to the symmetric space $U/H$, where $H$ is
the fixed point set of the involution $\sigma^m$ and $\cg_j$ is the eigenspace of $\sigma$ on
$\cg$.  
\eprop

\proof 
Note that the gauge transformation of $\o_\l$ by $E_1$ is
$$E_1\ast \o_\l= E_1\o_\l E_1^{-1} - dE_1 E_1^{-1} =(1- \l)  E_1aE_1^{-1}\ dx + (1-\l^{-1}) E_1
g^{-1}bgE_1^{-1}\ dt.$$
Let $\psi_\l= E_\l E_1^{-1}$.  Then 
$$\psi_\l^{-1} d\psi_\l = E_1\ast \o_\l.$$
By Corollary \ref{ft}, 
$s=\psi_{-1}= E_{-1}E_1^{-1}$ is a harmonic map from $\R^{1,1}$ to $U$ and   $s^{-1} s_x$,
$s^{-1}s_t$ are conjugate to $2a$, $2b$ respectively. 

If the order $k$ of $\sigma$ is $2$, then $\tau\sigma=\sigma^{-1}\tau=\sigma\tau$ and $\sigma$
leaves
$\cu$ invariant.   Let
$K$ denote the fixed point set of
$\sigma$ in $U$, and $\cu=\ck+\cp$ the decomposition of eigenspaces of $\sigma$ on $\cu$ with
eigenvalues
$1, -1$ respectively.  Note that the reality condition is
$$\tau(\o_{\bar \l})= \o_\l, \ \  \sigma(\o_\l)= \o_{-\l}.$$
So $E_\l$ satisfies the reality condition
$$\tau(E_{\bar\l})=\o_\l, \ \  \sigma(E_\l)=E_{-\l}.$$ This implies that 
$s=E_{-1}E_1^{-1}= E_{-1} \sigma(E_{-1})^{-1}$. So the image of $s$ lies in the totally geodesic
submanifold $M=\{g\sigma(g)^{-1}\n g\in U\}$ of $U$.  But $M$ is the Cartan embedding of the
symmetric space $U/U_0$ into $U$ as a totally geodesic submanifold.  It follows from Proposition
\ref{ep} that
$s$ is a harmonic map to the symmetric space $M=U/K$.  

   If the order of $\sigma$ is
$k=2m$, an even integer, then $E_\l$ satisfies the
$(G,\tau,\sigma)$-reality condition
$$\tau(E_{\bar \l})= E_\l, \ \  \sigma(E_\l)=E_{e^{2\pi i\over 2m}\l}.$$
Hence we have $$\sigma^m(E_\l)= E_{-\l}.$$
So $E_1= \sigma^m(E_{-1})$, and $s=E_{-1}\sigma^m(E_{-1})^{-1}$.  Therefore
$s$ is a harmonic map from
$\R^{1,1}$ into the symmetric space
$U/H$. \qed

\beg   Let
$G=SL(2,\C)$,
$\tau(\xi)=-\bar\xi^t$, and
$\sigma(\xi)=-\xi^t$. Note that 
$$SU(2) = \left\{\bpm
w&z\cr -\bar z& \bar w\cr
\epm\ \bigg| \ z, w\in \C, \n w\n^2+\n z\n^2=1\right\} = S^3,$$   and the Cartan
embedding of
$SU(2)/SO(2)$ is the totally geodesic
$2$-sphere 
$$\left\{\bpm 
w& ir\cr ir & -\bar w\cr
\epm\ \bigg|\  r\in \R,\  z\in \C, \ r^2+ \n w\n^2=1.\right\}$$
The $-1$-flow associated to $SU(2)/SO(2)$ is the SGE, so solutions of the SGE give rise to harmonic maps
from $\R^{1,1}$ to  $S^2$.  
\eeg

\beg \label{cw}
 Let $w:\R^2\to \R$ be a solution of the Tzitzeica equation \eqref{ck}, the
$-1$-flow associated to the homogeneous space
${SL(3,\R)\over \R^+}$ given in Example \ref{cm}.  For this example, the order of $\sigma$ is $6$. 
Let $\o_\l$ be the corresponding Lax pair \eqref{cj}, and $E_\l$ the frame of $\o_\l$.  A direct
computation shows that
 $$\sigma^3(\xi)= -PA^t P, \ \ {\rm where\ \ } P=\bpm
0&1&0\cr 1&0&0\cr 0&0&1\cr
\epm.$$
So the fixed point set of the involution $\sigma^3$ in $sl(3,\R)$ is $so(2,1)$, where $SO(2,1)$ is the
isometry group of the quadratic form $2 x_1x_2 + x_3^2$ on $\R^3$.  Hence $E_{-1} E_1^{-1}$ is a
harmonic map from
$\R^{1,1}$ to the symmetric space ${SL(3,\R)\over SO(2,1)}$.  
\eeg

\ms
The proof of Proposition \ref{fz} also implies

\bprop \label{gx}
  Let $(u_0, u_1, v_0, v_1):\R^2\to  \prod_{i=1}^4 \cu$ be a
solution of the hyperbolic system \eqref{ez} associated to $U$, and $E(x,t,\l)$ the frame of the
corresponding Lax pair  \eqref{he}.  Then $s=E(\cdot,\cdot, -1)E(\cdot,\cdot, 1)^{-1}$ is a harmonic map
from $\R^{1,1}$ to
$U$.  Moreover, if $\sigma$ is an involution on $\cu$, $\cu=\cu_0+\cu_1$ is the eigendecomposition,
and $u_0, v_0\in \cu_0$, $u_1,v_1\in \cu_1$, then the image of $s$ lies in the symmetric space
$U/U_0$ (embedded in $U$ via the Cartan embedding) and is a harmonic map from $\R^{1,1}$ to
$U/U_0$. 
\eprop 

\bs

\section{\bf  Dressing actions and factorizations}  \label{kt}

Suppose that $G$ is a Lie group, and that $G_+, G_-$ are subgroups of $G$ such that the multiplication maps
$G_+\times G_-\to G$ and $G_-\times G_+\to G$ defined by $(g_+,g_-)\mapsto g_+g_-$ and
$(g_-,g_+)\mapsto g_-g_+$ respectively are bijections.  Thus given any $g\in G$, there exist uniquely
$g_+\in G_+$ and
$g_-\in G_-$ so that
$g= g_+g_-$, and uniquely $h_+\in G_+$ and $h_-\in G_-$ such that $g= h_-h_+$.  
 The dressing action of $G_+$ on $G_-$ is defined as follows:  Factor
$g_+g_-$ as $\ti g_- \ti g_+$ with $\ti g_\pm\in G_\pm$.  Then the dressing action of $G_+$ on
$G_-$ is 
$g_+\ast g_-=
\ti g_-$. The dressing action of $G_-$ on $G_+$ is defined similarly.  

If the multiplication maps are one to one and the images are open dense
subsets of $G$, then the dressing actions are defined on an open neighborhood of the identity $e$
in $G_\pm$.  Moreover, the corresponding Lie algebra actions are well-defined.

There are two factorizations for a semi-simple Lie group, the Iwasawa and the Gauss factorizations.  The
analogous loop group factorizations are those given by Pressley and
Segal in \cite{PreSeg86} and the Birkhoff factorization respectively.  The dressing actions of these loop
group factorizations play important roles in finding solutions and explaining the hidden symmetries of
integrable systems.  We will review these loop group factorizations.  

\ms
\subsection{\bf Iwasawa and Gauss factorizations}  \label{ku}
\hfil\break

Let $G$ be a complex, semi-simple Lie group, $U$ a maximal compact subgroup, and $B$ a Borel
subgroup. {\it The Iwasawa factorization\/} of $G$ is $G=UB$, i.e., every $g\in G$ can be factored uniquely as
$ub$, where $u\in U$ and $b\in B$.  Let $A$ be a maximal abelian subgroup of $G$, and $N_+$ and
$N_-$ the nilpotent subgroups generated by the set of positive roots and negative roots with respect to a
fixed simple root system of $\ca$ respectively.  The multiplication map from $N_-\times A\times N_+$ to
$G$ is one to one and the image is an open dense subset of $G$, called a {\it big cell\/} of $G$. 
{\it The Gauss factorization\/} is the factorization of a big cell of $G$, i.e., every element $g$ in the big cell can
be factored uniquely as $n_-an_+$ with $n_\pm \in N_\pm$ and
$a\in A$.  Let $B_+=AN_+$.  We call the factorization of $g$ in the big cell as $n_-b_+$ with $n_-\in N_-$
and $b_+\in B_+$ again the Gauss factorization. 

\beg \label{hv}
 Let $G=SL(n,\C)$, $\D_+(n)$ the subgroup of upper triangular $g\in SL(n,\C)$,
and $\D_-(n)$ the subgroup of lower triangular matrix $g\in SL(n,\C)$ with $1$ on the the diagonal.   
The multiplication maps  from
$\D_+(n)\times \D_-(n)\to SL(n,\C)$ and
$\D_-(n)\times \D_+(n)\to SL(n,\C)$ are one to one and the images are open and dense.  Moreover, 
the factorization of $g\in SL(n,\C)$ can be carried out using the Gaussian elimination to rows
and columns of $g$.   This is the Gauss factorization of $SL(n,\C)$. 
\eeg

\beg \label{hw}
 Let $G=SL(n,\C)$, $U=SU(n)$, and $B_+(n)$ the subgroup of upper triangular
matrices with real diagonal entries.  The multiplication maps
$B_+(n)\times U(n)\to SL(n,\C)$ and
$U(n)\times B_+(n)\to SL(n,\C)$ are one to one and onto.  Moreover,  the factorization can be done
by applying the Gram-Schmidt process to rows and columns of $g$.  This is the Iwasawa factorization of
$SL(n,\C)$. 
\eeg

\subsection{\bf  Factorizations of loop groups}  \label{kv}
\hfil\break

Let $G$ be a complex, semi-simple Lie group, $\tau$ an involution of $G$ that gives the compact real
form $U$, and $\sigma$ an order $k$ automorphism of $G$.  Let $B$ be a Borel subgroup of $G$ such
that  $G=UB$ is the Iwasawa factorization.  Given an open subset
$\co$ of
$S^2$, let
$\hol(\co,G)$ denote the group all holomorphic maps $f:\co\to G$ with multiplication defined by
$(fg)(\l)=f(\l)g(\l)$.  
 Let $\e>0$,
\begin{align*}
&S^2=\C\cup\{\infty\},\ \ \  \C^*=\{\l\in \C\n \l\not=0\}, \cr
&\co_\e=\{\l\in \C \n \n\l\n<\e\}, \ \ \  \co_{1/\e}=\{\l\in S^2\n \n\l\n>1/\e\}.\cr
\end{align*}
To explain symmetries of soliton flows we need to consider the following groups:
\begin{align*}
\L(G)&=\hol(\C\cap \co_{1/\e},G),\cr
\L_+(G)&=\hol(\C, G), \cr
\L_-(G)&=\{f\in\hol(\co_{1/\e}, G)\n f(\infty)=e\}.\cr
\end{align*}

The following is the Gauss loop group factorization (the Birkhoff factorization):

\bthm \label{bv} {\bf The Gauss loop group factorization. } The multiplication maps from 
$\L_+(G)\times\L_-(G)$ and $\L_-(G)\times \L_+(G)$ to $\L(G)$
are $1-1$ and the images are open and dense.  In particular,  there exists an open dense subset
$\L(G)_0$ of
$\L(G)$ such that given $g\in \L(G)_0$, $g$ can be factored uniquely as $g=g_+g_-=h_-h_+$ with
$g_+, h_+\in \L_+(G)$ and  $g_-, h_-\in \L_-(G)$.  
\ethm

Suppose $\tau\sigma=\sigma^{-1}\tau$.  Let $\cg=\cg_0+\cdots + \cg_{k-1}$ be the
eigendecomposition of
$\sigma$, and $\cu_j=\cu\cap \cg_j$.   Let $\hat \tau$ and
$\ti \sigma$ be the involution of $\L(G)$ defined by \
\begin{align}\label{gh}
&\hat\tau(g)(\l)=\overline{\tau(g(\bar\l))},\cr
&\ti\sigma(g)(\l)=\sigma(g(e^{-2\pi i/k} \ \l)).\cr
\end{align}
Let $\L^\sigma(G)$ denote the fixed point set of $\ti \sigma$ on $\L(G)$.  
Since $\tau\sigma=\sigma^{-1}\tau$, a direct computation implies that $\ti\tau$ leaves
$\L^\sigma(G)$ invariant.  Let $\L^{\tau,\sigma}(G)$ denote the subgroup of $g\in \L(G)$ that is
fixed by
$\hat \tau$ and
$\ti \sigma$.    Let 
$\L^\tau(G)$ denote
the subgroups of $\L(G)$  fixed by $\hat \tau$.   Then 
\begin{align*}
\L^\tau(G)&=\{f\in \L(G)\n f\ {\rm satisfies\ } U{\rm -reality\  condition \ \eqref{fp}}\},\cr
\L^{\tau,\sigma}(G)&=\{f\in \L(G)\n f \ {\rm satisfies\ } U/U_0{\rm -reality \ condition\ 
\eqref{cb}}\},\cr
\L_\pm^\tau(G)&= \L^\tau(G)\cap \L_\pm(G),\cr
\L_\pm^{\tau,\sigma}(G)&=\L^{\tau,\sigma}(G)\cap
\L_\pm(G).\cr
\end{align*} 

\bcor \label{fi} Suppose $g\in \L(G)$  is
factored as $g=g_+g_-$ with $g_+\in \L_+(G)$ and $g_-\in \L_-(G)$.  If
$\tau\sigma=\sigma^{-1}\tau$, then
\begin{enumerate}
\item[(i)]  $g\in \L^\tau(G)$ implies that $g_\pm\in \L_\pm^{\tau}(G)$,
\item[(ii)] $g\in \L^{\tau,\sigma}(G)$ implies that $g_\pm \in \L_\pm^{\tau,\sigma}(G)$.
\end{enumerate}
\ecor 

To explain symmetries of the elliptic integrable systems, we need to consider the $(G,\tau)$-reality condition
\begin{equation} \label{fw}
\tau(g(1/\bar \l)= g(\l),
\end{equation}
and the
following groups:
\begin{align*}
&L(G)=C^\infty(S^1,G),\cr
&L_+(G)=\{f\in L(G)\n f\ {\rm extends \ holomorphically\  to\ }\n\l\n<1\},\cr
&L_e(U)=\{f\in C^\infty(S^1,U)\n f(1)=e\},\cr
&\cr
&\W^\tau(G)=\{f\in \hol((\co_\e\cup\co_{1/\e})\cap \C^*, G)\n f\  {\rm satisfies\ the \ }
(G,\tau){\rm -reality\/}\cr
&\qquad\qquad\ \  {\rm condition \ \eqref{fw}}\},\cr &\W^\tau_+(G)= \{f\in \W^\tau(G)\n f {\rm \
extends\ holomorphically\ to\ } \C^*\},\cr  
&\W^\tau_-(G)=\{f\in \W^\tau(G)\n f{\rm \ extends\ holomorphically\ to \
}\co_\e\cup\co_{1/\e}, f(\infty)=\I\}.\cr
\end{align*}

For the $(G,\tau,\sigma)$-system, we also need the subgroup of $g\in\W^\tau(G)$ that satisfies the {\it
$(G,\tau,\sigma)$-reality condition\/}:
\begin{equation} \label{fa}
\tau(g(1/\bar\l))=g(\l), \quad \sigma(g(e^{-{2\pi i\over k}}g(\l)).
\end{equation}
 
\bthm \label{gp} (\cite{McI94}).  The multiplication map from $\W_-^\tau(G)\times \W_+^\tau(G)$ to
$\W^\tau(G)$ is a bijection.  
\ethm

\bcor \label{hl}
 Given $g_+\in \W^\tau_+(G)$ and $g_-\in \W_-^\tau(G)$, $g_-g_+$ can
be factored uniquely as $\ti g_+\ti g_-$ with $\ti g_+\in \W^\tau_+(G)$ and $\ti g_-\in \W_-^\tau(G)$. 
Moreover, if $\tau\sigma=\sigma\tau$ and
$g_+, g_-$ satisfy the $(G,\tau,\sigma)$-reality condition \eqref{fa}, then $\ti g_+\in
\W^{\tau,\sigma}_+(G)$ and
$\ti g_-\in\W_-^{\tau,\sigma}(G)$. 
\ecor

\bthm \label{hn} {\bf The Iwasawa loop group factorization\/} (\cite{PreSeg86}).  The multiplication maps
$L_e(U)\times L_+(G)\to L(G)$ and 
$L_+(G)\times L_e(U)\to L(G)$ are bijections.   
\ethm

\bs
\section{\bf  Symmetries of the $U$-hierarchy}  \label{kw}

Let $\tau$ be a conjugate linear involution of $G$, $U$ its fixed point set.  In this Chapter, we assume $U$ is
compact.  Let
$\ca$ be a maximal abelian subalgebra of $\cu$.  The $U$-hierarchy has three types
of symmetries, i.e., three actions on the space of solutions of the $(b,j)$-flow in the $U$-hierarchy:
\begin{itemize}
\item An action of the infinite dimensional abelian algebra $\hat \ca_+$ of polynomial maps from
$\C$ to $\ca\otimes \C$. 
\item An action of $\L_-^\tau(G)$.
\item An action of the subgroup of  $f\in L^\tau_+(G)$ such that the infinite jet of $f_b-\I$ at $\l=-1$ is
$0$, where $f=f_uf_b$ with $f_u\in U$, $f_b\in B$, and $G=UB$ is the Iwasawa factorization.  
\end{itemize}

\ni 
The first two symmetries arise naturally from the dressing actions of the factorization theorems given in
section \ref{kv}  The third action comes from a new factorization.  

\ms
 \subsection{\bf The action of an infinite dimensional abelian group}  \label{kx}
\hfil\break

Let $j>0$ be an integer, $b\in \ca$, $\xi_{b,j}\in \L^\tau_+(\cg)$ defined by
$\xi_{b,j}(\l)= b\l^j$,  and $e_{b,j}(t)$ the one-parameter subgroup of $\L^\tau_+(G)$
defined by $\xi_{b,j}$, i.e.,
$$e_{b,j}(t)(\l)= e^{b\l^j t}.$$  
Let $\hat A_+$ be the subgroup of $\L^\tau_+(G)$ generated by $$\{e_{b,j}(t)\n b\in \ca, j>0 \ {\rm
integer,\ } t\in \R\}.$$  The Lie algebra $\hat\ca_+$ of $\hat A_+$ is the subalgebra of $\L^\tau_+(\cg)$
generated by
$$\{\xi_{b,j}\n b\in
\ca, j>0 \ {\rm integer\ }\}.$$    It follows from Corollary
\ref{fi} that given $f\in
\L^\tau_-(G)$, there is an open subset $\co$ of $0$ in $\R$ such that for all $x\in \co$ there exist $E(x)\in
\L^\tau_+(G)$ and $m(x)\in \L_-^\tau(G)$ such that 
$$f^{-1}e_{a,1}(x)= E(x) m^{-1}(x).$$
Expand $m(x)(\l)$ at $\l=\infty$ to get
$$m(x)(\l)= \I + m_1(x)\l^{-1}+ m_2(x)\l^{-2} + \cdots.$$
Define 
$$\cf(f)= u^f:=[a,m_1].$$
Then $\cf$ is a map from $\L_-^\tau(G)$ to the space $C_0^\infty(\co,\ca^\perp\cap
\cu)$ of germs of smooth maps from $\R$ to $\ca^\perp\cap\cu$ at $0$.  

Given $b\in \ca$, a positive integer $j$, and $f\in \L_-^\tau(G)$, it follows from the Gauss loop group
factorization Theorem \ref{bv} and Corollary
\ref{fi} that there exists  a  neighborhood of $(0,0)$ in $\R^2$ such that for all $(x,t)$ in this
neighborhood we have
$$f^{-1}e_{a,1}(x) e_{b,j}(t)= E(x,t)m(x,t)^{-1}$$ with
$E(x,t)\in\L^\tau_+(G)$ and $m(x,t)\in \L_-^\tau(G)$. A straightforward direct computation
implies that (cf. \cite{TerUhl98}):
\item {(i)} $E^{-1}E_x(x,t,\l)$ must be of the form $a\l + u(x,t)$ and $u=[a,m_1]$, where $m_1$ is
the coefficient of $\l^{-1}$ of the expansion of $m$ as 
$$m(x,t)(\l)=\I + m_1(x,t)\l^{-1} + m_2(x,t)\l^{-2} +\cdots.$$ 
\item {(ii)} $u$ is a solution of the $(b,j)$-flow in the
$U$-hierarchy.  

\ms\ni
By definition of the dressing action,  $m(x,t)= (e_{a,1}(x) e_{b,j}(t))\ast f$.  Hence the $(b,j)$-flow
arises naturally from the dressing action of $\hat A_+\subset \L_+^\tau(G)$ on
$\L_-^\tau(G)$.  Moreover, 
$$u^{e_{b,j}(t)\ast f}= \phi_{b,j}(t)(u^f),$$
where $\phi_{b,j}(t)$ is the one-parameter subgroup generated by the vector field $X_{b,j}$
corresponding the
$(b,j)$-flow in the $U$-hierarchy (i.e., $X_{b,j}$ defined by \eqref{dx}).   In other words, the dressing
action of the abelian group
$\hat A_+$ on
$\L_-^\tau(G)$ induces  an action of $\hat A_+$ on $\cf(\L_-^\tau(G))$, and the action of
$e_{b,j}(t)$ is the $(b,j)$-flow in the $U$-hierarchy. 

\ms
\subsection{\bf The action of $\L_-^\tau(G)$}  \label{ky}
\hfil\break

It follows from the Gauss loop group factorization \ref{bv} that the group $\L_-^\tau(G)$ acts on
$\L^\tau_+(G)$ by local dressing action.  This induces an action of $\L^\tau_-(G)$ on the space of germs of
solutions of  the $(b,j)$-flow \eqref{co} in the $U$-hierarchy  at the origin as follows:  Let
$u:\R^2\to \ca^\perp\cap \cu$ be a solution of the $(b,j)$-th flow \eqref{co}, $\o_\l$ the
corresponding Lax pair
\eqref{ex}, and $E(x,t,\l)$  the frame of
$u$, i.e.,
$E$ is the solution of 
$$\bca
E^{-1}E_{x} = a\l + u, &\cr
E^{-1}E_t=\sum_{i=0}^j Q_{b,j-i}(u)\l^i,&\cr 
E(0,\l)= e.&\cr
\eca$$
 Then $E(x,t)(\l)= E(x,t,\l)$ is holomorphic in $\l\in \C$, i.e., $E(x, t)\in \L_+(G)$ for all $(x,t)$.  Since
$\o_\l$ satisfies the $U$-reality condition
$$\tau(\o_{\bar\l}) = \o_\l,$$
 $E(x,t)$ satisfies 
$$\tau(E(x,t)(\bar \l))=E(x,t)(\l).$$ In other words, $E(x,t)\in
\L_+^\tau(G)$.   Given $g\in \L_-^\tau(G)$, by the Gauss loop group factorization  \ref{bv} and
Corollary
\ref{fi} there is an open subset $\co$ of the origin in $\R^2$ such that the dressing
action of $g$ at $E(x,t)$ is defined for all  $(x,t)\in \co$.  Let
$g\ast E(x,t)$ denote the dressing action of $g$ at $E(x,t)$.  This is obtained as follows:
Factor  $gE(x,t)$ as
$$gE(x,t)= \ti E(x,t) \ti g(x,t)$$
with $\ti E(x,t)\in \L^\tau_+(G)$, and $ \ti g(x,t)\in \L^\tau_-(G)$.  Then 
$$g\ast E(x,t)= \ti E(x,t).$$
Expand $\ti g(x)(\l)$ at  $\l=\infty$:
$$\ti g(x,t)(\l)= I + g_1(x,t)\l^{-1} + g_2(x,t)\l^{-2} + \cdots.$$
The following results are known (c.f. \cite{TerUhl00a}):
\ms
\item {(i)} $\ti u= u + [a,g_1]$ is again a solution of the $(b,j)$-flow.
\item {(ii)} $g\ast u=
\ti u$ defines an action of $\L_-^\tau(G)$ on the space of local solutions of the
$(b,j)$-flow.
\item {(iii)} $\ti E(x,\l)$ is the frame of $\ti u$.
\item {(iv)} Suppose $g\in \L^\tau(G)$ is a rational map with only simple poles, i.e., $g$
satisfies the $U$-reality condition and is of the form
$g(\l)= I +\sum_{j=1}^k {\xi_j\over
\l-\a_j}$ for some $\xi_j\in \cg$ and $\a_j\in \C\setminus\{0\}$.  Then $g\ast u$ can be
computed explicitly in terms of $g, u, E$.   
\item {(v)} If $U$ is compact and $u$ is a smooth solution defined on all $(x,t)\in \R^2$ that is
rapidly decaying in $x$, then $g\ast u$ is also defined on all $\R^2$ and is rapidly decaying in $x$. 

\ms
We claim that $\L_-^{\tau,\sigma}(G)$ acts on the space of solutions of flows in the $U/U_0$-hierarchy. 
To see this,  let $E$ be the frame of a solution $u$ of the $(b,mk+1)$-flow in the $U/U_0$-hierarchy. Then
$E(x,t,\cdot)\in\L^{\tau,\sigma}_+(G)$.  By Corollary \ref{fi}, $g\ast E$ satisfies the
$(G,\tau,\sigma)$-reality condition.  Hence $g\ast u$ is a solution of the $(b,mk+1)$-flow in the
$U/U_0$-hierarchy.  We give an explicit example next. 

\beg \label{dj}
 (\cite{BrDuPaTe02}).  Let $\tau(y)= \bar y$, and $\sigma(y) = \I_{n,n}\ y \ \I_{n,n}^{-1}$ be
the involutions of $O(2n,\C)$ that give the symmetric space $\onn$ as in Example \ref{ai}.  Let
$v=\bpm
0&F\cr -F^t&0\cr
\epm$ be a solution of the $\onn$-system
\eqref{ag},
$\o_\l$ the corresponding Lax $n$-tuple \eqref{dy}, and $E$ the frame of $v$.  We give an explicit
construction of the action of certain rational map with two poles in $\L_-^{\tau,\sigma}(G)$ on $v$
below.  Let  $W, Z$ be two unit vectors of
$\R^n$, $\pi$ the projection of $\C^{2n}$ onto the complex linear subspace spanned by $\bpm
W\cr iZ\cr\epm$, $s\in \R$ non-zero, and
\begin{equation} \label{dz}
h_{is, \pi}(\l)= \left(\pi+ {\l-is\over \l+is} (\I-\pi)\right)\left(\bar\pi + {\l+is\over \l-is}
(\I-\bar\pi)\right).
\end{equation} A direct computation shows that $h_{is,\pi}\in \L_-^{\tau,\sigma}(G)$. 
Then we have:
\begin{enumerate}
\item $h_{is, \pi}E(x)=\ti E(x) h_{is,\ti\pi(x)}$, where $\ti \pi(x)$ is the Hermitian
projection onto the complex linear subspace spanned by 
$$\bpm
\ti W(x)\cr i \ti Z(x)\cr\epm=E(x,-is)^{-1}\bpm 
W\cr iZ\cr\epm.$$
\item $(\ti W, i\ti Z)^t$ is a solution of the following first order system:
\begin{equation} \label{dh}
\bpm
\ti W\cr i\ti Z\cr
\epm_{x_j}= -(-is \ a_j + [a_j, v]) \bpm
\ti W\cr i\ti Z\cr
\epm.
\end{equation} 
\item  Given $v$, system \eqref{dh} is solvable for $\bpm
\ti W\cr i\ti Z\cr\epm$ if and only if $v$
is a solution of \eqref{ag}.  
\item Let  $\xi=(\xi_{ij})$,  $\phi(\xi)=\xi-\sum_i \xi_{ii} e_{ii}$, and $\ti F= F - 2s\
\phi(\ti\pi)$.  Then
$$h_{is,\pi}\ast \bpm 
0&F\cr -F^t&0\cr\epm=\bpm
0&\ti F\cr -\ti
F^t&0\cr\epm.$$  
\item  $\ti E= h_{is,\pi}\ast E$ is the frame of $\ti F= h_{is,\pi}\ast F$. 
\end{enumerate}
\eeg

\ms
\subsection{\bf The orbit $\L_-^\tau(G)\ast 0$} \label{kz}
\hfil\break

Note that $u=0$ is a trivial solution of the $(b,j)$-flow in the $U$-hierarchy and the corresponding Lax
pair is $\o_\l= a\l dx + b\l^j dt$.  So the frame $E^0$ of $u=0$  is 
$$E^0(x,t,\l)= \exp(ax\l  + b\l^j t).$$  
Given $g\in \L_-^\tau(G)$, to compute $g\ast 0$, the first step is to factor 
\begin{equation} \label{gy}
gE^0(x,t) = \ti E(x,t)\ti g(x,t), \ \  {\rm with \ } \ti E(x,t)\in \L_+^\tau(G), \ \ti g(x,t)\in
\L_-^\tau(G).
\end{equation}
 The second step is to expand $\ti g(x,t)$ as 
$$\ti g(x,t)(\l)= \I + \ti g_1(x,t)\l^{-1} + \ti g_2(x,t)\l^{-2} + \cdots.$$
Then $g\ast 0= [a, \ti g_1]$ is a solution of the $(b,j)$-flow.  

The orbit $\L_-^\tau(G)\ast 0$  contains several
interesting classes of solutions (cf. \cite{TerUhl98}):
\begin{enumerate}
\item If $g\in\L^\tau_-(G)$, then $g\ast 0$ is a local analytic solution of the $(b,j)$-flow.
\item  If $g\in \L_-^\tau(G)$ is a rational map with simple poles, then the factorization \ref{gy} can be
carried out using residue calculus and linear algebra.  In fact, $\ti g$ can be given explicitly in terms of a
rational function of exponentials.  Hence  
$g\ast 0$ can be written explicitly as a rational function of exponentials.  
Moreover,  $(g\ast 0)(x,t)$ is defined for all
$(x,t)\in\R^2$, is rapidly decaying as $\n x\n\to \infty$ for each $t\in \R$, and is a pure soliton solution. 
\item  If $g\in\L^\tau_-(G)$ such that $g^{-1}(\l)ag(\l)$ is a polynomial in $\l^{-1}$, then
$g\ast 0$ is a finite type solution, and $g\ast 0$ can be obtained either by solving a system of
compatible first order differential equations or by algebraic geometric method. 
\end{enumerate}
\ms
\subsection{\bf Rapidly decaying solutions}  \label{ma}
\hfil\break

A global solution $u(x,t)$ of the $(b,j)$-flow in the $U$-hierarchy is called a {\it Schwartz class
solution\/} if $u(x,t)$ is rapidly decaying as $\n x\n\to \infty$ for each $t\in \R$. 
The orbit $\L_-^\tau(G)\ast 0$
contains soliton solutions, which are Schwartz class solutions.  But most Schwartz class solutions of the
$(b,j)$-flow  do not belong to this orbit.  We need to use a different loop group factorization than the
ones given in section \ref{kv} to construct general Schwartz class solutions.  

Since we assume $U$ is the compact real form of $G$, there is a Borel subgroup
$B$ such that $G=UB$ (the Iwasawa factorization).    
Let $D_-^\tau$ denote the group of holomorphic  maps
$g:\C\setminus\R\to G$ that satisfying the following conditions:
\begin{enumerate}
\item[(a)] $g$ has an asymptotic expansion at $\l=\infty$ of the form 
$$g(\l) 
 \sim I + g_1\l^{-1} +g_2\l^{-2} +\cdots ,$$ 
\item[(b)] $g$ satisfies the $U$-reality condition \eqref{fp},
\item[(c)] $\lim_{s\to 0^{\pm}} g(r+is) = g_\pm(r)$ is smooth, 
\item[(d)] $h_+-I$ is in the Schwartz class, where $g_+= v_+ h_+$ with $v_+\in U$  and $h_+\in B$.
\end{enumerate}
\ni 
The $U$-reality condition implies that $g_-=\tau(g_+)$.   

The group $D_-^\tau$ is isomorphic to the subgroup of $f\in L_+(G)$
such that
$f-\I$ vanishes up to infinite order at $\l=-1$.  To see this, we consider the following linear
fractional transformation 
\begin{equation} \label{gk}
\l=\phi(z)= i(1-z)/(1+z).
\end{equation}
Note that  $\phi$ has the following properties:
\begin{enumerate}
\item[(i)] $\phi$ maps the unit circle $\n z\n=1$ to the real axis, 
\item[(ii)] $\phi(-1)=\infty$, 
\item[(iii)] $\phi$ maps the unit disk $\n z\n<1$ to
the upper half plane. 
\end{enumerate}

\bthm \label{ga}
(\cite{TerUhl98}).  Let $g\in D_-^\tau$, $\phi$ the linear fractional transformation
defined by \eqref{gk}, and  $\Phi(g)(z)= g(\phi(z))$ for  $\n z\n\not= 1$.  Then:
\begin{enumerate}
\item $\Phi$ is one to one, 
\item $\Phi(g_+)(e^{i\o})$ and $\Phi(g_-)(e^{i\o})$ are the limit of $\Phi(g)(z)$ as $z\to
e^{i\o}$ with $\n z\n<1$ and $\n z\n >1$ respectively.  Moreover, $\Phi(g_+)\in L_+^\tau(G)$.
\item $\Phi(g)$  satisfies the $(G,\tau)$-reality condition,
$$\tau(\Phi(g)(1/\bar\l))=\Phi(g)(\l).$$
\item Let $G=UB$ be the Iwasawa factorization of $G$.  Factor $\Phi(g)=f_uf_b$ with $f_u\in U$ and
$f_b\in B$. Then the infinite jet of
$f_b-\I$ at $z=-1$ is zero.
\end{enumerate}
\ethm

\ms
The above Theorem identifies $g\in D_-^\tau$ with $\Phi(g)\in L_+^\tau(G)$, whose $B$-component is equal
to the identity up to infinitely order at $z=-1$. Recall that the frame of the trivial solution $u=0$ of
the
$(b,j)$-flow in the
$U$-hierarchy is 
$E^0(x,t)(\l)=e^{ax\l+b\l^jt}$.  Note that $\Phi(E^0(x,t))(z)$ is smooth for all $z\in S^1$ except at
$z=-1$, where it has an essential singularity.  So we can {\it not\/} use the dressing action from
the Gauss factorization $L(G)=L_e(U)L_+(G)$ and the identification $\Phi$ to
induce an action of $D_-^\tau$ at 
$u=0$.  However, we can still factor
$gE^0$ as $E \ti g$ with $E\in \L^\tau_+(G)$ and $\ti g\in D_-^\tau$.  Intuitively speaking, the
essential singularity at $z=-1$ is compensated by the infinite flatness of
$\tau(\Phi(g_+))^{-1}\Phi(g_+)(z)$ at $z=-1$.  

\bthm \label{hj}
(\cite{TerUhl98}).  There is an open dense subset $\cd$ of $D_-^\tau$ such that if
$g\in\cd$ then for each $(x,t)\in \R^2$, $gE$ can be factored uniquely as 
$$g(\l)e^{ax\l+b\l^j t}= E(x,t,\l) g(x,t,\l)$$ such that $E(x,t,\cdot)\in \L^\tau_+(G)$ and
$g(x,t,\cdot)\in D_-^\tau$.  Moreover, 
\begin{enumerate}
\item[(i)] $E^{-1}E_x$ is of the form $a\l + u(x,t)$, 
\item[(ii)] $u(x,t)$ is a Schwartz class solution of the $(b,j)$-flow,
\item[(iii)] $E$ is the frame of $u$.  
\end{enumerate}
\ethm    

Let $g\sharp 0$ denote the solution $u$ constructed in Theorem \ref{hj}.  Then:

\bthm \label{gb}
(\cite{TerUhl98}). \begin{enumerate}
\item  $(D_-^\tau\sharp 0)\cap (\L_-^\tau(G)\ast 0) =\{0\}$.  
\item $\L_-^\tau(G)\ast(D_-^\tau\sharp 0)$ is open and dense in the space of Schwartz class
solutions of the $(b,j)$-flow in the $U$-hierarchy. 
\end{enumerate}
\ethm

\ms
\subsection{\bf Geometric transformations}  \label{mb}
\hfil\break

Let $g\in \L_-^{\tau,\sigma}(G)$, $v$ a solution of some integrable system associated to  $U/U_0$,
and
$E$ the frame  of $v$.  Factor $gE$ as $\ti E \ti g$ with $\ti E(x,t)\in \L_+^{\tau,\sigma}(G)$ and
$\ti g\in\L_-^{\tau,\sigma}(G)$.  Then $g\ast v= \ti v$, and $g\ast E= \ti E$ is the frame of
$g\ast v$.   We have seen in Chapter 2 that  geometries associated to solutions of an integrable system
can often be read from their frames at some special value $\l=\l_0$. So $E_{\l_0}\mapsto g\ast
E_{\l_0}$ gives rise to a geometric transformation for the corresponding geometries.  For example, it
is known that solutions of the SGE corresponding to surfaces in $\R^3$ with constant Gaussian
curvature $K=-1$. The Lax pair of the SGE satisfies the $SU(2)/SO(2)$-reality condition, and the dressing action
of
$g_{is,\pi}(\l)= \pi + {\l-is\over \l+is} (\I-\pi)$ on the space of solutions of the SGE corresponds to
the classical B\"acklund transformation of surfaces in $\R^3$ with $K=-1$ (cf. \cite{TerUhl00a}).  
 In this section, we
use another example to demonstrate this correspondence between the dressing action and geometric
transformations.  We describe the  geometric transformation of flat submanifolds that corresponds to
the action of the rational element
$h_{is,\pi}$ (defined by
\eqref{dz}) on the solutions of the $\onn$-system.

First, we need to recall the following definition given by Dajczer and Tojeiro in \cite{DajToj00, DajToj02}.

\bdefn \label{gc}
  Let $M^n$ and $\tilde M^n$ be
submanifolds of  $S^{2n-1}$ with flat normal bundle.  A  vector bundle isomorphism
$P:\nu(M)\to
\nu(\tilde M)$, which covers a diffeomorphism $\ell:M\to \tilde M$, is called a {\it
Ribaucour Transformation\/} if  $P$ satisfies the following properties:
\begin{enumerate}
\item[(a)] If $\xi$ is a parallel normal vector field of $M$,
then $P\circ\xi\circ\ell^{-1}$ is a parallel normal field of
$\tilde M$.
\item[(b)]   Let $\xi\in \nu_x(M)$, and $\g_{x,\xi}$ the normal geodesic with
$\xi$ as the tangent vector at $t=0$.  Then for each $\xi\in \nu(M)_x$, 
$\g_{x,\xi}$ and
$\g_{\ell(x), P(\xi)}$ intersect at a point that is equidistant from $x$ and $\ell(x)$
(the distance depends on $x$).
\item[(c)] If $\eta$ is an eigenvector of the shape
operator $A_{\xi}$ of $M$, then $\ell_{\ast} (\eta)$ is an
eigenvector of the shape operator $A_{P(\xi)}$ of $\tilde M$. Moreover, the geodesics $\g_{x,\eta}$
and $\g_{\ell(x), \ell_\ast(\eta)}$ intersect at a point equidistant to $x$ and $\ell(x)$. 
\end{enumerate}
\edefn

\ms
Dajczer and Tojeiro used geometric methods to prove the existence of Ribaucour transformations between
flat $n$-submanifolds of
$S^{2n-1}$ in \cite{DajToj00}.  These Ribaucour transformations are exactly the one obtained
from dressing actions of $h_{is,\pi}$ given in Example \ref{dj}, i.e.,   

\bthm \label{gd}
(\cite{BrDuPaTe02}). Let $F$ be a solution of the $\onn$-system \eqref{ag}, and $E$ the
frame of the corresponding Lax $n$-tuple \eqref{dy}.  Let
$h_{is,\pi}\in \L_-^{\tau,\sigma}(G)$ defined by \eqref{dz}, $\ti F=h_{is,\pi}\ast F$, and $\ti E=
h_{is,\pi}\ast E$ as in Example \ref{dj}.    Let
$M$ be a flat
$n$-submanifold in
$S^{2n-1}$ associated to $F$ as in Theorem \ref{au}.  Then there exist a flat $n$-submanifold $\ti M$
in
$S^{2n-1}$ and a Ribaucour transformation $P:\nu(M)\to \nu(\ti M)$ constructed from $\ti E=
h_{is,\pi}\ast E$ such that the solution of the
$\onn$-system for $\ti M$ is $\ti F=h_{z,\pi}\ast F$.  
\ethm

\subsection{\bf The characteristic initial value problem for the $-1$-flow}  \label{mc}
\hfil\break

Given $a, b\in \cu$ such that $[a,b]=0$, 
the $-1$-flow in the $U$-hierarchy defined by $a, b$ is the following equation for $g:\R^2\to U$ 
\begin{equation} \label{iv}
(g^{-1}g_x)_t= [a, g^{-1}bg]
\end{equation}
with the constraint $g^{-1}g_x\in [a, \cu]$.   It has a Lax pair 
$$\o_\l=(a\l + g^{-1}g_x) dx + \l^{-1} g^{-1}bg dt.$$
Equation \eqref{iv} is hyperbolic, and $x$-, $t$-curves are the characteristics. The {\it characteristic initial
value  problem\/} (or the {\it degenerate Goursat problem\/}) is the initial value problem with initial data
defined on two characteristic axes, i.e.,  given
$h_1, h_2:\R\to U$ satisfying 
$h_1^{-1}(h_1)_x\in [a,\cu]$ and
$h_1(0)=h_2(0)$, solve
\begin{equation} \label{ja}
\bca
(g^{-1}g_x)_t= [a, g^{-1}bg], &\cr 
g(x,0)=h_1(x), \quad g(0,t)= h_2(t).&\cr
\eca
\end{equation}

If we write 
$$u= g^{-1}g_x, \ \ v= g^{-1}bg,$$ then the $-1$-flow equation \eqref{iv} becomes the following system for
$(u,v)$,
\begin{equation} \label{ix}
u_t= [a, v], \quad v_x= -[u,v],
\end{equation}
The Lax pair is
\begin{equation} \label{iy}
\o_\l= (a\l + u)dx + \l^{-1} v dt.
\end{equation}
Let $M_b$ denote the Adjoint $U$-orbit in $\cu$ at $b$.
Since $u(x,0)= h_1^{-1} h_1'(x)$ and $v(0, t)= h_2(t)^{-1}b h_2(t)\in M_b$, 
the characteristic initial value problem \eqref{ja} become the following initial value problem for \eqref{ix}: 
given
$\xi:\R\to [a,
\cu]$ and
$\eta:\R\to M_b$, find $(u,v):\R^2\to [a,\cu]\times M_b$ so that 
\begin{equation} \label{iz}
\bca
u_t= [a, v], \quad v_x= -[u,v], &\cr 
u(x,0)=\xi(x), \quad v(0,t)= \eta(t).\cr
\eca
\end{equation}

In \cite{DorEit01}, Dorfmeister and Eitner use the Gauss loop group factorization to construct all local
solutions of the Tzitzeica equation \eqref{ck}.  Their construction in fact solves the characteristic initial value
problem
\eqref{iz} for the
$-1$ flow in the $U$-hierarchy:

\bthm \label{jb}
 (\cite{DorEit01}).  Let $\xi,\eta :\R\to [a,\cu]\times M_b$ be smooth maps, and $L_+(x,\l)$ and
$L_-(t,\l)$ solutions of
$$\bca
(L_+)^{-1} (L_+)_x= a\l + \xi(x), &\cr L_+(0,\l)= I,\cr
\eca \quad 
\bca
(L_-)^{-1}(L_-)_t= \l^{-1}
\eta(t), &\cr L_-(0,\l) =I,&\cr
\eca$$
respectively.  Factor 
\begin{equation} \label{jc}
L_-^{-1}(t,\l) L_+(x,\l)= V_+(x,t,\l) V_-^{-1}(x,t,\l)
\end{equation} 
with $V_\pm (x,t,\cdot)\in L_\pm (G)$
via the Gauss loop group factorization.  Set 
$$\phi(x,t,\l)= L_-(t,\l) V_+(x,t,\l)= L_+(x,\l)V_-(x,t,\l).$$
Then $\phi^{-1}\phi_x= a\l + u(x,t)$ and $\phi^{-1} \phi_t= \l^{-1} v(x,t)$ for some $u,v$, and $(u,v)$ solves the
initial value problem of the $-1$-flow  \eqref{iz} in the $U$-hierarchy.  
\ethm

\proof Differentiate $\phi= L_-V_+= L_+V_-$ to get
$$\phi^{-1}\phi_x= V_-^{-1} (a\l +\xi(x)) V_- + V_-^{-1} (V_-)_x = V_+^{-1} (V_+)_x.$$
  So  $\phi^{-1}\phi_x\in \cl_+(\cg)$ and
$$\phi^{-1}\phi_x = \pi_+(V_-^{-1} (a\l +\xi(x)) V_- ),$$
where $\pi_\pm$ is the projection of $\cl(\cg)$ onto $\cl_\pm(\cg)$ with respect to the decomposition $\cl(\cg)=
\cl_+(\cg)+\cl_-(\cg)$.    Expand 
$$V_-(x,t,\l)= I + m_1(x,t)\l^{-1} + \cdots.$$  A direct computation shows that
$$\pi_+(V_-^{-1}(a\l + \xi(x))V_-)= a\l + \xi(x) + [a, m_1(x,t)].$$  Hence
$\phi^{-1}\phi_x= a\l +u(x,t)$, where $u=\xi + [a, m_1]$.    Similar argument implies that 
$$\phi^{-1}\phi_t=\pi_-(V_+^{-1}\l^{-1} \eta(t) V_+) = \l^{-1} g_0(x,t)\eta(t) g_0^{-1}(x,t),$$
where $g_0(x,t)$ is the constant term in the expansion of $V_+(x,t,\l)$
$$V_+(x,t,\l)= \sum_{j=0}^\infty g_j(x,t)\l^j.$$
This implies that $(u,v)$ is a solution of \eqref{ix}, where $$u(x,t)=\xi(x) +[a, m_1(x,t)], \quad v(x,t)=
g_0(x,t)\eta(t) g_0(x,t)^{-1}.$$
It remains to prove $(u,v)$ satisfies the initial conditions.   This can be seen from the factorization \ref{jc}. 
Note that
$$L_-(0,\l)^{-1}L_+(x,\l)= V_+(x,0,\l) V_-(x,0,\l)^{-1}.$$
Since $L_-(0,\l)= I$, the right hand side lies in $L_+(G)$.  Hence $V_-(x,0,\l)= I$, which proves that $m_1(x,0)=0$. 
Therefore $u(x,0)=\xi (x)$.  Similarly argument implies that $v(0,t)= \eta(t)$. \qed

They also show that every local solution of \eqref{ix} can be constructed using suitable $\xi(x)$ and
$\eta(t)$.  To see this, let $(u,v)$ be a solution of \eqref{ix}, and $\phi(x,t,\l)$ the trivialization of the
corresponding Lax pair:
$$\phi^{-1} d\phi= (a\l+ u) dx + \l^{-1} v dt, \quad \phi(0,0,\l)=I.$$
Use Gauss loop group factorization to factor 
$$\phi(x,t,\cdot)= L_-(x,t,\cdot) V_+(x,t,\cdot)= L_+(x,t,\cdot) V_-(x,t,\cdot).$$
Differentiate the above equation to get
\begin{align*}
L_-^{-1}(L_-)_x &= \pi_-(V_+(a\l + u) V_+^{-1}) =0, \cr
L_-^{-1}(L_-)_t&= \pi_-(V_+\l^{-1}v  V_+^{-1})= \l^{-1} g_0vg_0^{-1}, \cr
L_+^{-1}(L_+)_x&= \pi_+(V_-(a\l + u) V_-^{-1})= a\l + u + [m_1, a],\cr
L_+^{-1} (L_+)_t &= \pi_-( V_-\l^{-1} v V_-^{-1})= 0,\cr
\end{align*}
where $g_0$ is the constant term in the power series expansion of $V_+(x,t,\l)$ in $\l$ and $m_1$ is the coefficient
of $\l^{-1}$ of $V_-$.   This implies that $(L_-)_x=0$, $(L_+)_t=0$, $L_+^{-1}(L_+)_x= a\l + u(x,0)$ and
$L_-^{-1}(L_-)_t= \l^{-1}v(0,t)$.

\ms
Let $\tau$ be the involution of $G$ with the real form $U$ as its  fixed point set,  and $\sigma$ an order $k$ 
automorphism of $G$ such that
$\sigma\tau=\tau^{-1}\sigma^{-1}$.  Let $\cg_j$ denote the eigenspace of $\sigma_\ast$ of $\cg$ with
eigenvalue
$e^{2\pi ij\over k}$, and $\cu_j=\cu\cap \cg_j$.  Let
$a\in \cg_1$, $b\in \cg_{-1}$ such that $[a,b]=0$.    The $-1$-flow in the $U/U_0$-hierarchy is the restriction of the
$-1$-flow in the $U$-hierarchy \eqref{ix} to the space of maps
$(u,v):\R^2\to [a,\cu_{-1}]\times
\Ad(U_0)b$.   It is
easy to see that the solution constructed for initial data $\xi:\R\to [a,\cu_{-1}]$ and $\eta:\R\to \Ad(U_0)(b)$ in
Theorem \ref{jb} is a solution of the $-1$-flow in the $U/U_0$-hierarchy.   In other words, the characteristic
initial value problem for the
$-1$ flow in the $U/U_0$-hierarchy can be solved by the algorithm given in Theorem \ref{jb}.

\bs

\section{\bf Elliptic systems associated to $G, \tau, \sigma$}  \label{md}

Let $G$ be a complex Lie group, and $\tau$ a conjugate linear involution of $\cg$, and 
 $\sigma$ an order $k$ complex linear automorphism of $\cg$ such that
$$\tau\sigma= \sigma\tau.$$     Let $\cg_j$ be the
eigenspace of $\sigma$ with eigenvalue $e^{2\pi j i\over k}$.  So we have $\cg_j=\cg_m$ if $j\equiv
m$
$({\rm mod\/}\ k)$,  and
$$\cg=\cg_0+ \cg_1 + \cdots + \cg_{k-1}, \quad [\cg_j,\cg_r]\subset \cg_{j+r}.$$ 
We claim that $\tau(\cg_j)\subset \cg_{-j}$.  To see this, let $\sigma(\xi_j)=\a^j \xi_j$, where
$\a=e^{2\pi i\over k}$.  Then 
$$\sigma(\tau(\xi_j))= \tau(\sigma(\xi_j))= \tau(\a^j \xi_j)= \bar\a^j \tau(\xi_j) =
\a^{-j}\tau(\xi_j).$$
Let $\cu$ be the fixed point set of $\tau$, and  $U_\sigma$ denote the fixed point set of $\sigma$ on $U$. 
Since
$\sigma\tau=\tau\sigma$, $\sigma(U)\subset U$ and $\sigma\n U$ is an order $k$ automorphism of $U$. 
The quotient space 
$U/U_\sigma$ is called a {\it $k$-symmetric space\/}. 

We will construct  the sequences of $U$- and $U/U_\sigma$- systems.   

\ms

\subsection{\bf The $m$-th $(G,\tau)$-system}  \label{me}
\hfil\break

{\it The $m$-th $(G,\tau)$-system\/} (also called the {\it $m$-th elliptic $U$-system\/}) is the  system  for 
$(u_0,\cdots, u_m):\C\to\prod_{i=0}^m \cg$:
\begin{equation} \label{fu}
\bca
(u_j)_{\bar z} =\sum_{i=0}^{m-j}[u_{i+j},\tau(u_i)],&\text{ if $1\leq j\leq m$,}\cr
(u_0)_{\bar z}-(\tau(u_0))_z = \sum_{i=0}^m[u_i,\tau(u_i)].&\cr
\eca
\end{equation}
It has a Lax pair:
\begin{equation} \label{fv}
\o_\l=\sum_{j=0}^m u_j\l^{-j}dz + \tau(u_j)\l^j d\bar z.
\end{equation}
Equation \eqref{fu} is also referred to as the {\it $m$-th elliptic $U$-system\/}, where $U$ is the fixed point
set of $\tau$. 

The Lax pair \eqref{fv} satisfies the
{\it
$(G,\tau)$-reality condition\/} \eqref{fw}, i.e.,
$$\tau(\o_{1/\bar \l})=\o_\l. $$
Note that $\xi=\sum_j \xi_j\l^j$ satisfies the $(G,\tau)$-reality condition if and only if
$\xi_{-j}=\tau(\xi_j)$ for all $j$.  

\ms
\subsection{\bf The $m$-th  $(G,\tau,\sigma)$-system}  \label{mf}
\hfil\break

 \ni {\it The $m$-th $(G, \tau,\sigma)$-system\/} is the equation for 
$(u_0, \cdots, u_m):\C\to\oplus_{j=0}^m \cg_{-j}$,
\begin{equation} \label{ef}
\bca
(u_j)_{\bar z}=\sum_{i=0}^{m- j} [u_{i+j}, \tau(u_i)], &\text{ if 
$1\leq j\leq m$,}\cr -(u_0)_{\bar z} + (\tau(u_0))_z+\sum_{j=0}^m\ [u_j,
\tau(u_j)]=0.&\cr
\eca
\end{equation}
It has a Lax pair
\begin{equation} \label{eg}
\o_\l=\sum_{i=0}^m  u_i\l^{-i} dz + \tau(u_i) \l^i d\bar z.
\end{equation}
Note that: 
\begin{enumerate}
\item[(i)] The $m$-th $(G,\tau,\sigma)$-system  is the restriction
of the $m$-th $(G,\tau)$-system to the space of maps
$(u_0, \cdots, u_m)$ with values in $\oplus_{j=0}^m \cg_{-j}$.
\item[(ii)] The Lax pair of the $m$-th $(G,\tau,\sigma)$-system satisfies the  {\it $(G, \tau,\sigma)$- reality
condition\/} \ref{fa}, i.e.,
$$\tau(\o_{1/ \bar \l})=\o_\l, \ \  \sigma(\o_\l)= \o_{e^{2\pi i\over k} \l}.$$
\item[(iii)] $\xi(\l)= \sum_j \xi_j \l^j$ satisfies the $(G,\tau,\sigma)$-reality condition if and only if 
$\xi_j\in \cg_j$ and $\xi_{-j}=\tau(\xi_j)$ for all $j$.
\end{enumerate}

Let $U$ be the fixed point set of $\tau$, and $U_\sigma$ the fixed point set of $\sigma$ on $U$.  
System \eqref{ef} will also referred to as the {\it $m$-th elliptic $U/U_\sigma$-system\/}. 

\ss
A direct computation gives the following Proposition: 

\bprop \label{jg}
 Let $\tau$ be an involution, and $\sigma$ an order $k$ automorphism of $G$ such that
$\sigma\tau=\tau\sigma$, and  $1\leq m<{k\over 2}$.  If
$\psi:\C\to U$ is a map such that 
\begin{equation} \label{jh}
\psi^{-1}\psi_z= u_o+\cdots + u_m\in \cg_0+\cg_{-1}+\cdots + \cg_{-m},
\end{equation} then $(u_0, \cdots, u_m)$ is a solution of the $m$-th $(G,\tau,\sigma)$-system.  
Conversely, if \goodbreak
\noindent $(u_0,\cdots, u_m)$ is a solution of the $m$-th $(G,\tau,\sigma)$-system \eqref{ef}, then
there exists $\psi:\C\to U$ such that 
$\psi^{-1}d\psi= \sum_{j=0}^m u_j dz + \tau(u_j) d\bar z$.
\eprop

\bdefn \label{jj}
 Let $\tau,\sigma, k$, and $U$ be as in Proposition \ref{jg}, and $1\leq m <{k\over 2}$.  
A map $\psi:\C\to U$ is called a $(\sigma,m)$-map if $\psi^{-1}\psi_z\in \oplus_{j=1}^m \cg_{-j}$. 
\edefn

\bdefn \label{jk}
 (\cite{BurPed94}).  
Let $U/U_\sigma$ denote the $k$-symmetric space (with $k\geq 3$) given by $\tau, \sigma$, and $\pi:U\to
U/U_\sigma$ the natural projection.  A map $\phi:\C\to U/U_0$ is called {\it primitive\/} if there is a lift
$\psi:\C\to U$ (i.e.,
$\pi\circ
\psi= \phi$) so that
$\psi^{-1}\psi_z\in \cg_0+\cg_{-1}$.  In other words, there is a lift $\psi$ that is a $(\sigma,1)$-map.   
\edefn

By Proposition \ref{jg}, the equation for $(\sigma,m)$-maps is the $m$-th $(G,\tau,\sigma)$-system, and the
equation for primitive maps is the first
$(G,\tau,\sigma)$-system.  We refer the readers to \cite{BurPed94} for more detailed study of primitive
maps.

\beg \label{eh}
 Let $\cg=sl(3,\C)$, $\tau(\xi)=-\bar \xi^t$, and $D=
\bpm
0&1&0\cr -1&0&0\cr 0&0&1\cr
\epm$.  Let   
$$\sigma(\xi) = -D\xi^t D^{-1}.$$
Note that $D^{-1}=D^t$, 
$$\sigma^2(\xi)=\diag(-1,-1, 1)\xi\diag(-1,-1, 1),$$
 $\sigma$ has order $4$, and $\sigma \tau= \tau \sigma$.  
A direct computation shows that the eigenspace $\cg_j$ with eigenvalue $(\sqrt{-1}\ )^j$ are:
\begin{align*}
\cg_0&= \left\{\bpm
\xi&0\cr 0& 0\cr\epm
\bigg| \ \xi\in sl(2,\C)\right\},\cr
\cg_1&=\left\{\bpm
0&0 & a\cr 0&0 & b\cr ib& -ia&0\cr\epm \bigg|\  a,b\in \C\right\},\cr
\cg_2&=\C\ \diag(1, 1, -2),\cr
\cg_3&= \left\{\bpm 
0&0 & a\cr 0&0 & b\cr -ib& ia&0\cr\epm \bigg|\ a,b\in \C\right\}.\cr
\end{align*}
The $2$nd $(G,\sigma, \tau)$-system (or the second elliptic $SU(3)/SU(2)$-system) is the system for $(u_0,
u_1, u_2):\C\to
\cg_0\times
\cg_{-1}\times
\cg_{-2}$:
\begin{equation} \label{ek}
\bca
(u_2)_{\bar z}= 0, &\cr
-(u_1)_{\bar z} + [u_1, \tau(u_0)] + [u_2, \tau(u_1)]=0,&\cr
-(u_0)_{\bar z} + (\tau(u_0))_z + [u_0, \tau(u_0)] + [u_1, \tau(u_1)] =0.&\cr
\eca
\end{equation}
\eeg

\beg \label{gq}  (\cite{Gue97, McI94}).  Let $\cg=sl(n,\C)$,
$\tau(\xi)=-\bar
\xi^t$, and
$\sigma(\xi)=C\xi C^{-1}$, where $C=\diag(1, \a,\cdots, \a^{n-1})$ and $\a=e^{2\pi i\over n}$.  Then
$\tau\sigma=\sigma\tau$, and the eigenspace $\cg_j$ of $\sigma$ is spanned by $\{e_{i,i+j}\n 1\leq i\leq
n\}$.  Here we use the notation that $e_{ij}=e_{i'j'}$ if $i\equiv i'$ and $j\equiv j'$ (mod $n$).  The first
$(G,\tau,\sigma)$-system is the equation for 
$A_0= \diag(u_1, \cdots, u_n)$ and $A_1=\sum_{i=1}^n v_i e_{i,i-1}$ so that $$\o_\l=(A_0+A_1\l^{-1})\ dz - 
(\bar A_0^t+ \l \bar A_1^t) \ d\bar z$$ is flat for all $\l$.  
If $v_i>0$ for all $1\leq i\leq n$ and $v_1\cdots v_n=1$, then flatness of $\o_\l$ implies that  we can
write $u_i= (w_i)_z$ and $v_i= e^{w_i-w_{i-1}}$ for some
$w_1, \cdots, w_n$.  The first $(G,\tau,\sigma)$-system written in terms of $w_i$'s is the {\it $2$-dimensional elliptic
periodic Toda lattice\/}:
$$2(w_i)_{z\bar z}= e^{2(w_{i+1}-w_i)}- e^{2(w_i-w_{i-1})}.$$  
\eeg

\ms
\subsection{\bf The normalized  system}  \label{mg}
\hfil\break

{\it The normalized $m$-th  $(G,\tau)$-system\/} or the {\it normalized $m$-th $U$-system\/} is the  system for
$v_1,
\cdots, v_m:\R^2\to \cg$: 
\begin{equation} \label{hq}
(v_j)_{\bar z} =\sum_{i=1}^{m-j} [v_{i+j}, \tau(v_i)] - \sum_{i=1}^m
[v_j,\tau(v_i)],  \ \  \ 1\leq j< m.
\end{equation}
 It has a Lax pair
\begin{equation} \label{hr}
\Theta_\l=\sum_{j=1}^m (\l^{-j}-1) v_j \ dz + (\l^j-1) \tau(v_j)\  d\bar z,
\end{equation}
which satisfies the $(G,\tau)$-reality condition \eqref{fw}.   

We  claim that the Lax pair \eqref{hr} is gauge equivalent to the Lax pair \eqref{fv} for the $m$-th 
$(G,\tau)$-system \eqref{fu}. Hence system \eqref{hq} and \eqref{fu} are gauge equivalent.   To see this, let
$(u_0, \cdots, u_m)$ be a solution of 
\eqref{fu}, 
$$\o_\l= \sum_{j=0}^m u_j\l^{-j} dz + \tau(u_j)\l^j d\bar z$$
its Lax pair, 
and $E(z,\bar z)(\l)$ the frame of $\o_\l$, i.e.,
$$E^{-1} E_z= \sum_{j=0}^m u_j \l^{-j}, \ \ E^{-1} E_{\bar z} = \sum_{j=0}^m \tau(u_j) \l^j, \ \
E(0,0)(\l)=e.$$ Let $g=E(\cdot, \cdot, 1)$.  Since $\o_\l$ satisfies the $(G,\tau)$-reality condition,
$E$ satisfies 
$$\tau(E(\cdot, \cdot, 1/\bar\l)) = E(\cdot, \cdot, \l).$$  Hence $E(\cdot, \cdot, \l)\in U$
if $\n\l\n=1$. In particular, $g\in U$.  The gauge transformation of
$\o_\l$ by $g$ is
\begin{align*}
\ti\o_\l&= g\o_\l g^{-1} - dg g^{-1} = \sum_{j=1}^m (\l^{-j}-1) gu_jg^{-1}dz + (\l^j-1)
g\tau(u_j)g^{-1} d\bar z\cr
&= \sum_{j=1}^m (\l^{-j}-1) gu_jg^{-1} dz + (\l^j-1) \tau(gu_jg^{-1})d\bar z.\cr
\end{align*}
So $(v_1, \cdots, v_m)$ is a solution of \eqref{hq}, where $v_i= gu_j g^{-1}$, and 
$$F(z,\bar z)(\l)= E(z,\bar z)(\l)(E(z,\bar z)( 1))^{-1}$$ is the frame of $\ti \o_\l$.

\bs

\section{\bf Geometries associated to integrable elliptic systems}  \label{mh}

Let $\tau$ be the involution of $G$ whose fixed point set is the maximal compact subgroup $U$ of $G$, and
$\sigma$ an order $k$ automorphism of $G$ such that $\sigma\tau=\tau\sigma$.  Let $\cg_j$ denote the
eigenspace of
$\sigma_\ast$ on $\cg$ with eigenvalue $e^{2\pi ji\over k}$.  
Since $\sigma\tau=\tau\sigma$, we have  $\sigma(U)\subset U$, and
$\sigma\n U$ is an order $k$ automorphism of $U$. Let $U_\sigma$ denote the fixed point set of $\sigma$ in
$U$. The quotient
$U/U_\sigma$ is a symmetric space if
$k=2$, is a {\it $k$-symmetric space\/} if $k>2$.

 It is known that the first $(G,\tau)$-system is the equation for harmonic maps from
$\C$ to $U$ (\cite{Uhl89}).  The first $(G,\tau,\sigma)$-system is the
equation for harmonic maps from
$\C$ to the symmetric space $U/U_\sigma$ if the order of $\sigma$ is $2$,  and is the equation for
primitive maps if $k>2$ (\cite{BurPed94}).   It is also known that a primitive map $\phi:\C\to
U/U_\sigma$ is harmonic  if
$U/U_\sigma$ is equipped with a $U$-invariant metric and $\cg_1$ is isotropic (\cite{BurPed94}). 

The first $(G,\tau,\sigma)$-system also arise naturally in the study of surfaces in symmetric spaces with
certain geometric properties.   For example, constant mean curvature surfaces of simply connected
$3$-dimensional space forms
$N^3(c)$ (\cite{PinSte89, Bob91}), minimal surfaces in $\C P^2$ (\cite{Bur95, BoPeWo95}),
minimal Lagrangian surfaces in
$\C P^2$ (\cite{MaMa01, McI02}), minimal Legendre surfaces in
$S^5$ (\cite{Sha91, Has00}), and special Lagrangian cone in $\R^6= \C^3$ (\cite{Has00, McI02}).  

The only known surface geometry associated to
the $m$-th $(G,\tau,\sigma)$-system for
$m> 1$ was given by H\'elein and Roman.  They showed that the equations for Hamiltonian
stationary surfaces in $4$-dimension Hermitian symmetric spaces are the second elliptic system associated to certain
$4$-symmetric spaces (cf. \cite{HelRom00}).  

If the equation for surfaces with special geometric properties is the $m$-th $(G,\tau,\sigma)$-system, then
the techniques developed for integrable systems can be applied to study the corresponding
surfaces.  In particular, the finite type (or finite gap) solutions give rise to tori with  given
geometric properties.  This has been done for constant mean curvature tori of $N^3(c)$ in \cite{PinSte89} and
\cite{Bob91}, for minimal tori of $\C P^2$ in \cite{Bur95, BoPeWo95}, and for minimal
Legendre tori of
$S^5$ in \cite{Sha91, McI02}. 

 We will give a very brief review of some of the results mentioned above.  For more details, we refer the
readers to \cite{Gue97, Gue01} for harmonic maps, to \cite{PinSte89, Bob91} for constant mean
curvature surfaces in
$3$-dimensional space forms, to \cite{Bur95, BoPeWo95} for minimal surfaces in $\C P^2$, and to
\cite{HelRom00} for Hamiltonian stationary surfaces in four dimensional Hermitian symmetric spaces.

\ms
\subsection{\bf  Harmonic maps from $\R^2$ to $U$ and the first $(G,\tau)$-system}  \label{mi}
\hfil\break

First we state some results of Uhlenbeck (\cite{Uhl89}) on harmonic maps from $\C$ or $S^2$ to $U(n)$.  

\bthm \label{hb}
(\cite{Uhl89}).  Let $G$ be a complex semi-simple Lie group, $\cu$ the real form
defined by the conjugate linear involution $\tau$, 
$s:\C\to U$ a smooth map, and
$A=-{1\over 2} s^{-1}s_z$.  Then the following statements are equivalent:
\begin{enumerate}
\item[(i)] $s$ is harmonic,
\item[(ii)] $A_{\bar z} = -[A, \tau (A)]$, 
\item[(iii)] $(\l^{-1}-1) A\ dz + (\l-1) \tau(A)\ d\bar z$ is flat for all $\l\in \C\setminus \{0\}$, i.e.,
$A$ is a solution of the normalized $1$st  $(G,\tau)$-system.
\end{enumerate}
\ethm

\bcor \label{hc}
(\cite{Uhl89}).  Suppose $$\o_\l = (\l^{-1}-1) A(z,\bar z) dz + (\l -1)
\tau(A(z,\bar z)) d\bar z$$ is a flat
$\cg$-valued $1$-form for all $\l\in \C\setminus 0$, and $E_\l$ the corresponding frame (i.e.,
$E_\l^{-1}dE_\l=\o_\l$ and $E_\l(0)= e$). Then $E_{-1}$ is harmonic from $\C$ to $U$.  
\ecor

Use ellipticity of the harmonic map equation, Uhlenbeck proved that there are trivializations of the Lax pair of
harmonic maps from
$S^2$ to $U(n)$ that are polynomials in the spectral parameter:

\bthm \label{ht}
 (\cite{Uhl89}). Let $s:S^2\to U(n)$ be a harmonic map, and $E$ the frame of the corresponding Lax pair 
$\o_\l=  -{\l^{-1}-1\over 2}\ s^{-1}s_z\ dz -{\l-1\over 2}\ s^{-1}s_{\bar z} \ d\bar z$.  Then there exist
$\g\in L_e(U)$ and smooth maps
$\pi_i:S^2\to
\text{Gr}(k_i,\C^n)$ such that 
$$\g(\l) E(\cdot, \cdot, \l)= (\pi_1+\l \pi_1^\perp) \cdots (\pi_r+\l \pi_r^\perp).$$
\ethm

A harmonic map from a domain of $\C$ to $U(n)$ is called a {\it finite uniton\/} if the corresponding Lax pair admits a
trivialization that is polynomial in $\l$ (\cite{Uhl89}).  The above theorem implies that all harmonic maps  from
$S^2$ to $U(n)$ are {\it finite unitons\/}. 

\ss
A proof similar to that of Proposition \ref{fz} gives 

\bprop \label{hk}
 Let $\tau$ be the involution of $G$ that defines the real form $U$,
 $(u_0,u_1):\C\to \cg\times\cg$ a solution of the first  $(G,\tau)$-system \eqref{fu}, and
$E_\l(z,\bar z)$ the frame of the corresponding Lax pair \eqref{fv}.  Then $s=E_{-1}E_1^{-1}$ is a
harmonic map from $\C$ to $U$. Moreover, let   $\sigma$ be an involution of $\cg$ that commutes
with
$\tau$, and $\cg=\cg_0+\cg_1$ the eigendecomposition of  $\sigma$.  If
$(u_0,u_1)\in \cg_0\times \cg_1$, then $s=E_{-1}E_1^{-1}$ is a
harmonic map from $\C$ to the symmetric space  $U/U_0$.
\eprop

\bcor \label{gr}
   Let $\tau$ be a conjugate liner involution, and $\sigma$ an order $k=2m$
automorphism of $\cg$ such that $\sigma\tau=\tau\sigma$.  Let $(u_1, u_0)$ be a solution of the
first  $(G,\tau,\sigma)$-system
\eqref{ek}, and $E$ the frame of the corresponding Lax pair $\o_\l$ (defined by \eqref{eg}).  Then
$s=E(\cdot,\cdot,{-1})E(\cdot, \cdot, 1)^{-1}$ is a harmonic map from $\C$ to the symmetric space $U/H$,
where $H$ is the fixed point set of the involution
$\sigma^m$ on $U$. 
\ecor

The following Theorem is proved by Burstall and Pedit. 

\bthm \label{jl}
 (\cite{BurPed94}).  Let $(u_0,u_1):\C\to \cg_0\times \cg_{-1}$ be a solution of the first
$(G,\tau,\sigma)$-system, $U/U_\sigma$ the $k$-symmetric space corresponding to $(\tau,\sigma)$, and $k>2$. 
Let
$\pi:U\to U/U_\sigma$ be the natural fibration.  If $E(z,\bar z, \l)$ is a trivialization of the Lax pair of the first
$(G,\tau,\sigma)$-system, then $\phi=\pi\circ E(\cdot, \cdot, 1)$ is primitive. Moreover, if  $U/U_0$ is
equipped with an invariant metric and $\cg_{-1}$ is isotropic, then $\phi$ is harmonic.  
\ethm

\ms
\subsection{\bf  The first $(G,\tau,\sigma)$-system and  surface geometry}  \label{mj}
\hfil\break

Let $N^n(c)$ denote the simply connected space form of constant sectional curvature $c$, i.e., $N^n(0)=\R^n$,
$N^(1)= S^n$ the unit sphere in $\R^{n+1}$, and 
$$N^n(-1)= \H^n=\{x\in \R^{n,1}\n \li x, x\ri =-1\}\subset \R^{n,1}.$$
The Gauss map $\phi$ of a $k$-dimensional submanifold $M$ in $N^n(c)$ is the map from $M$ to the symmetric
space $Y(k,c)$, where $$Y(n,0)= \text{Gr}(k, \R^n), \quad Y(n,1) = \text{Gr}(k, \R^{n+1}), \quad Y(n,-1)=
\text{Gr} (k,
\R^{n,1}).$$  
  A theorem of Ruh and Vilms states
that the Gauss map of a $k$-submanifold with parallel mean curvature vector in $N^n(c)$ is harmonic.  Moreover,  the
Gauss-Codazzi equation for constant mean curvature (CMC) surfaces in $N^3(c)$ is the first
$(G,\tau,\sigma)$-system, where
$\tau,\sigma$ are the involutions that define the symmetric space $Y(3,c)$.  Since the equation for harmonic maps
from
$\C$ to $Y(3,c)$ defined by $\tau, \sigma$ is the  first $(G,\tau,\sigma)$-system, techniques developed for the
first
$(G,\tau,\sigma)$-system (or harmonic maps) can be used to study CMC surfaces in $N^3(c)$ (cf.
\cite{PinSte89, Bob91}).

There are natural definitions of Gauss maps for surfaces and Lagrangian surfaces in $\C P^2$, for
 Legendre surfaces in $S^5$, and Lagrangian cones in $\R^6=\C^3$.   The target manifolds of these
Gauss maps are now $k$-symmetric spaces.  The minimality of surfaces is equivalent to  their Gauss maps being
primitive.  Hence equations of these surfaces are the  corresponding first $(G,\tau,\sigma)$-system.  

\beg \label{jm}
  {\bf Minimal surfaces in $\C P^2$\/}

Let $f:M\to \C P^2$ be an immersed surface, $L\to \C P^2$ the tautological complex line bundle, $z$ a local conformal
coordinate on $M$, and
$f_0$ a local cross section of $f^\ast(L)$.  Choose $f_1, f_2$ so that $(f_0, f_1, f_2)\in SU(3)$ and 
$$\C f_0+\C f_1= \C f_0 + \C {\p f_0\over \p z}.$$
The Gauss map of $M$ is the map $\phi$ from $M$ to the flag manifold $Fl(\C^3)$ of $\C^3$ defined by $\phi(f)=$
the flag $(\C f_0, \C f_0+\C f_1, \C^3)$.  Note that $Fl(\C^3)= SU(3)/T^2$ is a $3$-symmetric space given by
$\tau(g)=(g^t)^{-1}$ and $\sigma(g)=CgC^{-1}$, where $C=\diag(1, e^{2\pi i\over 3}, e^{4\pi i\over 3})$.  It is
proved in \cite{Bur95, BoPeWo95} that $M$ is minimal in $\C P^2$ if and only if  the Gauss map
$\phi:\C\to SU(3)/T^2$ is primitive. 
\eeg

\ms\ni 
\beg \label{jn}
  {\bf Minimal Legendre surfaces in $S^5$\/}

Let $v_1=(1,0,0)^t$, and $\R^3$ the real part of $\C^3$.  Let $Fl_1$ denote the $SU(3)$-orbit of $(v_1, \R^3)$,
i.e.,
\begin{align*}
Fl_1&=\{(gv_1, g(\R^3)\n g\in SU(3)\} \cr &=\{(v, V)\n v\in S^5, v\in V, V \ {\rm is\ Lagrangian\ 
linear\  subspace\ of\ }\C^3\}\cr
&={SU(3)\over1\times SO(2)}\cr
\end{align*}
  Note that $Fl_1$ is a $6$-symmetric space corresponding to automorphisms
$\tau,
\sigma$ of $SL(3,\C)$ defined by
$\tau(g) =(\bar g^t)^{-1}$ and $\sigma(g)=R(g^t)^{-1}R^{-1}$, where $R$ is the rotation that fixes the $x$-axis
and rotate ${\pi\over 3}$ in the $yz$-plane. 

Let $\a$ be the standard contact form on $S^5$. 
A surface $M$ in $S^5$ is {\it Legendre\/} if the restriction of $\a$ to $M$ is zero.  It is
easy to see that $M$ is Legendre if and only if the cone 
$$C(M)=\{tx\n t>0, x\in M\}$$ is Lagrangian in $\C^3$.  If $M\subset S^5$ is Legendre, then there is a natural map
$\phi$ from $M$ to $Fl_1$ defined by  
$\phi(x)= (x, V(x))$, where $V(x)$ is the real linear subspace $\R x+ TM_x$.  It is known that (cf. \cite{Has00,
McI02, MaMa01, Sha91}) that the following statements are equivalent:
\begin{enumerate}
\item[(i)] $M$ is minimal Legendre in $S^5$,
\item[(ii)] the cone $C(M)$ is minimal Lagrangian in $\R^6=\C^5$,
\item[(iii)] the Gauss map $\phi:M\to Fl_1$ is primitive. 
\end{enumerate}

\ms
Let $\pi:S^5\to \C P^2$ be the Hopf fibration, $N$ a surface in $\C P^2$, and $\ti N$ a horizontal lift 
of $N$ in  $S^5$ with respect to the connection $\a$ (the contact form).
Then $N$ is minimal Lagrangian in $\C P^2$ if and only $\ti N$ is minimal Legendre in $S^5$.   Hence there are three
surface geometries associated to the first $(G,\tau,\sigma)$-system associated to the $6$-symmetric space
$Fl_1$: minimal Lagrangian surfaces in
$\C P^2$, minimal Legendre surfaces in
$S^5$, and  minimal Lagrangian cones in $\R^6$.

\eeg

\ms
\beg \label{jo} {\bf Hamiltonian stationary surfaces in $\C P^2$}

Let $N$ be a K\"ahler manifold.  Given a smooth function $f$ on
$N$, let $X_f$ denote the Hamiltonian vector field associated to $f$.   A Lagrangian submanifold $M$
is called {\it Hamiltonian stationary\/} if it is a critical point of the area functional $A$ with respect to
any Hamiltonian deformation, i.e., 
$${\p\over \p t}\bigg|_{t=0} \ A(\phi_t(M))=0$$
for all $f$, where $\phi_t$ is the one-parameter subgroup generated by $X_f$.   This class of submanifolds was
studied by Schoen and Wolfson in \cite{SchWol99}.  When
$N$ is a four dimensional Hermitian symmetric space $U/H$,  H\'elein and Romon proved that the Gauss-Codazzi
equation for Hamiltonian stationary surfaces is the $2$nd $(G,\tau,\sigma)$-system, where $\tau$ is
the involution that gives $U$ and $\sigma$ is an order four automorphism such that $\sigma^2$ gives
rise to the natural complex structure of $U/H$.  In particular, they proved that if $M$ is a Hamiltonian stationary
Lagrangian surface of $\C P^2$, then locally the Gauss-Codazzi equation for $M$ is the $2$nd
$(G,\tau,\sigma)$-system \eqref{ek} given by Example
\ref{eh}. Conversely, if $(u_0, u_1, u_2)$ is a solution of \eqref{ek}, then for each non-zero $r\in \R$,
$E_rE_{-r}^{-1}$ is a Hamiltonian stationary Lagrangian surface of $\C P^2$, where $E_\l$ is the
frame of the Lax pair \eqref{eg} corresponding to $(u_0, u_1, u_2)$. 
\eeg
\bs

\section{\bf Symmetries of the $(G,\tau)$-systems}  \label{mk}

There have been extensive studies on harmonic maps from a
Riemann surface to a compact Lie group $U$.  For example, there are loop group actions, finite unitons, finite
type solutions, and a method of constructing all local harmonic maps from meromorphic data. The
equation for harmonic maps from $\C$ to $U$ is the first $(G,\tau)$-system.   Most results for the first
$(G,\tau)$-system hold for the $m$-th $(G,\tau)$-system as well.  
We will give a brief review here.  For more detail, see \cite{DPW98, Gue97, Gue01, Uhl89}. 

\ms
\subsection{\bf  The action of $\W_-^\tau(G)$}  \label{ml}
\hfil\break

Let $u=(u_0, \cdots, u_m):\C\to \prod_{i=0}^m\cg$ be a solution of the $m$-th 
$(G,\tau)$-system \eqref{fu}, $\o_\l$ the corresponding Lax pair \eqref{fv}, and $E$ the frame of
$\o_\l$, i.e.,
$$
\bca
E^{-1}E_z= \sum_{i=0}^m u_i\l^{-i}, &\cr
E^{-1}E_{\bar z}=\sum_{i=0}^m \tau(u_i) \l^i,&\cr
E(0,0,\l)=\I.&\cr
\eca$$
Let $E(z,\bar z)(\l)=E(z,\bar z, \l)$.  Since $\o_\l$ satisfies the $(G,\tau)$-reality condition,
$$\tau(E(z,\bar z)(1/\bar \l))= E(z,\bar z)(\l),$$ i.e., $E(z,\bar z)\in \W^\tau_+(G)$.  
 Given $g\in \W^\tau_-(G)$, we can use Theorem \ref{gp} to factor
$gE(z,\bar z)= \ti E(z,\bar z) \ti g(z,\bar z)$ with $\ti E(z,\bar z)\in \W^\tau_+(G)$ and $\ti
g(z,\bar z)\in \W^\tau_-(G)$ for $z$ in an open subset of the origin, i.e., $\ti E(z, \bar z)= g\ast E(z,\bar z)$
the dressing action.  A direct computation gives
\begin{equation} \label{hp}
\ti E^{-1}\ti E_z = -\ti g_z \ti g^{-1} + \ti g\left(\sum_{j=0}^m u_j\l^{-j}\right)\ti
g^{-1}.
\end{equation}
 Since $\ti g(z,\bar z)(\l)$ is holomorphic at $\l=0$, the right hand side of \eqref{hp}
has a pole of order at most $m$ at $\l=0$.  Hence there exist some $\ti u_0, \cdots, \ti u_m$ such
that 
$$\ti E^{-1}\ti E_z = \sum_{j=0}^m \ti u_j\l^{-j}.$$
Since $\ti E$ satisfies the $(G,\tau)$-reality condition $\tau(\ti E(1/\bar\l))=\ti E(\l)$, 
$\ti E^{-1}d\ti E$ satisfies the $(G,\tau)$-reality condition \eqref{fw}.  Hence 
$$\ti E^{-1} d\ti E= \sum_{j=0}^m \l^{-j}\ti u_j dz + \l^j\tau(\ti u_j)\ d\bar z.$$
In other words, $\ti u= (\ti u_0, \cdots, \ti u_m)$ is a solution of the $m$-th $(G,\tau)$-system.  Moreover,  $g\ast
u= \ti u$ defines an action of $\W_-^\tau(G)$ on the space of solutions of the $m$-th
$(G,\tau)$-system. This gives the following Theorem of Uhlenbeck \cite{Uhl89} (see also \cite{Gue97}).

\bthm \label{gs}
(\cite{Gue97, Uhl89}).  Let $E$ be the frame of a solution $u$ of the $m$-th 
$(G,\tau)$-system
\eqref{fu}, and $g\in \W^\tau_-(G)$.  Then the dressing action $\ti E(z,\bar z)= g\ast E(z,\bar z)$ is
the frame of another solution $\ti u=g\ast u$.  Moreover, $(g,u)\mapsto g\ast u$ defines an action of
$\W_-^\tau(G)$ on the space of solutions of the $m$-th  $(G,\tau)$-system. 
\ethm

 If $g\in \W^\tau_-(G)$ is a rational map with only simple poles, then the factorization of
$gE(z,\bar z)=\ti E(z,\bar z)
\ti g(z,\bar z)$ with $\ti E(z,\bar z)\in \W^\tau_+(G)$ and $\ti g(z,\bar z)\in \W^\tau_-(G)$ can be
computed by an explicit formula in terms of $g$ and $E$.  In fact, 
  if $g$ has only one simple pole at $\a\in \C\setminus S^1$, then the factorization  can be done by
one of the following methods:
\item {(i)} Equate the residues of both sides of $$g(\l)E(z,\bar z,\l)= \ti E(z,\bar z,\l) \ti g(z,\bar
z,\l)$$ at the pole
$\l=\a$ to get an algebraic formula for $g\ast u$ in terms of $g$ and $E$. 
\item {(ii)} Let $\ti\o_\l= \ti E^{-1}d\ti E$.  Equate the coefficient of $\l^j$ in $\ti \o_\l\ti g= d\ti g
+ \ti g \o_\l$ for each $j$ to get a system of compatible ordinary differential equations. Then  $g\ast u$ can be
obtained from the solution of this system of compatible ODEs.  

\beg \label{ge} 
 (\cite{Uhl89}).  Let $G=GL(n,\C)$, and $\tau(g)=(\bar g^t)^{-1}$.  The fixed point set $U$
of $\tau$ is $U(n)$.    Let $V$ be a complex linear subspace of $\C^n$, 
$\pi$ the Hermitian projection of $\C^n$ onto $V$, $\pi^\perp=\I-\pi$,
$\a\in \C\setminus S^1$, and
$$f_{\a,\pi}(\l)= \pi +\zeta_\a(\l) \pi^\perp,$$
where $\zeta_\a(\l) = {(\l-\a)(\bar \a-1)\over (\bar\a \l -1)(1-\a)}$.  
Note that that $f_{\a,\pi}$ satisfies the $(G,\tau)$-reality condition:
$$\overline{f(1/\bar \l))}^t f(\l)=\I.$$  
If $E$ is the frame of a solution $u$ of the $m$-th  $(G,\tau)$-system \eqref{fu}, then for each
$(z,\bar z)$ the factorization
$f_{\a,\pi}E(z,\bar z)$  must be of the form
\begin{equation} \label{gl}
f_{\a,\pi}E(z,\bar z)= \ti E(z,\bar z) \ti f_{\a, \ti \pi(z,\bar z)}
\end{equation}
 for some $\ti E(z,\bar z)\in \W^\tau_+(G)$ and projection
$\ti \pi(z,\bar z)$.  Use method (i) to conclude that the image
$\ti V(z,\bar z)$ of $\ti \pi(z,\bar z)$ is 
\begin{equation} \label{gm}
\ti V(z,\bar z)= \overline{(E(z,\bar z)(\a))}^t(V).
\end{equation}
 Moreover, 
\begin{equation} \label{hs}
f_{\a,\pi}\ast E= f_{\a,\pi} E f_{\a, \ti \pi(z,\bar z)}^{-1}= (\pi + \zeta_\a(\l)\pi^\perp) \
E \ (\ti \pi + \zeta_\a(\l)^{-1}\ti\pi^\perp)
\end{equation}
is the frame of $f_{\a,\pi}\ast u$.  
For example, if $a\in \cu$ is a constant, then $a$ is a constant solution of
the $1$st normalized   $(G,\tau)$-system
\eqref{fu} with Lax pair $\o_\l=a(\l^{-1}dz+ \l d\bar z)$ and  frame $E_\l(z)=\exp(a\l^{-1}z+ a\l \bar z)$. 
The corresponding harmonic map is 
$$s= E_{-1}(z)E_1^{-1}(z)=  \exp(-2a(z+ \bar z))= \exp(-4ax),\ \ z=x+iy,$$
which is a geodesic.  Since $E$ is given explicitly for the constant solution,  $f_{\a,\pi}\ast a$ is given
explicitly and so is the harmonic map $f_{\a,\pi}\ast s$.  
\eeg

\ms
\subsection{\bf The DPW method and harmonic maps with finite uniton
 number}  \label{mn}
\hfil\break

It is well-known that minimal surfaces in $\R^3$ have Weierstrass representations, i.e., they can  be
constructed from meromorphic functions.   Dorfmeister, Pedit, and Wu gave a construction (the DPW method) of
harmonic maps using meromorphic maps and the Iwasawa loop group factorization (Theorem
\ref{hn}). They call this construction of harmonic maps  the {\it Weierstrass representation of
harmonic maps\/}.  The equation for harmonic maps from $\C$ to
$U$ is the first normalized 
$(G,\tau)$-system.  The DPW method works for the $m$-th normalized 
$(G,\tau)$-system \eqref{hq} as well.  
To explain the DPW method, we need the Iwasawa loop group factorization $L(G)=L_e(U)\times L_+(G)$,
i.e.,  every $g\in L(G)$ can be factored uniquely as
$g=g_1g_2$ with $g_1\in L_e(U)$ and $g_2\in L_+(G)$.  Let $U$ denote the fixed point set of $\tau$.  Recall
that  $L_e(U)$ is the subgroup of $g\in L(U)$ such that
$g(1)=e$ and
$L_+(G)$ the space of smooth loops $g:S^1\to G$ that are boundary value of holomorphic maps defined in
$\n\l\n<1$.   The following Theorem was proved in \cite{DPW98} for the first $(G,\tau)$-system (the
harmonic map equation), but their proof works for the
$m$-th $(G,\tau)$-system as well.  

\bthm \label{hu}
 (\cite{DPW98}). Let $\co$ be a simply connected, open subset of $\C$, and
$\mu(z,\l)=\sum_{j\ge -m} h_j(z) \l^j$ holomorphic in $z\in \co$ and smooth in $\l\in S^1$.   Let
$H:\co\times S^1\to G$ be a solution of
$$\bca
H^{-1} H_z= \sum_{j\geq -m} h_j(z)\l^{j}, &\cr
  H^{-1}H_{\bar z}=0.&\cr
\eca$$ Then:
\begin{enumerate}
\item[(i)] $H$ can be factored as $H(z,\l)=F(z,\bar z,\l)\phi(z,\bar z, \l)$ such that $F(z,\bar z,\cdot)\in
L_e(U)$ and $\phi(z,\bar z,\cdot)\in L_+(G)$.  
\item[(ii)] $F^{-1}F_z$ is of the form $\sum_{i=1}^m (\l^{-i}-1)f_i$ and $f_\mu= (f_1, \cdots, f_m)$
is a solution of the normalized $m$-th $(G,\tau)$-system. 
\end{enumerate} 
\ni Moreover, every solution of the $m$-th $(G,\tau)$-system can be constructed from some $\mu$. 
\ethm  

\proof 
We give a sketch of the proof (for more detailed proof see \cite{DPW98}).  Statement (i) follows from the
Iwasawa loop group factorization \ref{hn}.

The Iwasawa loop group factorization $L(G)=L_e(U)L_+(G)$ implies that there is a Lie algebra
factorization  
\begin{equation} \label{gj}
\cL(\cg)= \cL_e(\cu) + \cL_+(\cg).
\end{equation}
In fact, we can use Fourier series to write
down the Lie algebra factorization easily:  Given
$\xi=\sum_{j\in \Z} \xi_j \l^j$, then $\xi=\eta+\zeta$, where
\begin{align*}
\eta&=\sum_{j=1}^\infty (\xi_{-j}(\l^{-j}-1) +\tau(\xi_{-j})(\l^j-1))\ \in \ \cl_e(\cu),\cr
\zeta&= b_0+ \sum_{j=1}^\infty (\xi_j-\tau(\xi_{-j}))\l^j \ \in \ \cl_+(\cg).\cr
\end{align*}
Note that 
 Since $F= H\phi^{-1}$,
$$F^{-1} dF = \phi H^{-1}dH \phi^{-1} - (d\phi) \phi^{-1}.$$
Let $p_1, p_2$ denote the projection of $\cl(\cg)$ onto $\cl_e(\cu)$ and $\cl_+(\cg)$ with
respect to \eqref{gj}.   Since $(d\phi) \phi^{-1}\in \cL_+(\cg)$ and $F^{-1}dF\in \cL(\cu)$,
$$F^{-1}dF= p_1(\phi (H^{-1} dH) \phi^{-1}).$$ 
 It follows from the fact that $\phi^{-1} d\phi \in \cl_+(\cg)$ and $$H^{-1}dH=\sum_{j\geq
-m} h_j(z)\l^{-j}dz$$ that we have
$$F^{-1} dF=\sum_{j=0}^m f_j(\l^{-j}-1) dz + \tau(f_j)(\l^j-1) d\bar z$$
for some $f_0, \cdots, f_m$.  
This proves (ii).

Let $u=(u_1, \cdots, u_m)$ be a solution of the $m$-th $(G,\tau)$-system, and $E$ a trivialization of the
corresponding Lax pair.   To prove (iii), it suffices to find $h(z,\l)$ so that $g=Eh^{-1}$ is holomorphic in $z\in
\co$.  Since we want
$$g^{-1}dg= h\left(\sum_{j=1}^m (\l^{-j}-1) u_j dz + (\l^j-1)\tau(u_j) d\bar z\right) h^{-1} - dh h^{-1}$$
has no $d\bar z$ term, we must solve $h$ from
$$h^{-1}h_{\bar z} = \sum_{j=1}^m (\l^j-1) \tau(u_j).$$
Since the right hand side lies in $\cl_+(\cg)$,  $h(z,\cdot)$ lies in $L_+(G)$.  Hence
$$g^{-1} dg= \sum_{j=1}^m (\l^{-j}-1) hu_jh^{-1} dz.$$
Hence $g$ is holomorphic in $z$ and $g^{-1}g_z$ is of the form $\sum_{j\ge -m} h_j(z)
\l^j$.  \qed

\ms
It is proved in \cite{DPW98} that finite type solutions arise from constant normalized potential.  We give a
brief explanation next.   Let
$\xi\in L(\cg)$, and $H=\exp(z \xi(\l))$.  So $H^{-1}H_z= \xi(\l)$, $H^{-1}H_{\bar z}=0$, and $H(0,\l)=e$.
Factor 
\begin{equation} \label{ib}
\exp(z\xi)= F(x,y)\phi(x,y), \ \  {\rm with\ } F(x,y)\in L_e(U), \phi(x,y)\in L(G),
\end{equation}
where $ z=x+ i y$. 
Then 
$$H\xi H^{-1}= \exp(z\xi)\xi \exp(-z\xi)=\xi = F\phi \xi \phi^{-1}F^{-1}.$$
This implies that 
\begin{equation} \label{ic}
F^{-1}\xi F= \phi \xi \phi^{-1}.
\end{equation}
Differentiate \eqref{ib} to get  $\xi dz = \phi^{-1} F^{-1}dF \phi + \phi^{-1}d\phi$.  So we have
$$\phi \xi \phi^{-1}dz= F^{-1}dF + d\phi \phi^{-1}.$$
Hence $$F^{-1}dF= p_1(\phi\xi \phi^{-1} (dx+ i\ dy)),$$
where $p_1$ is the projection of $\cl(\cg)$ to $\cl_e(\cu)$.    By \eqref{ic}, we get
$$F^{-1}dF= p_1(F^{-1}\xi F (dx+ i \ dy)).$$
But $d(F^{-1}\xi F)=[F^{-1}\xi F, F^{-1}dF]$.  So we have
$$d(F^{-1}\xi F)= [F^{-1}\xi F, p_1(F^{-1}\xi F (dx+ i\  dy))].$$
Or equivalently,
\begin{equation} \label{id}
\bca
(F^{-1}\xi F)_x=[F^{-1}\xi F, p_1(F^{-1}\xi F)],&\cr
(F^{-1}\xi F)_y=[F^{-1}\xi F, p_1(\sqrt{-1} \  F^{-1}\xi F)].&\cr
\eca
\end{equation}
Let $\xi(\l)=\l^{d-m}V(\l)$.   Then
\eqref{id} becomes
\begin{equation} \label{ie}
\bca
(F^{-1}VF)_x=[F^{-1}VF, p_1(\l^{d-m}F^{-1}VF)], &\cr 
(F^{-1}VF)_y =[F^{-1}VF, p_1(\sqrt{-1} \ \l^{d-m}F^{-1}VF)].&\cr
\eca
\end{equation}
Let $$\eta=F^{-1}VF.$$ Then \ref{ie} becomes
\begin{equation} \label{if}
\bca
\eta_x=[\eta,p_1(\l^{d-m}\eta)], &\cr 
\eta_y=[\eta, p_1(i\ \l^{d-m} \eta)].&\cr
\eca
\end{equation}
Note that this equation leaves the following finite dimensional submanifold of $\cl(\cg)$ invariant:
$$\cv_d=\left\{\eta\in \cl(\cg)\bigg|\  \eta(\l)=\sum_{\n j\n\leq d} \eta_j \l^j\right\}.$$
Hence given $V$ in $\cv_d$, we can solve the ODE system (\eqref{if})
to get
$\eta(x,y)$ such that $\eta(0,0)=V$.  System \eqref{if} is solvable if   
$p_1(\l^{d-m}\eta dz)$ is flat.  So there exists $F(x,y)\in L_e(U)$ such that 
$$F^{-1}dF=p_1(\l^{d-m}\eta dz),$$ i.e., $F$ is  a trivialization of the Lax pair of a solution of the
normalized $m$-th  $(G,\tau)$-system. This is the method of constructing finite type solutions developed by
Pinkall and Sterling in \cite{PinSte89} and  Burstall, Ferus, Pedit and Pinkall in \cite{BFPP93}. 

\ms
All local solutions can also be constructed from meromorphic data $\mu$ that are polynomial in $\l^{-1}$.  To
explain this, we need  

\bthm \label{ih}
 (\cite{DPW98}).  With the same notation as in Theorem \ref{hu}, then
 there exists a discrete set $S\subset \co$ such that for $z\in \co\setminus S$,  $H$ can be factored as
$$H(z,\l)= g_-(z, \l)g_+(z,\l)$$
with $g_-(z,\cdot)\in L^\tau_-(G)$ and $g_+(z,\cdot)\in L^\tau_+(G)$ via the Gauss loop group factorization.
Moreover, 
\begin{enumerate}
\item[(i)] $g_-(z,\l)$ is holomorphic in $z\in \co\setminus S$ and has poles at $z\in S$,
\item[(ii)] $g_-^{-1}dg_-=\sum_{j=1}^m\l^{-j}\eta_j(z) dz$ for some $\cg$-valued meromorphic map
$\eta_j$ on $\co$. 
\end{enumerate}
\ethm

Note that if we factor $g_-$ via the Iwasawa loop group factorization (Theorem \ref{hn}), then the $L_e(U)$
factor of $g_-$ is the same $E$ constructed in Theorem \ref{hu}.  This follows from 
$$g_-= Hg_+^{-1}= E\phi g_+^{-1}= E(\phi g_+^{-1}).$$ 
The converse is also true.  In fact, we have

\bcor \label{jp}
 Let $\mu(\l,z)=\sum_{j=1}^m \l^{-j} \eta_j(z)$ such that $\eta_j$ are meromorphic.  If
there exists
$h(z,\l)$ satisfying  $h^{-1}dh/dz= \mu$, then the $L_e(U)$-factor $E(z,\cdot)$ of
$h(z, \cdot)$ is a trivialization of some solution $f_\mu= (f_1, \cdots, f_m)$ of the $m$-th
$(G,\tau)$-system, i.e., 
$$E^{-1} dE= \sum_{j=1}^m (\l^{-j}-1) f_j(z) dz + (\l^j-1) \tau(f_j) d\bar z.$$
Moreover, every local solution of the $m$-th $(G,\tau)$-system can be constructed this  way. 
\ecor 

The $1$-form
$\mu(z,\l)=\sum_{j=1}^m\eta_j(z)\l^{-j}$ is called the {\it meromorphic potential\/} or the {\it normalized
potential\/}.  However, for general normalized potential $\mu$ the solution $f_\mu$ might have singularities.   An
important problem is to identify meromorphic potentials $\mu$ so that the corresponding solution
$f_\mu$ of the $m$-th $(G,\tau)$-system can be extended to a complete surface.  Burstall and Guest
\cite{BurGue97} have identified
$\mu$'s that give rise to harmonic maps with finite uniton number.  We explain some of their results next. 
  
  Burstall and Guest noted that if $\mu=\l^{-1}h(z)$ is nilpotent and $h(z)$ has no simple poles, then the equation 
\begin{equation} \label{iu}
\bca
H^{-1}H_z = \l^{-1} h(z), &\cr  H^{-1}H_{\bar z}=0&\cr
\eca
\end{equation}
can be solved by integrations.  We use $G=SL(n,\C)$ to explain this.  Let $\cn$ denote the strictly upper triangular
matrices in $sl(n,\C)$, and $h:\co\to \cn$ meromorphic.  To solve \eqref{iu}, we may assume 
$$H(z,\l)= I + b_1(z)\l^{-1} + b_2(z) \l^{-2} + \cdots$$ with meromorphic $b_j$'s.
Equate coefficients of $\l^j$ to get
$$(b_1)_z= h, \  (b_2)_z= b_1 h, \ (b_3)_z= b_2 h, \ \cdots.$$
Since $\cn^n=0$, if we assume that $h(z)$ has no simple poles and the initial data $b_j(0)=0$ then $b_1, \cdots,
b_{n-1}$ can be solved by integration,
 $b_j=0$ for all $j\geq n$, and  $H(z,\l)$ is a polynomial of degree $\leq n-1$ in $\l^{-1}$.   Motivated by this
computation and Uhlenbeck's finite uniton solutions, they make the following definition:

\bdefn \label{jq}
 A harmonic map $s$ from a Riemann surface $M$ to $U$ is said to have {\it finite uniton
number\/} if there is a meromorphic $h:M\to \cg$ such that \eqref{iu} has a solution $H(z,\l)$ satisfying the
following conditions: 
\begin{enumerate}
\item[(i)] $H(z,\l)$  is meromorphic in
$z\in M$ and a polynomial in $\l$ and $\l^{-1}$, 
\item[(ii)]  $s=E(\cdot, -1)$, where $E(z, \cdot)$ is the
$L_e(U)$-component of the Iwasawa factorization of $H(z,\cdot)$. 
\end{enumerate} 
\edefn

\ni In other words, $s$ is the harmonic map constructed from the normalized potential $\mu=\l^{-1} h(z)$.  

\bthm \label{is}
 (\cite{BurGue97, Gue01}).  If $M$ is a Riemann surface and $s:M\to U$ is harmonic map with finite
uniton number, then there exists a complex extended solution $H$ (associated to $s$) of the form 
$$H(z,\l)= \exp(\l^{-1} b_1(z) + \cdots + \l^{-r} b_r(z)),$$
where $b_1, \cdots, b_r$ are meromorphic maps from $M$ to the nilpotent subalgebra $\cn$ of the Iwasawa
decomposition $\cg=\ck +\ca + \cn$.  Moreover, 
\begin{enumerate}
\item[(i)] integer $r$ can be computed in terms of root system of $G$,
\item[(ii)] the maps $b_2, \cdots, b_r$ satisfies a meromorphic ordinary differential equation, which can be
solved by quadrature for any choice of $b_1$.  
\end{enumerate}
\ethm

In fact, the normalized potential $\mu$ corresponding to the harmonic map constructed by Theorem \ref{is}
is
$\mu= \l^{-1}(b_1)_z$.

\bs
\subsection{\bf Some comparisons}  \label{mo}
\hfil\break

Let $G$ be a complex, semi-simple Lie group, and $U$ the maximal compact subgroup of $G$, and $\tau$ the
corresponding involution with fixed point $U$.   We have
discussed the constructions of solutions of soliton equations in the
$U$-hierarchy in Chapter 6 and of equations in the 
$(G,\tau)$-hierarchy in section 7.2.  Loop group factorizations are used in both cases.  In this section, we give a
summary and some comparisons of these constructions of solutions for the two hierarchies.  To make the
exposition easier to follow, we will not give references in this section (for references see previous sections).    

Let $\ca$ be a maximal abelian subalgebra of $\cu$, and $a\in \ca$ a regular element.  
For the $U$-hierarchy defined by $a$, the data we use to construct solutions for the $(b,j)$-flow in the
$\cu$-hierarchy of soliton equations is one of the following types of maps:
\begin{enumerate}
\item[(i)] $f$ is a holomorphic map from a neighborhood of $\l=\infty$ in $S^2=\C\cup \{\infty\}$ to $G$ that
satisfies the $U$-reality condition, $\tau(f(\bar \l))= f(\l)$, and $f(\infty)=\I$, 
\item[(ii)] $f:\R\to G$ is smooth, has an asymptotic expansion at $\infty$, $f(\infty)=\I$,  $f$ is the
boundary value of a holomorphic map on the upper half plane, and $f_b$ is
rapidly decaying at infinity, where $f=f_uf_b$ is the pointwise Iwasawa factorization of $G=UB$, i.e., $f_u\in U$
and $f_b\in B$,
\item[(iii)] $f=f_1f_2$, where $f_1:S^2=\C\cup\{\infty\}\to G$ is a rational map of type (i) and 
$f_2$ is of type (ii),  
\item[(iv)] $f=f_1f_2$, where $f_1$ is of type (i) and $f_2$ is of type (ii). 
\end{enumerate}
\ss
\ni  
To construct solutions, we start with an $f$ of type (i), (ii), (iii), or (iv), then factor
$f^{-1}e_{a,1}(x)e_{b,j}(t)$ as
$ E(x,t)m(x,t)^{-1}$ with
$E(x,t)\in
\L^\tau_+(G)$ and $m(x,t)$ of type (i), (ii), (iii) or (iv) accordingly, where  $e_{\xi,j}(t)=
e^{\xi\l^j t}$.  Then 
$$u^f(x,t)=[a,m_1(x,t)]$$
is a solution of the $(b,j)$-flow in the $U$-hierarchy,
where $m_1(x,t)$ is the coefficient of $\l^{-1}$ in the expansion of $m(x,t)(\l)$ at $\l=\infty$:
$$m(x,t,\l)\ \sim \ \I + m_1(x,t)\l^{-1} + \cdots.$$
Moreover,  we know:
\begin{enumerate}
\item $u^f=u^g$ if and only if $f=hg$ for some $A$-valued map $h$. 
\item $u^f$ is a local real analytic solution if $f$ is of type (i).
\item If $f$ is of type (iii), then $u^f(x,t)$ is a solution defined for all $(x,t)\in\R^2$ and is rapidly decaying in $x$
for each fixed $t$.  The space of such solutions $u^f$ is open and dense in the space of all rapidly decaying
solutions.
\item  If $u$ is a finite gap solution (an algebraic geometric solution described by theta functions), then
there exists an
$f$ of type (i) such that
$f^{-1}af$ is a polynomial in $\l^{-1}$ and $u=u^f$.
\item If $f$ is of type (i) and is a rational map from $S^2$ to $G$, then $u^f$ is a pure soliton solution.
\end{enumerate}

\ms
For the normalized $m$-th  $(G,\tau)$-system, we start with meromorphic potential $\mu(z,\l)=\sum_{j=1}^m
\eta_j(z)\l^{-j}\ dz$.  There are two
steps to construct a solution:

\ni  Step 1.   Find a solution $H(z,\l)$ of
$H^{-1}dH= \mu$ that is smooth for all $\l\in S^1$ and
meromorphic in $z\in \co\subset \C$.

\ni Step  2.  Factor $H$ as $F\phi$ with $F\in L_e(U)$ and
$\phi\in L_+(G)$.  Then $F^{-1}F_z$ is of the form $\sum_{j=1}^m (\l^{-j}-1)v_i$ for some $v_1, \cdots,
v_m$.  Hence 
$$v_\mu=(v_1, \cdots, v_m), \qquad s_\mu=F(\cdot, -1)$$ 
are a solution of the normalized $m$-th  $(G,\tau)$-system and a harmonic map from $\co$ to  $U$
respectively.   

\ms
For the first normalized  $(G,\tau)$-system, to go beyond solutions with finite uniton numbers we note
that:
\begin{enumerate}
\item[\di]  There is no simple condition on $\mu$ to guarantee that Step 1 can be done.   
\item[\di] Every local smooth solution can be constructed from some $\mu$.  
However, in general, there is no canonical choice of $\mu$.
\item[\di] One of the main open problems is to identify the set of $\mu$ so that $s_\mu$ can be extended to
a harmonic map on a closed surface.
\end{enumerate}

\bigskip

\providecommand{\bysame}{\leavevmode\hbox to3em{\hrulefill}\thinspace}

\vskip 0.25in 


\begin{thebibliography}{10}

\bibitem{AblCla91}
Ablowitz, M.J.,Clarkson, P.A.,\emph{{S}olitons,
non-linear evolution equations and inverse scattering},  Cambridge Univ.
Press (1991)

\bibitem{AKNS74}
Ablowitz, M.J., Kaup, D.J., Newell, A.C. and Segur, H., \emph{{T}he inverse
scattering transform - Fourier analysis for nonlinear
problems}, Stud. Appl. Math. \textbf{53} (1974), 249--315

\bibitem{BeaCoi84}
Beals, R., Coifman, R.R.,\emph{{S}cattering and inverse scattering for
first order systems}, Commun. Pure Appl. Math. \textbf{37} (1984), 39--90.

\bibitem{Bob91}
Bobenko, A.I., \emph{{A}ll constant mean curvature tori in $R^3$, $S^3$, $H^3$ in terms of theta
functions}, Math. Ann. \textbf{290} (1991), 1--45

\bibitem{Bob99}
Bobenko, A.I., \emph{{D}iscrete indefinite affine spheres}, Discrete integrable geometry and physics,
Oxford Lecture Ser. Math. Appl., Oxford Univ. Press, New York, \textbf{16} (1999), 113--138

\bibitem{Bob00p}
Bobenko, A.I., \emph{{E}xploring surfaces through methods for the theory of integrable
systems}, lectures on Bonnet problem, preprint, math.DG/0006216

\bibitem{BoPeWo95}
Bolton, J., Pedit, F., and Woodward, L., \emph{{M}inimal surfaces and the affine field
model}, J. reine angew. Math. \textbf{459} (1995), 119--150

\bibitem{BrDuPaTe02}
Br\"uck, M., Du, X., Park, J., and Terng, C.L., \emph{{S}ubmanifold geometry of real 
Grassmannian systems}, The Memoirs, vol 155, AMS,  \textbf{735} (2002), 1--95

\bibitem{Bur95}
Burstall, F., \emph{{H}armonic tori in spheres and complex projective spaces}, J. reine angew.
Math. \textbf{469} (1995), 149--177

\bibitem{BFPP93} 
Burstall, F.E., Ferus, D., Pedit, F., Pinkall, U.\emph{{H}armonic
tori in symmetric spaces and commuting Hamiltonian systems on loop
algebras}, Annals of Math. \textbf{37} (1993), 173--212

\bibitem{BurGue97}
Burstall, F.E., Guest, M.A. \emph{{H}armonic two-spheres in compact symmetric
spaces}, Math. Ann.\textbf{309} (1997), 541--572

\bibitem{BurPed94}
Burstall, F.E., Pedit, F., \emph{{H}armonic maps via Adler-Kostant-Symes Theory}, Harmonic
maps and Integrable Systems, Vieweg (1994), 221--272

\bibitem{BurPed95}
Burstall, F.E., Pedit, F., \emph{{D}ressing orbits of harmonic maps}, Duke
Math. J. \textbf{80} (1995), 353--382

\bibitem{CieGolSym95}
Cie\'sli\'nski, J., Goldstein, P.,  and Sym, A., \emph{{I}sothermic surfaces in
  {${E}^3$} as soliton surfaces}, Phys.\ Lett.\ A \textbf{205} (1995), 37--43.

\bibitem{Cie97A}
Cie{\'s}li{\'n}ski, J., \emph{{T}he {D}arboux-{B}ianchi transformation for
  isothermic surfaces. {C}lassical results versus the soliton approach},
  Differential Geom.\ Appl. \textbf{7} (1997), 1--28.

\bibitem{DajToj95a}
Dajczer, M., Tojeiro, R., \emph{{I}sometric immersions and the generalized Laplace and elliptic
sinh-Gordon equations}, J. reine
angew Math, \textbf{467} (1995), 109--147

\bibitem{DajToj95b}
Dajczer, M., Tojeiro, R., \emph{{F}lat totally real submanifolds of $C P^n$ and the symmetric generalized
wave equation},  Tohoku Math. J., \textbf{47} (1995), 117--123

\bibitem{DajToj00}
Dajczer, M., Tojeiro, R., \emph{{T}he Ribaucour transformation for flat Lagrangian submanifolds}, J. Geom.
Anal. \textbf{10} (2000), 269--280

\bibitem{DajToj02}
Dajczer, M., Tojeiro, R., \emph{{A}n extension of the classical Ribaucour transformation}, Proc.
London Math. Soc. \textbf{85} (2002), 211--232.

\bibitem{DajToj02p}
Dajczer, M., and Tojeiro, R., \emph{{C}ommuting Codazzi tensors and the
  Ribaucour transformation for submanifolds}, preprint.

\bibitem{Dar10}
 Darboux, G., \emph{{L}econ sur les Syst\`emes
Orthogonaux}, Gauthier-Villars, (1910), 2nd edition

\bibitem{DPW98}
Dorfmeister, J., Pedit, F., and Wu, H., \emph{{W}eierstrass type
representation of harmonic maps into symmetric spaces}, Comm. Anal.
Geom., \textbf{6} (1998), 633--668

\bibitem{DorEit01}
Dorfmeister, J., and Eitner, U., \emph{{W}eierstrass type representation of affine spheres}, preprint, to
appear in Abh. Math. Univ. Hamburg

\bibitem{Dor02}
Dorfmeister, J., \emph{{G}eneralized Weierstrass representations of surfaces}, preprint

\bibitem{FerPed96a}
Ferus, D., Pedit, F., \emph{{C}urved flats in symmetric
spaces}, Manuscripta Math., \textbf{91} (1996), 445--454

\bibitem{FerPed96b}
Ferus, D., Pedit, F., \emph{{I}sometric immersions of space forms and
soliton theory}, Math. Ann., \textbf{305} (1996), 329--342

\bibitem{Gue97}
Guest, M., \emph{{H}armonic maps, loop groups, and integrable systems}, Cambridge
University  Press, (1997)

\bibitem{Gue01}
Guest, M., \emph{{A}n update on harmonic maps of finite uniton numbers, via the zero curvature
equation}, preprint,  math.dg/0105110

\bibitem{Has00}
Haskins, M., \emph{{S}pecial Lagrangian cones}, preprint, math.dg/0005164

\bibitem{HelRom00}
H\'elein, F., Romon, P., \emph{{H}amiltonian stationary Lagrangian surfaces in Hermitian symmetric
spaces}, preprint, math.dg/0010231

\bibitem{HerPed97}
Hertrich-Jeromin, U., Pedit, F., \emph{{R}emarks on the Darboux transform
of isothermic surfaces}, Doc. Math. \textbf{2} (1997), 313--333

\bibitem{KupWil81}
Kupershmidt, B.A., Wilson, G., \emph{{M}odifying Lax equations and the second
Hamiltonian structure}, Invent. Math. \textbf{62} (1981), 403--436

\bibitem{MaMa01}
Ma, H., and Ma, Y., \emph{{T}otally real minimal tori in $C P^2$}, preprint, math.dg/0106141

\bibitem{McI94}
McIntosh, I., \emph{{G}lobal solutions of the elliptic 2D periodic Toda
lattice}, Nonlinearity, \textbf{7} (1994), 85--108

\bibitem{McI02}
McIntosh, I., \emph{{S}pecial Lagrangian cones in $\C^ 3$ and primitive harmonic
maps}, preprint, math.dg/0201157

\bibitem{NomSas94}
Nomizu, K., Sasaki, T., \emph{{A}ffine differential geometry. Geometry of affine
immersions}, Cambridge Tracts in Mathematics, \textbf{111} (1994), Cambridge University Press

\bibitem{PinSte89}
Pinkall, U., Sterling, I., \emph{{O}n the classification of constant
mean curvature tori}, Ann. of Math. \textbf{130} (1989), 407--451 

\bibitem{PreSeg86}
Pressley, A., Segal, G. B., \emph{{L}oop Groups}, Oxford Science Publ., 
Clarendon Press, Oxford, (1986)

\bibitem{Sat84}
Sattinger, D.H., \emph{{H}amiltonian hierarchies on semi-simple Lie
algebras}, Stud. Appl. Math., \textbf{72} (1984), 65--86

\bibitem{SchWol99}
Schoen, R., Wolfson, J., \emph{{M}inimizing volume
among Lagrangian submanifolds}, Proc. Sympos. Pure Math.
AMS, \textbf{65}, (1999)

\bibitem{Sha91}
Sharipov, R., \emph{{M}inimal tori in the five dimensional sphere in $\C^3$}, Theoret. Math.
PHys., \textbf{87} (1991), 363--369

\bibitem{Ten85}
Tenenblat, K., \emph{{B}\"acklund's theorem for submanifolds of space
forms and a generalized wave equation}, Boll. Soc. Brasil.
Mat., \textbf{16} (1985), 67--92

\bibitem{Ter97}
Terng, C.L., \emph{{S}oliton equations and differential
geometry}, J. Differential Geometry, \textbf{45} (1997), 407--445

\bibitem{TerUhl98}
Terng, C.L., Uhlenbeck, K., \emph{{P}oisson actions and
scattering theory for integrable systems}, Surveys in Differential
Geometry: Integrable systems (A supplement to J. Differential
Geometry), \textbf{4} (1998), 315--402

\bibitem{TerUhl00a}
Terng, C.L., Uhlenbeck, K., \emph{{B}\"acklund transformations
and loop group actions}, Comm. Pure. Appl. Math., \textbf{53} (2000), 1--75

\bibitem{TerUhl00b}
Terng, C.L., Uhlenbeck, K., \emph{{G}eometry of
solitons}, Notice, A.M.S., \textbf{47} (2000), 17--25

\bibitem{Uhl89}
Uhlenbeck, K., \emph{{H}armonic maps into Lie group (classical solutions of the
Chiral model)}, J. Differential Geometry, \textbf{30} (1989), 1--50

\bibitem{Wil91}
Wilson, G., \emph{{T}he $\tau$-functions of the  $\cg$AKNS
equations}, Integrable systems, the Verdier Memorial,  Progress in Math.,
\textbf{115} (1991), 147--162

\bibitem{ZakSha79}
Zakharov, V.E., Shabat, A.B., \emph{{I}ntegration of non-linear equations of
mathematical physics by the inverse scattering method, II}, Funct. Anal.
Appl., \textbf{13} (1979), 166--174


\end{thebibliography}
\end{document}